\documentclass{amsart}

\usepackage[top=1.65in, bottom=1.65in, right=1.3in, left=1.3in]{geometry}

\usepackage{amsfonts}
\usepackage{amssymb}
\usepackage{amsmath}
\usepackage{fancyhdr}
\usepackage{eufrak}
\usepackage{latexsym}
\usepackage{amsbsy}
\usepackage[all]{xy}
\usepackage[latin1]{inputenc}
\usepackage[T1]{fontenc}        
\usepackage{amsmath}

\newtheorem{defn}{Definition}[section]
\newtheorem{thm}[defn]{Theorem}
\newtheorem{prop}[defn]{Proposition}
\newtheorem{lem}[defn]{Lemma}
\newtheorem{cor}[defn]{Corollary}
\newtheorem{rem}[defn]{Remark}

\newcommand{\fin}{\hfill{\Large$\Box$}\\}
\newcommand{\finsec}{\hfill{\Large$\Box$}}

\newcommand{\Eg}{{\cal E}_{\sl Green}}

\newcommand{\eps}{\varepsilon}

\newcommand{\z}{\zeta}
\newcommand{\la}{\lambda}

\newcommand{\om}{\omega}

\newcommand{\Q}{\mathbb {Q}}

\newcommand{\C}{\mathbb {C}}
\newcommand{\D}{\mathbb {D}}

\newcommand{\R}{\mathbb {R}}
\newcommand{\N}{\mathbb {N}}

\newcommand{\Pj}{\mathbb {P}}

\newcommand{\cd}{\textsf {Card}}

\newcommand{\Id}{{\rm Id}}
\newcommand{\Lip}{{\rm Lip \, }}
\newcommand{\Jac}{{\rm Jac \,  }}

\newcommand{\supp}{{\rm supp \, }}

\def\norm#1{\left\|\, #1\,\right\|}

\def\FF{{\cal F}}

\def\JJ{{\cal J}}

\def\LL{{\cal L}}

\def\MM{{\cal M}}

\def\NN{{\cal N}}

\def\com{\ar@{}[rd]|{\circlearrowleft}}

\newcommand{\cal}{\mathcal}

\title {Dynamical stability and Lyapunov exponents for holomorphic endomorphisms of $\Pj^k$}

\begin{author}[F.~Berteloot]{Fran\c cois Berteloot}
\address{ 
 Universit\'e de Toulouse - IMT\\
 UMR CNRS 5219\\
 31062 Toulouse Cedex \\
  France }
  \email{francois.berteloot$@$math.univ-toulouse.fr}
\end{author}

\begin{author}[F.~Bianchi]{Fabrizio Bianchi}
\address{ 
 Universit\'e de Toulouse - IMT\\
 UMR CNRS 5219\\
 31062 Toulouse Cedex \\
  France}
\email{f.bianchi@imperial.ac.uk}
\end{author}

\begin{author}[C.~Dupont]{Christophe Dupont}
\address{ 
 Universit\'e de Rennes - IRMAR\\
 UMR CNRS 6625\\
 35042 Rennes Cedex \\
  France }
  \email{christophe.dupont$@$univ-rennes1.fr}
\end{author}

\thanks{This research was partially supported by the ANR project LAMBDA,
ANR-13-BS01-0002. The work of the second author was partially supported by
the FIRB2012 grant ``Differential Geometry and Geometric Function
Theory''.}

\begin{document}
\maketitle 
\begin{abstract} We introduce a notion of stability for equilibrium measures in holomorphic families of endomorphisms of $\Pj^k$
and prove that it is equivalent to the stability of repelling cycles and equivalent to the existence of some  measurable holomorphic motion of Julia sets which we call equilibrium lamination. We characterize the corresponding bifurcations by the strict subharmonicity of the sum of Lyapunov exponents or the instability  of critical dynamics and analyze how repelling cycles may bifurcate. Our methods deeply exploit the properties of
Lyapunov exponents and are based on ergodic  and pluripotential theory. 

\end{abstract}
\vskip .4 cm 
\small{\noindent \emph{Key Words}: holomorphic dynamics, dynamical stability, positive currents, Lyapunov exponents.

\vskip .4 cm

}
 
\normalsize

%%%%%%%%%%%%%%
\section{Introduction}
%%%%%%%%%%%%%%%

%%%%%%%%%%%%%%%
\subsection{Main  definitions and results }
%%%%%%%%%%%%%%%%

 In the early 1980's, Ma\~{n}\'e, Sad and Sullivan \cite{MSS} and Lyubich \cite{L1, L2}  independently obtained fundamental results on the stability of holomorphic families 
$\left(f_{\la}\right)_{\la\in M}$ of rational maps of the Riemann sphere $\Pj^{1}$. They proved that the parameter space $M$ splits into an open and dense
\emph{stability locus} and its complement, the \emph{bifurcation locus}. They also obtained precise information on the distribution of hyperbolic parameters which led
to the so-called hyperbolic conjecture. This conjecture asserts that hyperbolic maps are dense in the space of rational maps.
The works of Douady and Hubbard on the Mandelbrot set provide a deeper understanding of these questions for
the quadratic polynomial family.  

In this theory, the finiteness of the critical set and the Picard-Montel theorem play a crucial role.
They allow to characterize the stability
of a parameter $\la_0 \in M$ by the stability of the critical orbits of the map $f_{\la_0}$. Equivalently, $\la_0$ is in the bifurcation locus if, after an arbitrarily small perturbation, there exists a repelling cycle capturing a critical orbit. The one-dimensional setting also permits, by means of the $\la$-lemma, to build holomorphic motions of Julia sets which conjugate the dynamics on connected components of the stability locus. The bifurcation locus also coincides with the closure
of the parameters $\la\in M$ for which $f_\la$ admits an unpersistent neutral cycle. \\

This article deals with bifurcations within holomorphic families of endomorphisms of $\Pj^k$ for $k\ge 1$.  Let $M$ be a connected complex manifold of dimension $m$.  A  holomorphic family of endomorphisms of $\Pj^k$ can be seen as a holomorphic mapping
$$f : M\times \Pj^{k}\to M\times \Pj^{k} \ \ , \ \ (\lambda,z) \mapsto (\lambda, f_\lambda(z))  $$ 
where the algebraic degree $d$ of $f_\lambda$ is larger than or equal to $2$ and does not depend on $\la$. For instance, $M$ can be the space ${\cal H}_d(\Pj^k)$ of all degree $d$ holomorphic endomorphisms of $\Pj^k$, which is a Zariski open subset in some $\Pj^N$. 

Our main result is Theorem \ref{main2} below, it asserts that different natural notions of stability are equivalent and leads to  a coherent notion of bifurcation for holomorphic families $f$ in $\Pj^k$. 
Our arguments exploit some  ergodic and pluripotential tools as those 
developped in the works of  Bedford-Lyubich-Smillie,
Fornaess-Sibony, Briend-Duval, Dinh-Sibony on holomorphic dynamics on $\Pj^k$ or $\C^k$ (see the survey \cite{DS2} for precise references). 
Let us recall that, for each $\la\in M$, we have an ergodic dynamical system $(J_\la,f_\la,\mu_\la)$ where $\mu_\la$ is the equilibrium measure of $f_\la$
and $J_\la$ is the topological  support of $\mu_\la$ called the \emph{Julia set}.
The measure $\mu_\la$ enjoys a potential interpretation
$$\mu_\la= \left( dd^c_{z} \,  g(\la,z)  +\omega_{FS} \right)^k,$$ 
where $g$ is the Green function of $f$ and $\omega_{FS}$ the Fubini-Study form on $\Pj^k$.
The repelling cycles of $f_\la$ equidistribute the measure $\mu_\la$ and hence  are dense in $J_\lambda$. However, in higher dimension, some repelling cycles may belong 
to the complement of $J_\la$. We denote by $L(\lambda) := \int_{\Pj^k} \log \vert \Jac f \vert\,  d\mu_\lambda$ the sum of the Lyapunov exponents of $\mu_\lambda$. This is a  plurisubharmonic function on $M$ which satisfies $L(\lambda) \geq k {\log d \over 2}$.  Let $[C_f]$ denote the current of integration on the critical set $C_f$ of $f$ taking into account the multiplicities of $f$. \\

Our main result is as follows. The  definitions occuring in  (A), (C), (D) and (F) are explained below. 
 
\begin{thm}\label{main2}  
 Let $f:M\times \Pj^{k}\to M\times \Pj^{k}$
be a holomorphic family of endomorphisms where $M$ is a simply connected open subset of the space ${\cal H}_d(\Pj^k)$ of
 endomorphisms of $\Pj^k$ of degree $d \geq 2$. 
 Then the following assertions are equivalent:
\begin{enumerate}
\item[(A)] the repelling $J$-cycles move holomorphically over $M$,
\item[(B)]  the function $L$ is pluriharmonic on $M$,
\item[(C)]  $f$ admits an equilibrium  web,
\item[(D)] $f$ admits an equilibrium lamination,
\item[(E)] any $\lambda_0 \in M$ admits a neighbourhood $U$ such that
  $\liminf_n d^{-k n}  \vert { (f^n)_{*}  [C_f]}\vert_{U\times\Pj^k} = 0$,
\item[(F)] there are no Misiurewicz parameters in $M$.
\end{enumerate}
When $k=2$, these equivalences remain true  for every simply connected manifold $M$. 
If one of these conditions is satisfied, $f$ admits a unique equilibrium  web $\cal M$ and ${\cal M}\left({\cal L}_1 \Delta {\cal L}_2\right)=0$
for any pair of equilibrium laminations ${\cal L}_1,{\cal L}_2$ of $f$.
\end{thm}

Theorem \ref{main2} leads us  to define the \emph{bifurcation current}  of a holomorphic family of endomorphisms of $\Pj^k$ as the closed positive current $dd^c_\la  L$, and the \emph{bifurcation locus} as the support of this current.  The family is \emph{stable} if its bifurcation locus is empty, stability is clearly a local notion. This is coherent with the  one-dimensional definition of the bifurcation current, due to DeMarco \cite{dM}. We stress that Theorem \ref{main2} stays partially true for general families (see Theorem \ref{main}).\\

Let us now specify the definitions. A central notion  is the set
$${\cal J} := \left \{ \gamma : M \to \Pj^k \, \colon  \, \gamma \textrm{ is holomorphic and } \gamma(\la) \in J_\la \; \textrm{for every} \;  \la \in M \right \}  .$$
The graph $\{\left(\la,\gamma(\la)\right)\;\la\in M\}$ of any element $\gamma\in {\cal J}$ is denoted $\Gamma_\gamma$.
We endow $\cal J$  with the topology of local uniform convergence and note that $f$ induces a continuous self-map
$${\cal F} : {\cal J}\to {\cal J}\;\textrm{ given by}\; {\cal F}\cdot\gamma (\la):= f_\la(\gamma(\la)).$$

\begin{defn} For every $\lambda \in M$, a repelling $J$-cycle of $f_\lambda$ is a repelling cycle which belongs to $J_\lambda$. We say that these cycles \emph{move holomorphically over} $M$ if, for every period $n$, there exists a finite subset $\{\rho_{n,j},\; 1\le j\le N_n\}$ of ${\cal J}$ such that $\{\rho_{n,j} (\la),\;1\le j\le N_n\}$ is precisely the set of $n$ periodic repelling $J$-cycles of $f_\la$ for every $\la \in M$.
\end{defn}

The holomorphic motion of repelling $J$-cycles over $M$ also means that for every repelling periodic point $z_0 \in J_{\la_0}$ of $f_{\la_0}$ there exists $\gamma\in \JJ$ such that $\gamma(\la_0)=z_0$ and $\gamma(\la)$ is a periodic repelling point of $f_\la$ for every $\la\in M$.\\

Our notions of equilibrium webs and  laminations are as follows.

\begin{defn}\label{defiMotMea} An \emph{equilibrium web} is a probability measure $\cal M$ on $\JJ$ such that 
\begin{enumerate}
\item $\cal M$ is $\FF$-invariant and its support is a compact subset of $\cal J$, 
\item for every $\lambda \in M$ the probability measure ${\cal M}_\la :=\int_\JJ \delta_{\gamma(\la)}\;d\cal M(\gamma)$ is equal to $\mu_\lambda$.
\end{enumerate}
 \end{defn}
 
 This notion is related to Dinh's theory of woven currents and somehow means that the measures $\left(\mu_\la\right)_{\la\in M}$ are holomorphically glued together.  In this article we shall also say that the $\left(\mu_\la\right)_{\la\in M}$ move holomorphically when such a web exists.

 \begin{defn} 
 An \emph{equilibrium lamination} is a relatively compact subset $\LL$ of $\JJ$ such that 
 
 \begin{enumerate}
\item $\Gamma_\gamma \cap \Gamma_{\gamma'} = \emptyset$ for every distinct $\gamma,\gamma' \in \LL$,
\item $\mu_\lambda \{ \gamma(\lambda) , \gamma \in \LL\} =1 $ for every $\lambda \in M$,
\item $\Gamma_\gamma$ does not meet the grand orbit of the critical set of $f$ for every $\gamma \in \LL$,
\item the map $\FF : \LL \to \LL$ is $d^k$ to $1$.
\end{enumerate}
\end{defn}
 
The existence of an equilibrium lamination corresponds to the property  of structural stability for sets of full $\mu_\lambda$-measure.
 It is rather easy to show that the existence of an equilibrium lamination 
implies that of an equilibrium web. The converse is much more delicate,  
equilibrium laminations will be extracted from the support of equilibrium webs by using ergodic theory for
 the dynamical system $(\JJ,\FF,\MM)$.\\
 
 Misiurewicz parameters  play a central role
 for proving that  the vanishing of $dd^c_\la  L$ is a sufficient condition for stability. They are also very useful to study the structure of bifurcation loci. 

\begin{defn}
One says that $\lambda_0 \in M$ is a \emph{Misiurewicz parameter} if   there exists a holomorphic map $\gamma $ from a neighbourhood of $\lambda_0$ into $\Pj^{k}$ such that:
\begin{itemize}
\item[1)] $\gamma(\lambda)\in J_\lambda$ and is a repelling $p_0$-periodic point of $f_{\lambda}$ for some $p_0 \geq 1$ and every $\la$,
\item[2)] $(\lambda_0,\gamma(\lambda_0))\in f^{n_0}(C_f)$ for some $n_0 \geq 1$,
\item[3)] the graph $\Gamma_\gamma$ of $\gamma$ is not contained in $f^{n_0}(C_f)$.
\end{itemize}
\end{defn}

%%%%%%%%%%%%%%
\subsection{Sketch of the proofs and further results}
%%%%%%%%%%%%%%

The novelty of our approach relies  on two specific features. 
The first one is the use of a formula for the $dd^c$ of the sum $L$ of the Lyapunov exponents to read the interplay between bifurcations and critical dynamics. This formula is  
\begin{equation}\label{BaBe}
dd^c_\la  L=\pi_{M\star}\left( \left( dd^c_{\la,z} \,  g(\la,z)  +\omega_{FS} \right)^k\wedge [C_f] \right) ,
\end{equation}
it was proved by Bassanelli and the first author \cite{BB1}, see also Pham's formula in Theorem \ref{thmPham}.
Like in dimension one, our proofs crucially rely on the links between  bifurcations and instability  in the critical dynamics. In higher dimension these interactions cannot be detected by a simple application of Picard-Montel theorem and Formula (\ref{BaBe}) aims to replace this theorem.
 
The second feature is the introduction of  equilibrium webs to overcome the  lack of $\la$-lemma and build holomorphic motions of Julia sets. This is a  weaker, but natural, notion  dealing with  the measures $\mu_\la$ rather than with their supports $J_\la$.  
It should be stressed that equilibrium webs are actually obtained as limits of discrete measures 
by mean of a compactness statement which may be considered as a measurable version of the $\la$-lemma (see Lemma \ref{limMS}). Finally we stress that our arguments will also rely on Misiuriwicz parameters, some kind of transversality and perturbations of Siegel discs. \\

We now specify our approach and summarize the proof of Theorem \ref{main2}. Simultaneously we 
state some related results. The implication (A)$\Rightarrow$(B) is proved in Proposition \ref{thmHMR}. We actually establish a stronger statement:  we show that $dd^c_\la  L$  vanishes if $f$ admits an equilibrium  web which is a limit of discrete measures supported on graphs avoiding the critical set of $f$. We obtain (B)$\Leftrightarrow$(E) by using Formula (\ref{BaBe}) and the $f$-invariance of the Green function $g$ (see Proposition \ref{growth}).\\ 

 % This is done in subsection \ref{ddcl}. 

To show that the vanishing of $dd^c_\la  L$ is a sufficient condition for stability, we exploit the dynamics of the critical set
and, more specifically,  the notion of Misiurewicz parameters.
We first prove that the pluriharmonicity of $L$ prevents the apparition of such  parameters. To do this, we use again Formula (\ref{BaBe}) and a dynamical rescaling argument. This is done in subsection \ref{SecMis}.
To prove that the absence of Misiurewicz parameters implies the existence of an equilibrium  web, we apply our measurable version of the $\lambda$-lemma to sequences of discrete measures on pull-backs by $f^n$ of a
graph of repelling $J$-cycles avoiding the post-critical set of $f$, see Proposition \ref{propgraph}. The proof of the existence of such a graph is rather involved and requires entropy arguments. These results, which are valid in arbitrary families, are summarized in the following theorem.

\begin{thm}\label{main}  
Let $f:M\times \Pj^{k}\to M\times \Pj^{k}$
be a holomorphic family of endomorphisms of $\Pj^k$ of degree $d \geq 2$. Then the following assertions are equivalent:
\begin{enumerate}
\item[(a)]  the function $L$ is pluriharmonic on $M$,
\item[(b)]  there are no Misiurewicz parameters in $M$,
\item[(c)]  the restriction $f\vert_{B\times \Pj^k}$, where $B$ is any sufficiently small ball,  admits an equilibrium  web ${\cal M}=\lim_n{\cal M}_n$ and the graph of  any $\gamma \in \cup_n\supp{\cal M}_n$ avoids the critical set of $f$.
\end{enumerate}
\end{thm} 

Among equilibrium  webs, those giving no mass to the subset of 
$\gamma$'s in $\cal J$ whose graphs meet the grand orbit of the critical set of $f$ will play an essential role in the construction  of equilibrium laminations. Such  webs are called \emph{acritical} (see Definition \ref{defiASW}).
Both   Theorem \ref{main} and the implication (A)$\Rightarrow$(B) in Theorem \ref{main2} are used to get the following important fact.

\begin{cor}\label{CorASW}
Let $f:M\times \Pj^{k}\to M\times \Pj^{k}$
be a holomorphic family of endomorphisms of $\Pj^k$ of degree $d \geq 2$.
If the  repelling $J$-cycles move holomorphically over $M$ then $f$ admits an ergodic and acritical equilibrium  web.
\end{cor}

In section \ref{secFrom} we prove that (A)$\Rightarrow$(D).
We use  there the Corollary  \ref{CorASW} and
 exploit the  properties of the dynamical system $\left({\cal J},{\cal F},{\cal M}\right)$ where $\cal M$ is an acritical and ergodic equilibrium  web.
We show that the iterated inverse branches of $f$  are exponentialy contracting near  the graph $\Gamma_\gamma$ of $\cal M$-almost every $\gamma\in {\cal J}$ (see Proposition \ref{LemSF}).  This implies that for $\cal M$-almost every  $\gamma\in {\cal J}$
the graph $\Gamma_\gamma$  does not intersect any other graph $\Gamma_{\gamma'}$ where $\gamma\ne \gamma'\in \supp {\cal M}$ and
 allows us to
build equilibrium laminations (see Theorem \ref{PropExHM}).\\

 So far we have established that  (A)$\Rightarrow$(B), (B)$\Leftrightarrow$(E), (B)$\Leftrightarrow$(F),
 (A)$\Rightarrow$(D) and that (B)$\Rightarrow$(C') where (C') is a local version of (C) (see Theorem \ref{main}). We prove simultaneously that (C')$\Rightarrow$(C)$\Rightarrow$(A).
To this purpose, we investigate how the apparition of Siegel
discs may affect the continuity of $\la\mapsto J_\la$ in the Hausdorff topology (see Proposition \ref{SiegJul}).
The section \ref{ContJ}  is mainly devoted to this study,  it is the only place where a specific assumption on the parameter space $M$ is used 
(see Proposition \ref{thbouc}). Finally, one easily gets (D)$\Rightarrow$(C) by using our basic tool for constructing equilibrium webs (see Proposition \ref{propgraph}). This completes the proof of Theorem \ref{main2}.\\

In the last section, we investigate a few properties of bifurcation loci.
We first show that bifurcation loci contain some remarkable elements. Theorem \ref{main} says that Misiurewicz parameters are dense in any bifurcation locus.
In the same vein, we prove in Theorem \ref{theoAdS} that  the bifurcation locus in ${\cal H}_d(\Pj^k)$ coincides with the closure of the set of endomorphisms which admit  repelling $J$-cycles which bifurcate either by giving Siegel periodic cycles or repelling cycles outside the Julia set. 
We also show in Theorem \ref{thLa}  that in any stable family, all elements are Latt\`es maps as soon as one element is a Latt\`es map. This follows from the characterization of such maps by their Lyapunov exponents, see \cite{BL, BDu, Du2}.
We finally consider the possibility for a bifurcation locus to have  a non-empty interior and establish the following result.
 
\begin{thm} \label{thmopbif}
Let $f:M\times \Pj^{k}\to M\times \Pj^{k}$ be a holomorphic family  of endomorphisms of $\Pj^k$.  The set of parameters $\la$ for which $\Pj^k$ coincides with the closure of the post-critical set of $f_\la$ is  dense in any open
subset of the bifurcation locus  of $f$.
\end{thm}

Let us finally mention that bifurcation phenomena in families of polynomial automorphisms are studied by a totally different approach, the sharpest achievements are due to Dujardin and Lyubich in their recent work on the dissipative H\'enon maps of $\C^2$ \cite{DL}.\\

\emph{Acknowledgements}: We thank the referee for his careful reading as well as for his questions which enabled us to improve the exposition.

%%%%%%%%%%%%%%%%%%%%%
\section{Equilibrium  webs} \label{HolMot}
%%%%%%%%%%%%%%%%%%%%%
%%%%%%%%%%%%%%%%%%%%%
\subsection{Sufficient conditions for the existence of an equilibrium web} \label{HolMot}
%%%%%%%%%%%%%%%%%%%%%

\smallskip 
%Let $f:M\times \Pj^{k}\to M\times \Pj^{k}$
%be a holomorphic family of endomorphisms of $\Pj^k$ of degree $d \geq 2$. We recall that $M$ is a connected complex manifold %of dimension $m$ and that $f(\la,z)=\left(\la,f_\la(z)\right)$.
%Let $\mu_\la$ denote the equilibrium measure of $f_\la$ and let $J_\la$ denote the support of $\mu_\lambda$, this is the Julia set %of
%$f_\la$. We want here to define a notion of holomorphic motion for the family $\left(\mu_\la\right)_{\la\in M}$. 
Let us consider the set ${\cal O} \left(M,\Pj^k\right)$ of  holomorphic maps from $M$ to $\Pj^k$, endowed with the metric space topology of local uniform convergence, and the closed subspace
$${\cal J} := \left \{ \gamma \in {\cal O} \left(M,\Pj^k\right) \; \colon \; \gamma(\la) \in J_\la \; \textrm{for every} \;  \la \in M \right \} . $$

For any probability measure $\cal M$ on ${\cal O} \left(M,\Pj^k\right)$ and every $\la \in M$ we define the measure

$$ {\cal M}_\la:=\int \delta_{\gamma(\la)}\;d\cal M(\gamma).$$
This  is a probability measure on $\Pj^k$ which is actually equal to $p_{\la\star} \cal M$, where the mapping $p_\la: {\cal O} \left(M,\Pj^k\right) \to \Pj^k$ is given by 
$p_\la(\gamma):=\gamma(\la)$.\\
 
%Let us recall that an equilibrium  web for $f$ 
%is  a $\cal F$-invariant and compactly supported probability measure $\cal M$ on $\cal J$ such that
%${\cal M}_\lambda = \mu_\la$ for every $\la \in M$. 
We shall sometimes say that the measures $\mu_\la$
\emph{move holomorphically over} $M$ when $f$ admits an equilibrium strucural web.

Equilibrium  webs will  be obtained as limits of discrete measures on ${\cal O} \left(M,\Pj^k\right)$.
To this purpose we shall use the following simple tool which somehow plays the role of the classical $\la$-lemma.
%We refer to Lemma \ref{LemMSS} for a more general statement.

\begin{lem}\label{limMS}  Let $f:M\times \Pj^{k}\to M\times \Pj^{k}$
be a holomorphic family of endomorphisms of $\Pj^k$. 
Let $({\cal M}_n)_{n \geq 1}$ be a sequence of Borel probability measures on ${\cal O} \left(M,\Pj^k\right)$ such that:
\begin{itemize}
\item [1)] $\lim_n ({\cal M}_n)_\la = \mu_\la$ for every $\la \in M$,
\item [2)] ${\cal F}_\star {\cal M}_{n+1} ={\cal M}_n$ or ${\cal F}_\star {\cal M}_{n} ={\cal M}_n$ for every $n \geq 1$,
\item [3)] there exists a compact subset ${\cal K} \subset {\cal O} \left( M, \Pj^k \right)$ such that $\supp {\cal M}_n\subset{\cal K}$.
\end{itemize}
Then  any limit of $(\frac{1}{n}\sum_{l=1}^n{\cal M}_l)_n$ is an equilibrium web. When ${\cal F}_\star {\cal M}_{n} ={\cal M}_n$, any limit of
$({\cal M}_n)_n$ is an equilibrium web.
\end{lem}

\proof Let ${\cal N}_n:=\frac{1}{n}\sum_{l=1}^n {\cal M}_l$. By Assertion 3) $\left({\cal N}_n\right)_{n\ge 1}$ is a sequence of Radon probability measures on the compact metric space $\cal K$.  Banach-Alaoglu and Riesz-Markov theorems ensure that there exists a subsequence $\left({\cal N}_{n_k}\right)_{k\ge 1}$  converging weakly to a Radon probability measure $\cal M$ on $\cal K$. By Assertion 2), we have ${\cal F}_{\star}  {\cal N}_{n_k} = {\cal N}_{n_k} + {\cal E}_k $ where the mass of ${\cal E}_k$ is   less than $2/n_k$. This implies that ${\cal F}_{\star}{\cal M} = \cal M$ as measures on $\cal K$. Let us extend $\cal M$ to a Borel probability measure $\widetilde{\cal M}$ on ${\cal O}(M,\Pj^k)$ by setting $\widetilde{\cal M}(A):={\cal M}(A\cap {\cal K})$. Let us verify that $\widetilde{\cal M}$ is an equilibrium web. We  still have ${\cal F}_{\star}\widetilde{\cal M} = \widetilde{\cal M}$. Indeed,
\begin{eqnarray*}
{\cal F}_{\star}\widetilde{\cal M}\left(A\right)= {\cal M}\left( {\cal F}^{-1}(A)\cap {\cal K}\right) \ge {\cal M}\left( {\cal F}^{-1}(A) \cap {\cal K} \cap {\cal F}^{-1}(\cal K) \right) =\\
 {\cal M}\left( {\cal F}^{-1}(A\cap \cal K) \cap {\cal K} \right)=
 {\cal M}\left( {\cal F}^{-1}(A\cap \cal K) \right)
=\cal M (A \cap \cal K) = \widetilde {\cal M} (A)
\end{eqnarray*}
and the identity follows since ${\cal F}_{\star}\widetilde{\cal M}$ and $\widetilde{\cal M}$ are probability measures. From  $p_{\la\star} \widetilde{\cal M} 
=p_{\la\star} {\cal M}$ and $p_{\la\star} {\cal M} = \lim_k p_{\la\star} {\cal N}_{n_k}=\mu_\la$ provided by Assertion 1), we deduce $p_{\la\star} \widetilde{\cal M} =\mu_\la$. It remains to check $\supp \widetilde{\cal M}\subset {\cal J}$. If $\gamma_0\notin {\cal J}$ then $\gamma_0(\la_0)\notin \supp \mu_{\la_0}$ for some $\la_0\in M$. Let $V_0$ be a neighbourhood of $\gamma_0$ in ${\cal O}(M,\Pj^k)$ such that $p_{\la_0}(V_0)\subset \Pj^k\setminus \supp \mu_{\la_0}$. 
Then $$\widetilde{\cal M}(V_0) \le \widetilde{\cal M}\left( p_{\la_0}^{-1} (p_{\la_0}(V_0))\right)=p_{\la_0 \star}\widetilde{\cal M}(p_{\la_0}(V_0))=\mu_{\la_0}(p_{\la_0}(V_0))=0$$
implies that $\gamma_0 \notin \supp \widetilde{\cal M}$.\fin

We now explain how Lemma \ref{limMS} is concretely used to produce equilibrium  webs. The proof relies on the equidistribution of preimages of points, see the articles \cite{FS1, BD2, DS} and on the equidistribution of repelling cycles, see \cite{BD1}.
 
\begin{prop}\label{propgraph}
Let $f:M\times \Pj^{k}\to M\times \Pj^{k}$ be a  holomorphic family of  endomorphisms of $\Pj^k$ of degree $d$.
\begin{itemize}
\item[1)] Assume that $M$ is simply connected and that there exists  $\gamma \in {\cal O} \left(M,\Pj^k\right)$  such that the graph
$\Gamma_\gamma$  does not intersect the post-critical set of $f$. 
Then an equilibrium  web is given by 
any limit of $\left( \frac{1}{n} \sum_{i=1}^n \frac{1}{d^{ki}} \sum_{{\cal F}^i\cdot\sigma=\gamma} \delta_{\sigma}\right)_n$.
\item[2)] Assume that the repelling $J$-cycles of $f$ move holomorphically over $M$. Let $\left({\rho_{n,j}}\right)_{1\le j\le N_n}$ be the elements of $\cal J$ given by the motions of these $n$-periodic  cycles.
Then  an equilibrium  web is given by any limit of
$\left( \frac{1}{d^{kn}} \sum_{j=1}^{N_n} \delta_{\rho_{n,j}}\right)_n$.
\end{itemize}
In both cases, $f$ admits an equilibrium web $\cal M=\lim_n{\cal M}_n$ such that $\Gamma_\gamma\cap C_f =\emptyset$ for every $\gamma\in \cup_n
\supp {\cal M}_n$.
\end{prop}

\proof  1) The map $f^n:\left(M\times \Pj^k\right)\setminus f^{-n}\left(\cup_{1\le p\le n} f^p(C_f)\right)\to \left(M\times \Pj^{k}\right)\setminus \left(\cup_{1\le p\le n} f^p(C_f)\right)$
is a  covering of degree $d^{kn}$. Hence, there exist $d^{kn}$ holomorphic graphs $\Gamma_{\sigma_{j,n}}$ such that $f^n \left(\Gamma_{\sigma_{j,n}}\right)=\Gamma_\gamma$
i.e. ${\cal F}^n\cdot\sigma_{j,n}=\gamma$. Let us set ${\cal M}_n:=\frac{1}{d^{kn}}\sum_{j=1}^{d^{kn}} \delta_{\sigma_{j,n}}$. By construction ${\cal F}_\star{\cal M}_{n+1}={\cal M}_n$ and, for every $\la \in M$, one has $({\cal M}_n)_\la=\frac{1}{d^{kn}}\sum_{j=1}^{d^{kn}} \delta_{\sigma_{j,n}(\la)}= \sum_{f^n_\la (x)=\gamma(\la)} \delta_{x} \to \mu_\la$, where the limit comes from the fact that $\gamma(\la) \notin \cup_{p\ge 1} f_\lambda^p (C_{f_\lambda})$. The family $\left(\sigma_{j,n}\right)_{j,n}$ is normal, by a theorem of Ueda \cite[Theorem 2.1]{U}, and therefore
the closure $\cal K$ of $\cup_{n\ge 1} \supp {\cal M}_n$ is a compact subset of ${\cal O} \left(M,\Pj^k\right)$.  The conclusion immediately follows from Lemma \ref{limMS}.\\

2) Let us set ${\cal M}_n:=\frac{1}{d^{kn}}\sum_{j=1}^{N_n} \delta_{\rho_{j,n}}$. The convergence of $({\cal M}_n)_\la$ towards $\mu_\la$ follows from  the equidistribution of repelling periodic points with respect to the  equilibrium measure, see \cite{BD2} (note that the repelling cycles produced there are $J$-cycles). The normality of the family
$\left(\rho_{j,n}\right)_{j,n}$ can be seen by lifting these curves to curves of periodic points of a lift of $f$ to $\C^{k+1}$.
Again, one concludes by using Lemma \ref{limMS}. \finsec

%%%%%%%%%%%%%%%%
\subsection{Acritical webs}
%%%%%%%%%%%%%%%%
To construct equilibrium laminations, it will be crucial to deal with
 equilibrium  webs  giving no mass to the subset ${\cal J}_s$ of $ {\cal J}$
whose elements have a graph intersecting the grand orbit of the critical set of $f$. 

 \begin{defn}\label{defiASW} An equilibrium  web $\cal M$ is said \emph{acritical} if ${\cal M}\left({\cal J}_s\right)=0$ where
 ${\cal J}_s$ is the singular part of $\JJ$ given by ${\cal J}_s:=\{\gamma\in {\cal J}\;\colon\; \Gamma_\gamma\cap \left(\cup_{m\ge 0}f^{-m}\left(\cup_{n\ge 0} f^n \left(C_f\right)\right) \right) \ne \emptyset\}.$
\end{defn}
As it will turn out,  equilibrium  webs given by Proposition \ref{propgraph}  are acritical and
this property,  combined with  ergodicity,  will allow us to build equilibrium laminations. This motivates the following result.
 
\begin{prop}\label{PropGE} Let $f:M\times\Pj^k\to M\times\Pj^k$ be a holomorphic family of  endomorphisms of $\Pj^k$. If $f$ admits an acritical equilibrium  web ${\cal M}_0$ then
$f$ admits an acritical equilibrium  web ${\cal M}'_0$ which is ergodic and such that  $\supp {\cal M}'_0 \subset \supp {\cal M}_0$.
\end{prop}

\proof Let us consider the convex set ${\cal P}_{web}\left(\cal K\right)$ of equilibrium  webs of $f$ which are supported in $\cal K$,
where ${\cal K}:=\supp \left({\cal M}_0\right)$. Note that ${\cal F}(\cal K)\subset \cal K$ since ${\cal M}_0$ is $\cal F$-invariant. The set  ${\cal P}_{web}\left(\cal K\right)$ is a compact metric space for the topology of weak convergence of measures. It is actually closed in the unit ball $B_{C({\cal K})'}$ where 
$C(\cal K)$ is the separable Banach space of continuous functions on $\cal K$ endowed with the norm of uniform convergence.

%We will use Choquet decomposition theorem to find extremal points ${\cal M}'$ in ${\cal P}_{web}\left(\cal K\right)$ for %which ${\cal M}'\left({\cal J}_s\right)=0$
%and then prove the ergodicity of ${\cal M}'$ by showing  that these points are also extremal in the set ${\cal P}_{inv}\left(\cal %K\right)$ of $\cal F$-invariant probability measures on $\cal K$.
Let us denote by $\mbox{Ext}\left({\cal P}_{web}\left(\cal K\right)\right)$ the set of extremal points of the compact metric space ${\cal P}_{web}\left(\cal K\right)$.
By Choquet theorem, there exists a probability measure $\nu_0$ on $\mbox{Ext}\left({\cal P}_{web}\left(\cal K\right)\right)$ such that
\begin{eqnarray*}
{\cal M}_0=\int_{\mbox{Ext}\left({\cal P}_{web}\left(\cal K\right)\right)} {\cal E} \;d\nu_0\left(\cal E\right). 
\end{eqnarray*}
Then $$0={\cal M}_0\left({\cal J}_s\right)=\int_{\mbox{Ext}\left({\cal P}_{web}\left(\cal K\right)\right)} {\cal E}\left({\cal J}_s\right)\;d\nu_0\left(\cal E\right)$$ and
 the set of equilibrium  webs ${\cal E}\in \mbox{Ext}\left({\cal P}_{web}\left(\cal K\right)\right)$ for which
${\cal E}\left({\cal J}_s\right) =0$ has full $\nu_0$-measure.

To conclude, we are left to check that any ${\cal M}'\in \mbox{Ext}\left({\cal P}_{web}\left(\cal K\right)\right)$ is also extremal in ${\cal P}_{inv}\left(\cal K\right)$ and therefore ergodic. Assume that
${\cal M}'=\frac{1}{2} {\cal M}_1+\frac{1}{2} {\cal M}_2$ where ${\cal M}_j\in {\cal P}_{inv}\left(\cal K\right)$. Then, as ${\cal M}'$ is an equilibrium  web  we have
$\mu_\la=p_{\la\star} \left({\cal M}'\right)= \frac{1}{2} p_{\la\star}\left({\cal M}_1\right)+\frac{1}{2} p_{\la\star}\left({\cal M}_2\right)$ for every $\la \in M$.
Since $p_\la\circ {\cal F}=f_\la\circ p_\la$, the probability measures  $p_{\la\star}\left({\cal M}_j\right)$ are $f_\la$-invariant and therefore the ergodicity of $\mu_\la$ implies that $p_{\la\star}\left({\cal M}_1\right) =p_{\la\star}\left({\cal M}_2\right)=\mu_\la$. This shows that ${\cal M}_1$ and ${\cal M}_2$ actually belong to 
${\cal P}_{web}\left(\cal K\right)$ and the identity ${\cal M}'={\cal M}_1={\cal M}_2$ then follows from the fact that ${\cal M}'$ is extremal in 
${\cal P}_{web}\left(\cal K\right)$. \fin

The following simple dynamical properties of  the support of an equilibrium  web will be very useful. We thank R. Dujardin for pointing us this fact. 

\begin{lem}\label{lemMCR}
Let $M$ be a connected complex manifold and 
$f:M\times \Pj^{k}\to M\times \Pj^{k}$ be a holomorphic family  of endomorphisms of $\Pj^k$ which admits an equilibrium   web $\cal M$.
 Then:
\begin{enumerate}
\item[1)] the sequence $\left(f^p_\la(\gamma(\la))\right)_{p\ge 1}$ is normal for every $\gamma\in \supp {\cal M}$,
\item[2)] for every $(\la_0,z_0)\in M\times J_{\la_0}$ there exists
$\gamma\in \supp {\cal M}$ such that $z_0=\gamma(\la_0)$,
\item[3)] for every $(\la_0,z_0)\in M\times J_{\la_0}$ such that $z_0$ is $n$-periodic and repelling for 
$f_{\la_0}$, there exists a unique 
$\gamma \in \supp {\cal M}$ such that $z_0=\gamma(\la_0)$ and $\gamma(\la)$ is $n$-periodic for 
$f_{\la}$ for every $\la\in M$. 
\end{enumerate}
\end{lem}

\proof (1) We use $f^p_\la(\gamma(\la)) =\left({\cal F}^p\cdot \gamma\right) (\la)$
and the fact that $\cal M$ is  compactly supported and $\cal F$-invariant.

\smallskip

(2)  As $z_0\in J_{\la_0}$ and $J_{\la_0}=\supp \mu_{\la_0}=\supp {\cal M}_{\la_0}$,  there exist  
$(\gamma_n)_n\subset \supp {\cal M}$ such that $\gamma_n(\la_0)\to z_0$. Then, since $\cal M$ is compactly supported, we can take for $\gamma$ any limit of $\left(\gamma_n\right)_n$.
 
\smallskip

(3)  By the implicit function theorem, there exists a neighbourhood $V_{\la_0}$ of $\la_0$ and a holomorphic
map $w : V_{\la_0}\to \Pj^k$ such that $w(\la_0)=z_0$ and $w(\la)$ is $n$-periodic  for $f_\la$. We will show that  $w$ coincides on $V_{\la_0}$ with the map $\gamma$ given by the previous item; the conclusion then follows by analytic continuation. Our argument is local, so we can choose a chart and work on $\C^k$. Since $z_0$ is repelling, we can shrink $V_{\la_0}$ and find $A>1$, $r>0$ such that
\begin{equation}\label{eqnRep}
\Vert w(\la) - f^n_\la(z)\Vert = \Vert f^n_\la (w(\la)) - f^n_\la(z)\Vert    \ge A  \Vert w(\la) - z\Vert
\end{equation}    
when  $\la\in V_{\la_0}$ and $\Vert w(\la) - z \Vert <r$. On the other hand the first item ensures that $\left(f^{pn}_\la(\gamma(\la))\right)_p$ is a normal family, hence we can shrink again $V_{\la_0}$ so that $\Vert w(\la) - f^{pn}_\la(\gamma(\la))\Vert <r$ for every $p \geq 1$ and $\la \in V_{\la_0}$. Combining this with Equation (\ref{eqnRep}) we obtain 
$r >\Vert w(\la) - f^{pn}_\la(\gamma(\la))\Vert \ge A^p \Vert w(\la) - \gamma(\la)\Vert$ for every $p\ge 1$ and $\la \in V_{\la_0}$. This  implies $w(\la)=\gamma(\la)$ on $V_{\la_0}$ since $A>1$. \fin

%%%%%%%%%%%%%%%%
\subsection{Webs and currents}
%%%%%%%%%%%%%%%%

For every probability measure $\cal M$ on ${\cal O} \left(M,\Pj^k\right)$ and in particular for any equilibrium  web  we can define the current
$ W_{\cal M}:=\int [\Gamma_{\gamma}]\;d{\cal M}(\gamma)$.
It has  bidimension $(m,m)$ on $M \times \Pj^k$ and is a woven current following Dinh's terminology \cite{Di2}.
To perform certain computations, we will have to explicitely relate  equilibrium  webs with positive horizontal currents
(see Lemma \ref{lemLift} below).
Before doing this, we recall some basic facts about horizontal currents.

\begin{defn}
Let $M$ be a complex connected manifold.
A current $\cal R$  on $M\times \C^{k+1}$ is horizontal if $\supp {\cal R} \subset M\times K$ for some compact subset $K \subset \C^{k+1}$.
 \end{defn}

Let us assume that $\cal R$ is a closed, positive, horizontal current of bidimension $(m,m)$ on $M\times \C^{k+1}$ where $m$ is the complex dimension of $M$.
Then the slices $\langle {\cal R}, \pi_{M}, \la\rangle$ exist for Lebesgue-almost every $\la \in M$ and are positive measures on $M\times \C^{k+1}$ supported 
on $\{\la\}\times \C^{k+1}$. The following \emph{ basic slicing formula}
 holds for every continuous test function $\psi$ on $M\times \C^{k+1}$ and every continuous $(m,m)$-test form $\omega$ on $M$:
\begin{equation}\label{BSF}
\int_{M} \langle {\cal R}, \pi_{M}, \la\rangle\; \psi\vert_{\{\la\}\times\C^{k-1}}\; \omega(\la) = \langle {\cal R} \wedge \pi_{M}^{\star}(\omega), \psi\rangle.
\end{equation}

Dinh and Sibony  have shown that the slices of such currents do actually exist for \emph{every} $\la\in M$, see  \cite[theorem 2.1]{DS1}. Their basic result is as follows.

\begin{thm}\label{DSthm}{\bf (Dinh-Sibony)} Let $M$ be a $m$-dimensional 
 complex connected manifold and ${\cal R}$ be a closed, positive, horizontal current of bidimension $(m,m)$ on $M\times \C^{k+1}$. Then
 the following properties occur:
\begin{enumerate}
\item the slice $\langle {\cal R}, \pi_{M}, \la\rangle$ exists for every $\la \in M$ and its mass does not depend on $\la \in M$,
\item the function $\la \mapsto  \int_{\C^{k+1}} \psi(\la,z) \, \langle {\cal R}, \pi_{M}, \la \rangle$ is \emph{psh} or  $\equiv -\infty$ on $M$
for any  \emph{psh} function  $\psi$ defined  on a neighborhood of $\supp {\cal R}$.
\end{enumerate}
\end{thm}

Let us now state the announced lemma. Let $\pi : \C^{k+1}\setminus \{0\} \to \Pj^k$ be the canonical projection.

\begin{lem}\label{lemLift}
Let $B$ be a ball in $\C^m$ and let $f : B\times\Pj^{k}\to B\times\Pj^{k}$ be a holomorphic family of endomorphisms of $\Pj^k$.
 Let $\cal K$ be a compact subset of ${\cal O} \left(B , \Pj^k\right)$. Then, after shrinking $B$, one may associate to any probability measure $\NN$ supported on $\cal K$ a positive, horizontal $(m,m)$-bidimensional current $\widetilde{W}_{\NN}$ on 
 $B\times \C^{k+1}$ such that $\pi_{\star} \langle \widetilde{W}_{\NN},\pi_B,\la\rangle = \NN_\la$ for every $\la \in B$.  Moreover,  $\widetilde{W}_\NN$ depends continuously on $\NN$.
\end{lem}

\proof
Let $\left(\sigma_i\right)_{1\le i\le N}$ be holomorphic sections of $\pi$ whose domains of definition $\Omega_i$ cover $\Pj^k$. 
Since $\cal K$ is a normal family, we may shrink $B$ so that for each $\gamma\in {\cal K}$ there exists at least one $1\le i\le N$ such that $\Gamma_\gamma\subset
B\times\Omega_i$. This allows to define a map
\begin{center}
$\sigma:{\cal K}\to {\cal O}\left(B,\C^{k+1}\right)$\\
$\gamma\mapsto \sigma(\gamma):=\sigma_l \circ \gamma$
\end{center} 
where $l:= \min \{1 \leq i \leq N \;\textrm {such that}\; \Gamma_{\gamma}\subset B\times\Omega_i\}$.
Now, for any probability measure $\NN$ supported on $\cal K$ we set
\begin{center}
$\widetilde{W}_{\NN}:=\int_{\cal J} [\Gamma_{\sigma(\gamma)}]\;d\NN(\gamma)$.
\end{center}
Then $\pi_{\star} \langle \widetilde{W}_{\NN},\pi_B,\la\rangle = \NN_\la$ for every $\la \in B$ by construction.\finsec

%%%%%%%%%%%
\subsection{Continuity of Julia sets and equilibrium  webs}\label{ssContJ}
%%%%%%%%%%%%

In Section \ref{ContJ}, we will want to compare the holomorphic motions of the measures $(\mu_\la)_{\la\in M}$ with the continuity of their supports $J_\la$ in the Hausdorff sense. To this purpose, we recall a few definitions. Let ${\textsc Comp}^{\star}\left(\Pj^k\right)$ be the set of non-empty compact subsets of $\Pj^k$ endowed with the Hausdorff distance and let $K_\epsilon$ denote the $\epsilon$-neighbourhood of  $K \in {\textsc Comp}^{\star}\left(\Pj^k\right)$. A map $E : M \to {\textsc Comp}^{\star}\left(\Pj^k\right)$
is said upper semi continuous ($u.s.c$) at $\la_0\in M$  if for every $\epsilon >0$, one has $E(\la)\subset \left(E(\la_0)\right)_\epsilon$ when $\la$ is close enough to $\la_0$. It is lower semi continuous ($l.s.c$) 
at $\la_0$ if for every $\epsilon >0$, one has $E(\la_0)\subset \left(E(\la)\right)_\epsilon$  when $\la$ is close enough to $\la_0$. For every $A \subset M \times \Pj^k$ we define $(A)_\la := A \cap (\{ \la \} \times \Pj^k)$.

The starting point about continuity of Julia sets relies on  the following observations, see also \cite[exercises 2.52 and 2.53]{DS2}.

\begin{prop}\label{PropBH}
Let $f:M\times \Pj^{k}\to M\times \Pj^{k}$ be a  holomorphic family of  endomorphisms of $\Pj^k$.
The map  $\la \mapsto J_{\la}$ from $M$ to ${\textsc Comp}^{\star}(\Pj^k)$  is $l.s.c$. If $f$ admits an equilibrium  web $\cal M$  and  $W_{\cal M}$ is the woven current $\int_{\cal J} [\Gamma_{\gamma}]\;d{\cal M}(\gamma)$, then $J_\la\subset \big(\supp W_{\cal M}\big)_\la$ and the map  $\la \mapsto \big(\supp W_{\cal M}\big)_\la$ from $M$ to ${\textsc Comp}^{\star}(\Pj^k)$ is $u.s.c$. 
\end{prop}

\proof 
The lower semi continuity of $J_\la$ follows from the existence of continuous local potentials for $\mu_\la$. 
Assume indeed that $\la\mapsto J_\la$ is not $l.s.c$ at $\la_0$. Then we may find
$\epsilon >0$ and sequences $\la_n \in M$, $z_n \in J_{\la_0}$ such that $d_{\Pj^k}(z_n,J_{\la_n})\ge \epsilon$. After taking a subsequence we may assume that
$z_n\to z_0\in J_{\la_0}$ and $B(z_0,\frac{\epsilon}{4})\subset B(z_n,\frac{\epsilon}{2}) \subset B(z_0,\epsilon)$. If $\epsilon$ is small enough, then
$\pi : \C^{k+1} \setminus \{ 0\} \to \Pj^k$ admits a section $\sigma$ on $B(z_0,2\epsilon)$ and the functions $u_\la(z):=G(\la,\sigma(z))$ are local potentials for the equilibrium measures, which means that the restriction of 
$\mu_\la$  to $B(z_0,2\epsilon)$ is the Monge-Amp\`ere mass $\left(dd^c_{z}  u_\la(z)\right)^k$
(see the beginning of Subsection \ref{someFF}). Observe that, by the continuity of $G$, 
the potentials $u_{\la_n}$ converge locally uniformly to $u_{\la_0}$. This implies that $\liminf_n \mu_{\la_n} \big(B(z_0,\frac{\epsilon}{4})\big) \ge 
\mu_{\la_0} \big(B(z_0,\frac{\epsilon}{8})\big)$. A contradiction follows:
$0 < \mu_{\la_0} \big(B(z_0,\frac{\epsilon}{8})\big) \le \liminf_n \mu_{\la_n} \big(B(z_0,\frac{\epsilon}{4})\big) \le \liminf_n \mu_{\la_n}   
\big(B(z_n,\frac{\epsilon}{2})\big) =0$.\\
The inclusion $J_\la\subset \big(\supp W_{\cal M}\big)_\la$ follows from  $J_\la =\supp \mu_\la$ and 
$\mu_\la ={\cal M}_\la=\int_{\cal J} \delta_{\gamma(\la)}\; d{\cal M}(\gamma)$. The upper semi continuity of 
$ \big(\supp W_{\cal M}\big)_\la$ is elementary  topology, see \cite[Proposition 2.1]{Do}.  \fin

It is now easy to see  that the existence of an equilibrium  web  implies that the Julia sets depend continuously on the parameter.

\begin{prop}\label{propHC}
Let $f:M\times \Pj^{k}\to M\times \Pj^{k}$ be a  holomorphic family of endomorphisms of $\Pj^k$. If $f$ admits an equilibrium  web then the map $\la\mapsto J_\la$ 
from $M$ to ${\textsc Comp}^{\star}\left(\Pj^k\right)$ is continuous.
\end{prop}

\proof  According to Proposition \ref{PropBH}, it suffices to show that $\left(\supp W_{\cal M}\right)_\la\subset J_{\la}$. This  follows from the following lemma. \finsec

\begin{lem}\label{lemHMH}
Let $f:M\times \Pj^{k}\to M\times \Pj^{k}$ be a  holomorphic family of endomorphisms of $\Pj^k$. Assume that  $f$ admits an equilibrium  web $\cal M$. If $z_0 \notin J_{\la_0}$ then there exist $\epsilon >0$ and $r_0>0$ such that 
${\cal M} \{\gamma\in {\cal J}\;\colon\;\Gamma_{\gamma} \cap\left[B(\la_0,\epsilon)\times B(z_0,r_0)\right] \ne \emptyset\}=0.$
Moreover $\mu_\la\left(B(z_0,r_0)\right) =0$ for every $\la\in B(\la_0,\epsilon)$.
\end{lem}

\proof 
Pick $r_0>0$ such that $\mu_{\la_0}\left(B(z_0,2r_0)\right) =0$. As $\supp {\cal M}$ is a normal family, there exists $\epsilon >0$ such that for any $\gamma \in\supp {\cal M}$:
$$  \Gamma_{\gamma} \cap\left[B(\la_0,\epsilon)\times B(z_0,r_0)\right] \ne \emptyset \Rightarrow \gamma(\la)\in B(z_0,2r_0)\;\textrm{for any}\; \la\in B(\la_0,\epsilon). $$
Let $\alpha:={\cal M} \{\gamma\in {\cal J}\;\colon\;\Gamma_{\gamma} \cap\left[B(\la_0,\epsilon)\times B(z_0,r_0)\right] \ne \emptyset\} $. Then, for any $\la\in B(\la_0,\epsilon)$, we have
$$ \alpha \le {\cal M} \{\gamma\in {\cal J}\;\colon\;\gamma(\la)\in B(z_0,2r_0)\} = \mu_\la \left( B(z_0,2r_0)\right).$$
Applying this to $\la_0$ yields $\alpha=0$ as desired. For every $\la\in B(\la_0,\epsilon)$ we have $\mu_{\la} \left(B(z_0,r_0)\right) ={\cal M} \{\gamma\in {\cal J}\;\colon\;\gamma(\la)\in B(z_0,r_0)\} \le \alpha =0$. This completes the proof. \finsec

%%%%%%%%%%%%%%%%%%%
\section{Stability and the sum of Lyapunov exponents}
%%%%%%%%%%%%%%%%%%%%%

In this section we establish (A)$\Rightarrow$(B)$\Leftrightarrow$(E) of  Theorem \ref{main2} and prove Theorem \ref{main} and  Corollary \ref{CorASW}. 
Formulas relating the critical dynamics with the sum of Lyapunov exponents are at the heart of our approach. For a polynomial $P$ of degree $d$, Przytycki  \cite{Pr} proved that the Lyapunov exponent of the equilibrium measure satisfies 
$$L(P)=\sum_{c\in C_P} G_P(c) + \log d $$
where $G_P (z) =\lim_n d^{-n} \log^+ \vert P^n(z)\vert$ is the dynamical Green function of $P$. This formula was generalized by DeMarco \cite{dM} for the Lyapunov exponent $L(f)$ of a rational map $f$. In several complex variables, Bedford-Jonsson \cite{BJ} established an analogous formula for the sum of the Lyapunov exponents of polynomials mappings. We use here an extended formula for holomorphic endomorphisms of $\Pj^k$ obtained by Bassanelli-Berteloot, see \cite[Theorem 4.1]{BB1}.

%%%%%%%%%%%%%%%%%%%%%%%
\subsection{Formulas for the sum of Lyapunov exponents} \label{someFF}
%%%%%%%%%%%%%%%%%%%%%%%

To deal with this kind of formulas, a suitable framework is that of equilibrium currents for holomorphic families of $d$-homogeneous non-degenerate maps.
It has been introduced by Pham \cite{Ph} in the more general context
of polynomial like mappings (see also the lecture notes by Dinh and Sibony \cite[section 2.5]{DS2}).
 
\begin{defn}  
Let $F:M\times \C^{k+1}\to M\times \C^{k+1}$ be a holomorphic family of $d$-homogeneous non-degenerate maps where $M$ is some $m$-dimensional complex connected manifold. Let $\cal E$ be a closed, positive, horizontal current of bidimension $(m,m)$ on $M\times \C^{k+1}$. We say that $\cal E$ is an \emph{equilibrium current } for $F$ if the slice 
$\langle {\cal E}, \pi_{M}, \la\rangle$ is equal to the equilibrium measure of $F_\la$ for every $\lambda \in M$.
\end{defn}

Contrary to equilibrium webs, equilibrium currents always exist. One may dynamically produce them and they do not detect bifurcations.  For instance, Pham proved that the sequence of smooth forms $\big(\frac{1}{d^{(k+1)n}} F^{n\star} \left(\pi_{\C^{k+1}}^{\star} \theta\right)\big)_n$ converges to such a current for any smooth probability 
measure $\theta$ on $\C^{k+1}$. Note that such currents are not unique when $k>1$.\\

It is also possible to define equilibrium currents for families of endomorphisms of $\Pj^k$ by means of Green functions.
Let us briefly recall their construction. We consider a holomorphic family  $f : M \times \Pj^k \to M \times \Pj^k$
which admits a lift $F : M \times \C^{k+1} \to M \times \C^{k+1}$. The sequence
$$G_n(\la,\tilde z):=\frac{1}{d^n}\log \Vert F_{\la}^{n}(\tilde z)\Vert$$
converges locally uniformly on $M \times \C^{k+1}\setminus\{0\}$ to a function $G$ which we call the \emph{Green function of $F$}. 
The norm $\Vert\;\Vert$ is the euclidean one. The function $G$ is \emph{psh} and H\" older continuous, see \cite[section 1.2]{BB1}.  Let $\pi : \C^{k+1} \setminus \{ 0 \} \to \Pj^k$ be the canonical projection and $\omega_{FS}$ be the Fubini-Study form of mass $1$ on $\Pj^k$. The functions $G_n$ induce  functions $g_n : M \times \Pj^k\to \R$ by setting $g_n(\lambda, z) :=   G_n(\lambda, \tilde z) - \log \norm {\tilde z}$, for every $\tilde z$ satisfying $\pi(\tilde z)=z$.  We have: 
$${1 \over d}  f^\star\left(dd_{\lambda,z}^c \, g_n+\omega_{FS}  \right)  = dd^c_{\lambda,z} \, g_{n+1}+\omega_{FS}. $$
We define similarly  $g(\lambda, z) := \lim_n g_n(\lambda,z)$, which is equal to $G(\lambda, \tilde z) - \log \norm {\tilde z}$, and set
$$\Eg:= \left( dd^c_{\la,z} \,  g  +\omega_{FS} \right)^k.$$
This is a current of bidimension $(m,m)$ and, since slicing commutes with the operators $d$, $d^c$, the measure $\langle \Eg,\pi_M, \la\rangle$ is equal to the equilibrium measure $\mu_\la$ of $f_\la$ for almost every $\lambda \in M$, using the horizontality of this current one sees that $\langle \Eg,\pi_M, \la\rangle$ actually equals $\mu_\la$ for all $\la\in M$ (see \cite{DS1}). The current ${\cal E}_{\sl Green}$ will play an important role in our study (see Proposition \ref{thmmisiu}). We call it the \emph{Green equilibrum current} of $f$.\\

Before stating the results of this subsection, we fix a few notations. Let us set $D:=(k+1)(d-1)$. The line bundle ${\cal O}_{\Pj^k} (D)$ over $\Pj^k$ is seen as the quotient of $\left(\C^{k+1}\setminus \{0\} \right) \times \C$ by the relation $(\tilde z,x)\equiv (u\tilde z,u^D x)$ for every $u\in \C^{\star}$ and its elements are denoted by $[\tilde z,x]$. We endow ${\cal O}_{\Pj^k} (D)$ with the canonical metric
\begin{eqnarray*}
\Vert [\tilde z,x]\Vert_0 := e^{-D\cdot\log \Vert \tilde z\Vert } \vert x\vert
\end{eqnarray*}
or, for any $\la \in M$, with the metric 
\begin{eqnarray*}
\Vert [\tilde z,x]\Vert_{\la}:= e^{-D\cdot G(\la,\tilde z) } \vert x\vert.
\end{eqnarray*}
Let us set $J_F (\lambda, \tilde z) :=  \det d_{\tilde z} F_\lambda $. Then we obtain a family of holomorphic sections of ${\cal O}_{\Pj^k} (D)$ by setting, for every $\tilde z \in \C^{k+1} \setminus \{ 0 \}$:
 \begin{eqnarray*}
J_F^s\left(\la,\pi(\tilde z)\right):=[\tilde z,J_F(\la,\tilde z)].
\end{eqnarray*}
Observe that 
\begin{eqnarray}\label{For13}
\Vert J_F^s\left(\la,\pi(\tilde z)\right) \Vert_\la = e^{-D\cdot G(\la,\tilde z) } \vert J_F(\la,\tilde z) \vert.
\end{eqnarray}
The current $[C_f]:=dd^c_{\la,z} \log\Vert J_F^s\left(\la, z\right)\Vert_0$ is the current of integration on $C_f$ taking account the topological multiplicities of $f$, its bidimension is equal to $(\kappa,\kappa)$ where $\kappa := k+m-1$.\\

\begin{thm}{\bf (Bassanelli-Berteloot)}\label{thmFBB}
Let $f:M\times \Pj^{k} \to M\times \Pj^{k} $ be a holomorphic family of endomorphisms of $\Pj^k$. Let $L(\la)$ be the sum of the Lyapunov exponents of $\mu_\la$. Then
$$dd^c_\lambda L = \pi_{M\star} \left( \Eg\wedge [C_f] \right) .$$
\end{thm}

We end this subsection by explaining how Pham \cite{Ph} obtained a more general formula. His result holds for any equilibrium current of any family of polynomial-like maps, we state it in the special case of  non-degenerate homogeneous maps for sake of simplicity.
Let us consider a holomorphic family  $F : M \times \C^{k+1} \to M \times \C^{k+1}$ which is the lift of $f : M \times \Pj^k \to M \times \Pj^k$.
Then the function $\log \vert J_F (\lambda, \tilde z) \vert$ is \emph{psh}  on $M \times \C^{k+1}$ and the sum of Lyapunov exponents of $F_\la$ with respect to its equilibrium measure $\nu_\la$ is given by $\int_{\C^{k+1}} \log \vert J_F (\lambda, \tilde z )\vert \  d\nu_\lambda(\tilde z)$ and is equal to $L (\lambda) + \log d$ where $L(\la)$ is the sum of Lyapunov exponents of ($f_\la,\mu_\la)$.

\begin{thm}{\bf (Pham)}\label{thmPham}
Let $F:M\times \C^{k+1}\to M\times \C^{k+1}$ be a holomorphic family of non-degenerate $d$-homogeneous  maps
and let  $\cal E$ be an equilibrium current for $F$. Then:
\begin{enumerate}
\item the current $\log \vert J_F\vert  \cdot  {\cal E}$ has locally finite mass,
\item $dd^c_\lambda L = \pi_{M\star} \left({\cal E}\wedge dd^c_{\la , \tilde z}  \log \vert J_F\vert \right)$.
\end{enumerate}
\end{thm}

To prove that $dd^c_\la L$ vanishes when  repelling $J$-cycles move holomorphically (subsection \ref{ddcl}),   
we shall actually need the following  formula for $dd^c_\lambda L$ whose  proof follows Pham's arguments.

\begin{prop}\label{specform} 
Let $B$ be an open ball in $\C^m$ and let $f:B\times \Pj^{k}\to B\times \Pj^{k}$
be a holomorphic family of endomorphisms of $\Pj^k$. Assume that $f$ admits  an equilibrium  web $\cal M$. Then
$$ dd^c_\lambda L = \pi_{B\star}\big( \widetilde{W}_{\cal M} \wedge dd^c_{\lambda,z} \log \Vert J_F^s \left(\la,\pi(\tilde z)\right)\Vert_{\la}\big) $$
where $\widetilde{W}_{\cal M}$ is the $(m,m)$-bidimensional current on $M \times \C^{k+1}$ associated to $\cal M$ by Lemma \ref{lemLift}.
\end{prop} 

\proof  We first check that for every $\lambda \in B$ we have
\begin{eqnarray}\label{fich}
\int_{\C^{k+1}} \log \Vert J_F^s \left(\la,\pi(\tilde z)\right)\Vert_{\la}\;\langle \widetilde{W}_{\cal M},\pi_{B},\la\rangle =L(\la)+\log d.
\end{eqnarray}
Indeed, since $\pi_{\star} \langle \widetilde{W}_{\cal M},\pi_{B},\la\rangle =\mu_\la$, we get
\begin{eqnarray*}
\int_{\C^{k+1}} \log \Vert J_F^s \left(\la,\pi(\tilde z)\right)\Vert_{\la}\;\langle \widetilde{W}_{\cal M},\pi_{B},\la\rangle = 
\int_{\Pj^{k}} \log \Vert J_F^s \left(\la,z\right)\Vert_{\la}\;\mu_\la.
\end{eqnarray*}
On the other hand, by Formula (\ref{For13}) and since $G_\la$ identically vanishes on the support of the equilibrium measure $\nu_\la$ of $F_\la$ and  $\pi_{\star} \nu_\la = \mu_\la$, we have
\begin{eqnarray*}
\int_{\Pj^{k}} \log \Vert J_F^s \left(\la,z\right)\Vert_{\la}\;\mu_\la = \int_{\C^{k+1}} \log \Vert J_F^s \left(\la,\pi(\tilde z)\right)\Vert_{\la}\; \nu_\la \\
= \int_{\C^{k+1}} \log \vert J_F(\la,\tilde z) \vert\; \nu_\la =L(\la) +\log d ,
\end{eqnarray*}
and the identity (\ref{fich}) follows.

Pham proved that $u \cdot {\cal R}$ has locally finite mass for every \emph{psh} function $u$ and every horizontal current $\cal R$ as soon as $\int_{\C^{k+1}} u(\la,\cdot) \, \langle {\cal R},\pi_{B}, \la\rangle \ne -\infty$ for some $\la\in M$, see \cite[theorem A.2]{Ph}. It thus follows from (\ref{fich}) that the current $\log \Vert J_F^s \left(\la,\pi(\tilde z)\right) \Vert_{\la} \cdot  \widetilde{W}_{\cal M}$ is well defined and that its $dd^c_{\la , \tilde z} $ is equal to
 $ \widetilde{W}_{\cal M} \wedge dd^c_{\la , \tilde z}  
 \log \Vert J_F^s \left(\la,\pi(\tilde z)\right)\Vert_{\la} $.
  
We conclude by simple computation which relies on integration by parts (to make it rigourous one should approximate $\log \Vert J_F^s \left(\la,\pi(\tilde z)\right)\Vert_{\la} $ by smooth functions). Let $\varphi$ be a $(m-1,m-1)$ test form on $B$. Then
\begin{eqnarray*}
\langle \pi_{B\star}\left( \widetilde{W}_{\cal M} \wedge dd^c_{\la , \tilde z}  \log \Vert J_F^s \left(\la,\pi(\tilde z)\right)\Vert_{\la} \right) , \varphi \rangle = \langle \, \log  \Vert J_F^s \left(\la,\pi(\tilde z)\right)\Vert_{\la} \cdot \widetilde{W}_{\cal M} \, , \,  dd^c_{\lambda, \tilde z}  \left(\pi_B^{\star}\varphi\right)\, \rangle\\
= \langle \, \widetilde{W}_{\cal M}\wedge \pi_{B}^{\star} \left( dd^c_{\lambda}  \varphi \right) \, , \, \log  \Vert J_F^s \left(\la,\pi(\tilde z)\right)\Vert_{\la} \, \rangle . 
\end{eqnarray*}
By the basic slicing formula (\ref{BSF}) and the identity (\ref{fich}), this is equal to 
\begin{eqnarray*}
 \int_{B} \left(\langle \widetilde{W}_{\cal M}, \pi_{B}, \la\rangle \log \Vert J_F^s \left(\la,\pi(\tilde z)\right)\Vert_{\la} \right) dd^c_{\lambda} \varphi=
\int_{B} L \,  dd^c_{\lambda} \varphi = \langle dd^c_{\lambda} L, \varphi\rangle.
\end{eqnarray*}
This completes the proof. \finsec

%%%%%%%%%%%%%%%%%%%%%%%%%%%%%%%%%%%%%%%%%%%%%%%%
\subsection{Repelling cycles do not move holomorphically on $\supp dd^c_\la L$}\label{ddcl}
%%%%%%%%%%%%%%%%%%%%%%%%%%%%%%%%%%%%%%%%%%%%%%%%

The following proposition will be used in subsection \ref{sspartof} to prove (A)$\Rightarrow$(B) of Theorem \ref{main2}, namely that $dd^c_\la L=0$ on $M$
if the repelling $J$-cycles move holomorphically.

\begin{prop}\label{thmHMR}
 Let  $f:M\times \Pj^{k}\to M\times \Pj^{k}$
be a holomorphic family of  endomorphisms of $\Pj^k$ which admits an equilibrium  web $\cal M$ which is given by ${\cal M}=\lim_n {\cal M}_n$ where $\Gamma_\gamma \cap C_f=\emptyset$ for any $\gamma\in \cup_n \supp {\cal M}_n$. Then $dd^c_\la L=0$ on $M$.
\end{prop}     

The fact that the holomorphic motion of the repelling $J$-cycles over $M$ imply the pluriharmonicity of $L$ on $M$ was proved in \cite[Theorem 2.2]{BB1} or \cite[Theorem 1.5]{BDM}. Proposition \ref{thmHMR} actually provides a more general result which will also be used for establishing the density of Misiurewicz parameters in bifurcation loci (see subsection \ref{SecMis}). The proof  requires  the following technical lemma.

\begin{lem}\label{lemEstiZ}
 Let $B$ be an open ball in $\C^m$ and let $f : B\times \Pj^{k}\to B\times \Pj^{k}$
be a holomorphic family of endomorphisms of $\Pj^k$. Let $Z$ be a codimension $1$ analytic subset of $B\times \Pj^k$ which does not contain any fiber $\{\la\}\times\Pj^k$. Assume that  there exists an equilibrium  web  satisfying
 ${\cal M}=\lim_n {\cal M}_n$, where $\Gamma_\gamma \cap Z=\emptyset$ for  every $\gamma\in \cup_n \supp {\cal M}_n$ and every $n\ge n_0$. 
 Then, after shrinking $B$, the following estimate occur. For any relatively compact ball $B'$ in $B$, there exist $A>0$ and $0<a<1$ such that  
$${\cal M}\left(\{\gamma\in{\cal J}\;\colon\; \Gamma_{\gamma\vert_{B'}} \cap Z_\epsilon \ne \emptyset\}\right)\le A\epsilon^a  $$
for every sufficiently small $\epsilon>0$, where $Z_\epsilon$ is the $\epsilon$-neighbourhood of $Z$. 
\end{lem}

\proof 
We can assume that both $B$ and $B'$ are centered at some $\la_0$.
After maybe shrinking $B$ we may find a finite collection $(\Omega_i,h_i)_{1\le i\le N}$ where the $\Omega_i$ are open and cover $B\times\Pj^k$, the functions
 $h_i$ are holomorphic and bounded by $1$ on $\Omega_i$ and $Z\cap\Omega_i=\{h_i=0\}$ for any $1\le i\le N$. If $\epsilon$ is small enough, we may also assume that $Z_\epsilon\cap \Omega_i \subset \{\vert h_i\vert <C_1\epsilon\}$
and, by \L ojasiewicz inequality, that $\{\vert h_i\vert <\epsilon\}\subset Z_{C_2\epsilon^\tau}$  for some constants $C_1, C_2, \tau >0$. Similarly, one has $Z_\epsilon \cap\left(\{\la_0\}\times \Pj^k\right) \subset \left(Z \cap\left(\{\la_0\}\times \Pj^k\right)\right)_{C_3\epsilon^{\tau_0}}$ for some constants $C_3, \tau_0 >0$. \\

Since $\cal M$ has compact support in $\cal J$, we may shrink $B$ again so that for any $\gamma\in \supp {\cal M}$ there exists at least one $1\le i\le N$ such that $\Gamma_\gamma\subset \Omega_i$. We shall use the following claim.\\ 

{\bf Claim}: {\it there exists $0< \alpha \le 1$ such that $\sup_{B'} \vert \phi \vert \le \vert \phi(t_0)\vert ^{\alpha}$ for every  $t_0\in B'$ and for every holomorphic function $\phi:B\to \D^*$}. \\

Let $\gamma\in \supp {\cal M}$ such that $\Gamma_{\gamma}\cap Z =\emptyset$ and  $\Gamma_{\gamma\vert_{B'}} \cap Z_\epsilon \ne \emptyset$.
Applying the Claim to $h_i\circ \gamma$ with $\Gamma_\gamma\subset \Omega_i$ we obtain that 
$\Gamma_{\gamma\vert_{B'}}\subset Z_{C_4\epsilon^{\tau\alpha}}$ for some constant $C_4>0$.\\
On the other hand, by our assumption on the approximation of $\cal M$ by ${\cal M}_n$, Hurwitz lemma implies that either $\Gamma_\gamma\subset Z$ or 
$\Gamma_{\gamma}\cap Z =\emptyset$ for any $\gamma\in \supp {\cal M}$. We thus have
\begin{eqnarray*}
{\cal M}\left(\{\gamma\in {\cal J}\;\colon\; \Gamma_{\gamma\vert_{B'}}\cap Z_{\epsilon} \ne\emptyset\}\right)\le
{\cal M}\left(\{\gamma\in {\cal J}\;\colon\; \Gamma_{\gamma\vert_{B'}}\subset Z_{C_4\epsilon^{\tau\alpha}}\}\right)\le\\
{\cal M}\left(\{\gamma\in {\cal J}\;\colon\; (\la_0,\gamma(\la_0))\in  Z_{C_4\epsilon^{\tau\alpha}}\}\right)=
\mu_{\la_0} \left( Z_{C_4\epsilon^{\tau\alpha}}\cap \left(\{\la_0\}\times\Pj^k\right)\right) \le\\
\mu_{\la_0} \left[\left( Z \cap \left(\{\la_0\}\times\Pj^k\right)\right)_{C_3 (C_4\epsilon^{\alpha\tau})^ {\tau_0}}\right] \le A\epsilon^a 
\end{eqnarray*}
where the last estimate is due to the fact that $\mu_{\la_0}$ has H\"{o}lder-continuous local potentials (see \cite{DS2} Proposition 1.18) and 
$Z\cap \left(\{\la_0\}\times\Pj^k\right)$ is a proper analytic subset of $\Pj^k$. \\

It remains to prove the Claim. Let $\mathcal G := \{ \varphi \in{\cal O}(B,H)\;\colon\; \varphi(s) = -1 \textrm{ for some } s \in \overline{B'}\}$ where $H := \{ \Re z < 0 \}$ is the left half plane. Then $\mathcal G$ is compact for the topology of local uniform convergence, and thus the quantity $(-\alpha) := \sup _{\varphi \in \mathcal G} \sup_{s \in \overline{B'}} \Re \varphi(s)$ satisfies $-1 \le -\alpha < 0$.  Let $t_0 \in B'$ and $\phi : B \to \mathbb D^*$ be holomorphic.  After a rotation in $\mathbb D^*$ we may assume that $\vert \phi(t_0) \vert = \phi(t_0) \in ]0,1[.$ Let $\varphi : B \to H$ be the lift of $\phi$ by the exponential map, which satisfies $\varphi(t_0) = \log \phi(t_0) \in ]-\infty, 0[$. Then $\varphi_0(t) := - \varphi(t) / \varphi(t_0)$ belongs to $\mathcal G$ and thus  $\Re(\varphi_0) \leq -\alpha$ on $B'$. This is the desired estimate since $\vert \phi\vert = e^{\Re \varphi }\le e^{\alpha \log \phi(t_0)}=\vert \phi(t_0)\vert^{\alpha}$.\fin

\textsc{Proof of Proposition \ref{thmHMR}:} The problem is local and we may therefore take for $M$ a ball $B\subset \C^m$ and assume that  $f : B\times \Pj^{k}\to B\times \Pj^{k}$ admits a lifted family $F : B\times \C^{k+1}\to B \times \C^{k+1}$ of $d$-homogeneous non-degenerate maps. We will apply Lemma \ref{lemEstiZ} with
$Z=C_f$. Let $B'$ be any  relatively compact ball contained in $B$. \\

After shrinking $B$ we may use Lemma \ref{lemLift} and associate to $\cal M$ the following horizontal current on $B\times \C^{k+1}$
$$\widetilde{W}_{\cal M}=\int_{\cal J} [\Gamma_{\sigma(\gamma)}]\;d{\cal M}(\gamma).$$
According to Proposition \ref{specform}, one has 
$$ dd^c_\la L = \pi_{B\star}\big( \widetilde{W}_{\cal M} \wedge dd^c_{\la , \tilde z}  \log \Vert J_F^s \left(\la,\pi(\tilde z)\right)\Vert_{\la}\big) .$$
Using $\Vert J_F^s\left(\la,\pi(\tilde z)\right) \Vert_\la = e^{-D\cdot G(\la,\tilde z) } \vert J_F(\la,\tilde z) \vert$ (see Formula (\ref{For13})),  and the fact that the functions $L$ and $G$ are \emph{psh}, we obtain
 \begin{eqnarray*}
 0 \le dd^c_\la L = \pi_{B\star}\big( \widetilde{W}_{\cal M} \wedge dd^c_{\la, \tilde z} \log \vert J_F\vert\big)  -D \pi_{B\star}\big( \widetilde{W}_{\cal M} \wedge dd^c_{\la, \tilde z} G\big)\\
 \le \pi_{B\star}\big( \widetilde{W}_{\cal M} \wedge dd^c_{\la, \tilde z} \log \vert J_F\vert\big).
 \end{eqnarray*}
Hence it suffices to show that the current $\log \vert J_F\vert \; \widetilde{W}_{\cal M}$ restricted to $B'\times \C^{k+1}$ is $dd^c_{\la , \tilde z}$ closed.\\

For $\epsilon <1$ we set $\log _\epsilon:= \chi_\epsilon\circ \log$ where $\chi_\epsilon$ is a convex, smooth, increasing function on $\R$ such that $\chi_\epsilon(x) =x$ if $x\ge \log \epsilon$ and $\chi_{\epsilon}(-\infty)=2\log \epsilon$. Then $\log_\epsilon \vert J_F\vert $ is a decreasing family (when $\epsilon\to 0$) of smooth \emph{psh} functions which converges to $\log \vert J_F\vert $. As $\lim_{\epsilon\to 0} \log_\epsilon \vert J_F\vert  \; \widetilde{W}_{\cal M} =
\log  \vert J_F\vert  \; \widetilde{W}_{\cal M}$ we will actually deal with $\log_\epsilon \vert J_F\vert  \; \widetilde{W}_{\cal M}$.\\
To this purpose we set $U_\epsilon:=\{\vert J_F\vert<\epsilon\}$,
$S_{{\cal M},\epsilon}:=\{\gamma\in \supp {\cal M}\;\colon\; \Gamma_{\sigma(\gamma)\vert_{B'}}\cap U_\epsilon \ne \emptyset\}$ and decompose 
$\widetilde{W}_{\cal M}$ as:
$$\widetilde{W}_{\cal M}=\widetilde{W}_{{\cal M},\epsilon}+\widetilde{W}_{{\cal M},\epsilon}^{\star}$$
where $\widetilde{W}_{{\cal M},\epsilon}:=\int_{\cal J} [\Gamma_{\sigma(\gamma)}]1_{S_{{\cal M},\epsilon}}\;d{\cal M}(\gamma)$ and 
$\widetilde{W}_{{\cal M},\epsilon}^{\star}:=\widetilde{W}_{\cal M}-\widetilde{W}_{{\cal M},\epsilon}$.
Then
\begin{eqnarray*}
\log_\epsilon \vert J_F\vert \; \widetilde{W}_{\cal M}= \log_\epsilon  \vert J_F \vert \; \widetilde{W}_{{\cal M},\epsilon}+  \log_\epsilon  \vert J_F\vert \; \widetilde{W}_{{\cal M},\epsilon}^{\star}
\end{eqnarray*}
and, by construction, the current $\log_\epsilon  \vert J_F\vert  \; \widetilde{W}_{{\cal M},\epsilon}^{\star}\vert_{B'\times \C^{k+1}}$ is $dd^c_{\la , \tilde z} $-closed since   $\log_\epsilon  \vert J_F\vert=\log \vert J_F\vert$ is pluriharmonic on the graphs $\Gamma_{\gamma}$ which do not intersect $U_\epsilon$. It thus remains to check that 
$\lim_{\epsilon} \log_\epsilon  \vert J_F\vert \; \widetilde{W}_{{\cal M},\epsilon}=0$. This follows from the estimate
\begin{eqnarray*}
\Vert \log_\epsilon  \vert J_F\vert  \; \widetilde{W}_{{\cal M},\epsilon}\Vert \; \lesssim\; \vert \log \epsilon\vert {\cal M} \left(S_{{\cal M},\epsilon}\right) \; \lesssim\;
\epsilon^a \vert \log \epsilon\vert
\end{eqnarray*}
where the last inequality is obtained after having observed that there exist $b,\beta>0$ such that  $S_{{\cal M},\epsilon}\subset \{\gamma\in {\cal J}\;\colon\; \Gamma_{\gamma\vert_{B'}}\cap \left(C_f\right)_{b\epsilon^\beta} \ne \emptyset\}$  and applying Lemma \ref{lemEstiZ}. \finsec

%%%%%%%%%%%%%%%%%%%%%%%%%%%%%%%%
\subsection{Misiurewicz parameters belong to $\supp dd^c_\la  L$}\label{SecMis}
%%%%%%%%%%%%%%%%%%%%%%%%%%%%%%%%

We establish here the following result.

\begin{prop}\label{thmmisiu}
Let $f:M\times \Pj^{k}\to M\times \Pj^{k}$ be a holomorphic family
of  endomorphisms of $\Pj^k$. Then the Misiurewicz parameters belong to the support of $dd^c_\la  L$.
\end{prop}

The proof relies on an infinitesimal transfer mechanism which was first used by Buff and Epstein \cite{BuEp} in the context of rational maps.

\proof
If $\lambda_0 \in M$ is a Misiurewicz parameter then, by definition, there exists a holomorphic map $\gamma $ from a neighbourhood of $\lambda_0$ into $\Pj^{k}$ such that:
\begin{itemize}
\item[1)] $\gamma(\lambda)\in J_\lambda$ and is a repelling $p_0$-periodic point of $f_{\lambda}$ for some $p_0 \geq 1$,
\item[2)] $(\lambda_0,\gamma(\lambda_0))\in f^{n_0}(C_f)$ for some $n_0 \geq 1$,
\item[3)] the graph $\Gamma_\gamma$ of $\gamma$ is not contained in $f^{n_0}(C_f)$.
\end{itemize}

Without loss of generality  we may assume that $p_0=1$ and that $M$ is a disc $D_\rho \subset \C$ centered at $\la_0=0$ with radius $\rho$.
Moreover, conjugating by  $(\la,z)\mapsto (\la,T_{\gamma(\la)}(z))$ where $T_{\gamma(\la)}$ is a suitable family of linear automorphisms of $\Pj^k$  ensures that $\gamma$ is constant equal to $z_1 := \gamma(0)$. Let us denote by $B_r$ a ball centered at $z_1$ and of radius $r$. Taking $\rho$ and $r$ sufficiently small finally allows us to suppose that:
\begin{itemize}
\item[(i)] $f$ is injective and uniformly expanding on $D_{\rho}\times B_r$: there exists  $K >1$ such that $$\forall 
(\la,z)\in D_{\rho}\times B_r ,\;  d_{\Pj^k}\left(f (\la, z),f(\la, z_1)\right)\ge Kd_{\Pj^k}(z,z_1) $$ 
\item[(ii)] $(\la,z_1)\in f^{n_0}(C_f)\Leftrightarrow \la=0.$
\end{itemize}
The fact that $\gamma(\lambda)\in J_\lambda$ is crucial but will only be used at the very end of the proof.\\

We have to show that $\langle dd^c_\la L, 1_{D_{\epsilon}}\rangle > 0$ for some $0<\epsilon<\rho$. To this purpose, we will use the formula $dd^c_\la L= (\pi_{D_{\rho}})_{\star} \left(  \left(  dd^c_{\la , z}  g +\omega  \right) ^{k} \wedge [C_f]\right)$ given by Theorem \ref{thmFBB}, where $\omega := \omega _{FS}$. Let   $(g_n)_n$ be a sequence of smooth functions on $\Pj^k$ which converges uniformly to $g$ and satisfies $\frac{1}{d}f^{\star}(dd^c_{\la , z}  g_n + \omega)= dd^c_{\la , z}  g_{n+1} + \omega$ (see subsection \ref{someFF}). We shall proceed in three steps. \\

{\bf First step:} $\langle dd^c_\la  L, 1_{D_{\epsilon}}\rangle\ge d^{-n_0 k}\langle [f^{n_0}(C_f)]\wedge \left(dd^c_{\la , z}  g+\omega\right)^{k}, 1_{D_{\epsilon}\times B_r}\rangle.$\\

Pick $(0,z_0)\in C_f$ such that $f^{n_0}(0,z_0)=(0,z_1)$. After reducing $\epsilon$ and $r$, we may find a neighbourhood $U$ of $(0,z_0)$ such that the map
$f^{n_0}: U\to D_{\epsilon}\times B_r$ is proper. According to Theorem \ref{thmFBB}, we have
\begin{eqnarray*}
\langle dd^c_\la  L, 1_{D_{\epsilon}}\rangle = \langle \left(dd^c_{\la , z}  g +\omega\right)^{k}\wedge [C_f] , 1_{D_{\epsilon}}\circ \pi_{D_{\rho}}\rangle
\ge \langle \left(dd^c_{\la , z}  g + \omega\right)^{k}\wedge [C_f], 1_{U}\rangle.
\end{eqnarray*}
Using the smooth approximations $g_n$, we get 
\begin{align*}
 \langle \left(dd^c_{\la , z}  g +\omega\right)^{k}\wedge [C_f], 1_{U}\rangle & = \lim_n \langle \left(dd^c_{\la , z}  g_{n+n_0}+\omega\right)^{k} \wedge [C_f], 1_{U}\rangle \\
																					& =\lim_n d^{-n_0 k}\langle 1_U \cdot [C_f], (f^{n_0})^{\star}
																					\left(dd^c_{\la , z}  g_{n}+\omega \right)^{k}\rangle \\
																					& =\lim_n d^{-n_0 k}\langle (f^{n_0})_{\star}\left(1_U \cdot [C_f] \right), \left(dd^c_{\la , z}  g_{n} +\omega\right)^{k}\rangle.
\end{align*}
Now, as $f^{n_0}: U\to D_{\epsilon}\times B_r$ is proper, one has $(f^{n_0})_{\star}\left(1_U  [C_f]\right) \geq 1_{D_{\epsilon}\times B_r}  [f^{n_0}(C_f)]$ which, since $dd^c_{\la , z}  g_n +\omega$ is positive, yields
\begin{align*} \label{ineg}
\langle dd^c_\la  L, 1_{D_{\epsilon}}\rangle  & \geq \lim_n d^{-n_0 k}\langle 1_{D_{\epsilon}\times B_r}  [f^{n_0}(C_f)], \left(dd^c_{\la , z}  g_{n} +\omega\right)^{k}\rangle \\
                                               & \geq  \lim_n d^{-n_0 k} \langle \left(dd^c_{\la , z}  g_{n}+\omega\right)^{k} \wedge [f^{n_0}(C_f)], 1_{D_{\epsilon}\times B_r} \rangle.
\end{align*}
The desired estimate follows by uniform convergence of $g_n$ to $g$.\\

{\bf Second step:} We set $A_0:= 1_{D_{\epsilon}\times B_r} [f^{n_0}(C_f)]$ and  $A_{p+1}:=1_{D_{\epsilon}\times B_r} f_{\star}(A_p)$. Then 
\begin{eqnarray*}
\Vert A_p\wedge \left(dd^c_{\la , z}  g +\omega\right)^{k}\Vert=d^{p k}\Vert \big(1_{D_{\epsilon}\times B_r}\circ f^p\big)A_0\wedge \left(dd^c_{\la , z}  g +\omega \right)^{k}\Vert \\
 \le d^{p k}\Vert A_0\wedge \left(dd^c_{\la , z}  g +\omega \right)^{k}\Vert.
 \end{eqnarray*}

To prove this mass-estimate, we use again the smooth approximations $g_n$. 

\begin{align*}
\Vert A_{p+1}\wedge \left(dd^c_{\la , z}  g_n +\omega\right)^{k}\Vert & = \langle 1_{D_{\epsilon}\times B_r} f_{\star}(A_p), \left(dd^c_{\la , z}  g_n+\omega\right)^{k}\rangle \\
																		& = \langle A_p, f^{\star} \left(1_{D_{\epsilon}\times B_r}\left(dd^c_{\la , z}  g_n+\omega\right)^{k}\right)\rangle \\
																& = d^{k}\langle A_p, 1_{D_{\epsilon}\times B_r}\circ f \left(dd^c_{\la , z}  g_{n+1}+\omega\right)^{k}\rangle \\
																& = d^{k}\langle A_p\wedge \left(dd^c_{\la , z}  g_{n+1}+\omega\right)^{k}, 1_{D_{\epsilon}\times B_r}\circ f \rangle \\
																& = d^{k} \Vert \big(1_{D_{\epsilon}\times B_r}\circ f\big) A_p\wedge \left(dd^c_{\la , z}  g_{n+1}+\omega\right)^{k}  \Vert.
\end{align*}
Taking the limits when $n$ tends to infinity yields the conclusion.\\

{\bf Third step:} $\langle dd^c_\la  L, 1_{D_{\epsilon}}\rangle >0$.\\

By combining the two former steps, one gets:

\begin{eqnarray}\label{trtr}
d^{(p+n_0) k}\langle dd^c_\la  L, 1_{D_{\epsilon}}\rangle \ge \Vert A_p\wedge \left(dd^c_{\la , z}  g +\omega\right)^{k}\Vert. 
\end{eqnarray}

By (i) and (ii), $f$ is uniformly expanding on $D_{\rho}\times B_r$  and $\left(\supp A_0\right) \cap \left(D_{\rho} \times \{z_1\}\right) =\{(0,z_1)\}$. Thus $\supp A_p \subset D_{\epsilon_p}\times B_r$ for some $\epsilon_p \to 0$. Let us momentarily admit that there exists $m >0$ such that
\begin{equation}\label{exim}
A_p \to  m \, [\{0\}\times B_r].
\end{equation}
We then deduce from (\ref{trtr}) that, for $p$ large enough, one has:
$$d^{(p+n_0) k} \langle dd^c_\la  L, 1_{D_{\epsilon}}\rangle \ge {m\over 2} \left\Vert [\{0\}\times B_r]  \wedge \left(dd^c_{\la , z}  g +\omega\right)^{k}\right\Vert .$$
We conclude by using the fact that $z_1\in  J_0$: the right hand side is equal to 
$$ {m\over 2} \int_{B_r} \left(dd^c_{z}  g(0,z) +\omega \right)^{k} = {m\over 2} \mu_0 (B_r) > 0. $$
To complete the proof it remains to establish (\ref{exim}). Let us denote $V:=D_{\rho}\times B_r$ and $V':=f(V)$. By assumption 
$f:V\to V'$ is a biholomorphism whose inverse will be denoted by $h:V'\to V$. According to (i),  $V \subset V'$ and $\left(h\vert_{V}\right)^p$ converges to $(\la,z)\mapsto (\la, z_1)$. We now use (ii). After shrinking $\rho$ and $r$, we may find a Weierstrass polynomial 
$$ \psi(\la,z):=\la^m+\alpha_{m-1}(z)\la^{m-1}+\cdots+\alpha_0(z) $$
such that $\alpha_j(z_1)=0$ for $0\le j\le m-1$ and 
$f^{n_0}\left(C_f\right) \cap \left(D_{\rho}\times B_r\right) = \{\psi=0\}$. Observe now that $A_0=1_V dd^c_{\la, z}  \log \vert \psi \vert$ and that
$$ A_1= 1_V  f_{\star} A_0  =1_V h^{\star} A_0=
1_V (1_V\circ h)  dd^c_{\la, z}  \log \vert \psi\circ h\vert=1_V dd^c_{\la, z}  \log \vert \psi\circ h\vert , $$ 
where the last equality comes from $h(V) \subset V$.
Similarly we have $A_p= 1_V dd^c_{\la, z}  \log \vert \psi \circ \left(h\vert _V \right)^p\vert$  and the conclusion follows since
$\psi\circ \left(h\vert _V\right)^p(\lambda,z) \to \la^m$.\finsec

%%%%%%%%%%%%%%%%%%%%%%%%%%%%%%%%%%%%%%%%%%%%%
\subsection{ Misiurewicz parameters are dense in $\supp dd^c_\la L$}
%%%%%%%%%%%%%%%%%%%%%%%%%%%%%%%%%%%%%%%%%%%%%

We start with the following proposition; the statement is local since it is based on holomorphic motion of hyperbolic sets.

\begin{prop}\label{PropHmotMis}
Let  $f : B \times \Pj^k \to B \times \Pj^k$ be a holomorphic family of endomorphisms of $\Pj^k$ where 
$B$ is a ball centered at the origin in $\C^m$. If there is no  Misiurewicz parameter in $B$ then, after shrinking $B$,  there exists $\gamma \in {\cal J}$ for which $\Gamma_{\gamma}$ does not intersect the post-critical set of $f$. 
\end{prop}

Every hyperbolic set admits a holomorphic motion which preserves repelling cycles (see subsection \ref{HMHS}). We need a more precise result
concerning the size of such sets and the position of their motions with respect to Julia sets.
Here $B_r$ denotes a ball centered at the origin in $\C^m$ and of radius $r$.

\begin{thm}\label{gtg}
Let $f : B \times \Pj^k \to B \times \Pj^k$ be a holomorphic family of endomorphisms of $\Pj^k$. 
There exist an integer $N$, a compact hyperbolic $f^N_0$-invariant set $E_0 \subset J_0$  and a holomorphic motion $h :B_r\times E_0\to \Pj^k$ 
for some $0 < r < 1$ such that:
\begin{enumerate}
\item the repelling periodic points of $f^N_0$ are dense in $E_0$ and 
$E_0$ is not contained in the post-critical set of $f^N_0$,
 \item $h_\la(z)\in J_\la$ for every $\la\in B_r$ and every $z\in E_0$,
\item if $z$ is periodic repelling for $f^N_0$ then $h_\lambda(z)$ is periodic repelling for $f^N_\lambda$.
\end{enumerate}
 \end{thm}

 The proof of this result requires a few tools.
 To create hyperbolic sets, we use a classical device based on  the following proposition which is a consequence of \cite{BD1} (see also \cite{BDM}). For any  endomorphism $f_0$ of $\Pj^k$ and  every $A \subset \Pj^k$, $n \geq 1$ and $\rho > 0$, we denote by $C_n(A,\rho)$ the set of inverse branches $g_i$ of $f_0^{n}$ defined on $A$ and satisfying $g_i (A) \subset A$ and $\Lip g_i \leq \rho$. 

\begin{prop}\label{cab}
Let $f_0$ be an endomorphism of $\Pj^k$ of degree $d$.
For every $\rho > 0$ there exist a closed ball $A \subset \Pj^k$ centered on $ J_{f_0}$ and $\alpha > 0$ such that $\cd \, C_n(A,\rho) \geq \alpha d^{kn}$.
\end{prop}

To control the size of hyperbolic sets, we use an entropy argument. 
 Our key tool is the following result which is due to Briend-Duval \cite{BD2}, de Th\'elin \cite{dT} and Dinh \cite{Di3}
 (see also \cite{DS2} Corollary 1.117).

\begin{thm}\label{ggh}
Let $g$ be an endomorphism of $\Pj^k$  of degree $d$.
Let $\kappa$ be an ergodic $g$-invariant measure with entropy $h_\kappa > (k-1)\log d$.
Then $\kappa$ gives no mass 
to analytic subsets of dimension $\leq k-1$ and 
the support of $\kappa$ is included in the Julia set of $g$.
\end{thm}

\textsc{Proof of Theorem \ref{gtg}:} Let $\rho < 1$ and $A$ be a closed ball provided by Proposition \ref{cab}. Let us fix $N$ large enough such that $N':=\cd \, C_N (A,\rho) > d^{(k-1)N}$. We denote by $g_1,\ldots,g_{N'}$ the elements of $C_N(A,\rho)$. Let  $E_0:= \cap_{k\geq 1} E_k$, where
\[ E_k := \big \{ \ g_{i_1} \circ \ldots \circ g_{i_k} (A)\;\colon\; (i_1, \ldots , i_k) \in \{ 1 , \ldots , N' \}^k \ \big  \} .  \]
 Let $\Sigma := \{ 1 , \ldots , N' \}^{\N^*}$ endowed with the product metric and $z$ be a fixed point in $A\cap J_0$, for instance the center of $A$. The map $\omega : \Sigma \to E_0$ defined by $(i_1,i_2,\ldots) \mapsto \lim_{k \to \infty} g_{i_1} \circ \ldots \circ g_{i_k} (z)$ is a homeomorphism satisfying $f^N \circ \omega = \omega \circ s$, where $s$ is the left shift acting on $\Sigma$.  We take for $\kappa$ the image by $\omega$ of the uniform product measure on $\Sigma$:  this is a $f^N$-invariant ergodic measure with entropy $h_\kappa = \log N'>(k-1)\log d^n$, with support $E_0$.
 
 By construction $E_0\subset J_{f_0}$.  Indeed,  $E_0 =  \{ \, \lim_{k \to \infty} g_{i_1} \circ \ldots \circ g_{i_k} (z)\;\colon\; (i_1, i_2, \ldots  ) \in \Sigma \, \}$ and $ J_{f_0}$ is a closed $f_0^N$-invariant set. Also, repelling cycles of $f_0^N$ are dense in $E_0$.
 According to Theorem \ref{ggh}, $E_0=\supp \kappa$ is not contained in the countable union of analytic subsets $\cup_{n \ge 1} f_0^n(C_{f_0})$. 
 The set $E_0$ is hyperbolic for $f_0^N$ since $\vert (d f_0 ^N)^{-1} \vert^{-1} > \frac{1}{\rho} > 1$ on $E_0$ and thus there exists a holomorphic motion $h :B_r\times E_0\to \Pj^k$ 
 which preserves repelling cycles
 (see Theorem \ref{hypmot}). It remains to show $h_\lambda(E_0) \subset  J_{f_\lambda}$. For that purpose we use the fact that $h_\lambda : E_0 \to \Pj^k$ is a continuous injective mapping satisfying $h_\lambda \circ f_0^N =  f_\lambda^N \circ h_\lambda$ on 
 $E_0$. Then $(h_\lambda)_* \kappa$ is a $f_\lambda^N$-invariant ergodic measure whose support coincides with $h_\lambda (E_0)$ and whose metric entropy equals $h_\kappa$. Theorem  \ref{ggh} yields $h_\lambda(E_0) \subset  J_{f_\lambda}$ as desired. \fin

We  now use Theorem \ref{gtg} to establish Proposition \ref{PropHmotMis}.\\

\textsc{Proof of  Proposition \ref{PropHmotMis}:}  Let $E_0 \subset  J_{0}$ and $r \in ]0,1]$, $N\in \N$ provided by Theorem \ref{gtg}. Since $f^N_\la$ and $f_\la$ have same equilibrium measures and post-critical sets, we may assume that $N=1$.  Let us fix $z \in E_0 \setminus \cup_{n \geq 1} f_0^n (C_{f_0})$ (see item 1).\\
 Let us set $\gamma(\la):= h_\la(z)$.  By item 2 we have $\gamma\in {\cal J}$. Let us show that
\begin{eqnarray} \label{Cgam}
\Gamma_\gamma \cap \big( \cup_{n\ge 1} f^{n}(C_f)\big)  = \emptyset.
\end{eqnarray}
Assume to the contrary that there exists $n_0 \geq 1$ such that $\Gamma_\gamma \cap f^{n_0}(C_f) \ne \emptyset$. Note that $\gamma(0)\notin f^{n_0}\left(C_f\right)$. By item 1, there exists a sequence $(z_p)_p \subset E_0$ of $f_0$-periodic repelling points which converges to  $z$. Items 2 and 3 assert that $h_\la(z_p) \in  J_\la$ and  $h_\la(z_p)$ is a $f_\la$-periodic repelling point for every $\lambda \in B_r$. As $h$ is continuous,  $\lambda \mapsto h_\la (z_p)$ converges locally uniformly to $\lambda \mapsto h_\la (z)=\gamma(\la)$. Hence, for $p$ large enough, the graph $\{(\la, h_{\la}(z_p))\;\la\in B_r\}$
  is not contained in $f^{n_0}(C_f)$  (consider the parameter $\la=0$) and, by Hurwitz's lemma,  there exists $\la_p \in B_r$ such that $(\la_p, h_{\la_p}(z_p)) \in f^{n_0}(C_f)$.  The parameters $\la_p$ are Misiurewicz, contradicting our assumption. \fin

We can now prove Theorem \ref{main} which, in particular, says that Misiurewicz parameters are dense in the support of $dd^c L$.\\

\textsc{Proof of Theorem \ref{main}}: 
By Proposition \ref{thmmisiu} there are no Misiurewicz parameters in $M$ if $dd^c_\la  L \equiv 0$ on $M$ and thus (a)$\Rightarrow$(b).
If there are no Misiurewicz parameters in $M$ then, by Propositions \ref{PropHmotMis} and \ref{propgraph},  for any parameter $\la$
one an find an open open ball  $B$ centered at $\la$ such that the restriction $f\vert_{B\times \Pj^k}$  admits an equilibrium  web ${\cal M}=\lim_n {\cal M}_n$ satisfying   $\Gamma_\gamma \cap C_f=\emptyset$ for any $\gamma\in \cup_n \supp {\cal M}_n$. Thus (b)$\Rightarrow$(c). Finally, (c)$\Rightarrow$(a) follows from Proposition \ref{thmHMR}.\finsec

%%%%%%%%%%%%%%%%%%%%%
\subsection{Proofs of part of Theorem \ref{main2} and Corollary \ref{CorASW}}\label{sspartof}
%%%%%%%%%%%%%%%%%%%%%

Let  $f:M\times \Pj^{k}\to M\times \Pj^{k}$ be a holomorphic family of endomorphisms of $\Pj^k$.
We first establish the implication (A)$\Rightarrow$(B) in Theorem \ref{main2}.
If the repelling $J$-cycles of $f$ move holomorphically then, using  the second assertion of Proposition \ref{propgraph}, 
 one gets an equilibrium  web $\cal M$ for $f$
such that ${\cal M}=\lim_n {\cal M}_n$ and  $\Gamma_\gamma \cap C_f=\emptyset$ for any $\gamma\in \cup_n \supp {\cal M}_n$.
By Proposition \ref{thmHMR}, this implies that $dd^c_\la  L\equiv 0$ on $M$. This justifies (A)$\Rightarrow$(B). In the spirit of the proposition 1.26 of  \cite{DS2} concerning the Julia set of a single endomorphism of $\Pj^k$, we have the following proposition. It implies the equivalence (B)$\Leftrightarrow$(E).
 
\begin{prop}\label{growth} 
Let $B$ be an open ball in $\C^m$ and let $f : B \times \Pj^k \to B \times \Pj^k$ be a holomorphic family of endomorphisms of $\Pj^k$ of degree $d$. We endow $B \times \Pj^k$ with the metric $dd^c_\la \vert \lambda \vert^2 + \omega_{FS}$ and denote $\vert \cdot \vert_U$ the mass of currents in $U \times \Pj^k$. The following properties are equivalent.
\begin{enumerate}
\item $\lambda_0 \in \supp dd^c_\lambda L$.
%\item $\vert   {\cal E}_{Green} \wedge [C_f] \vert_U > 0$ for every neighbourhood $U$ of $\lambda_0$.
\item $\liminf_n d^{-k n}  \vert { (f^n)_*  [C_f]}\vert_U > 0$  for every  neighbourhood $U$ of $\lambda_0$.
\item $\limsup_n d^{-(k-1) n} \vert { (f^n)_*  [C_f]} \vert_U  = + \infty$ for every neighbourhood $U$ of $\lambda_0$.
\end{enumerate}
\end{prop}

\proof The equivalence between 1. 2. and 3.
is a direct consequence of the following Lemma. \finsec

\begin{lem}\label{grgr}
There exists $\alpha = \alpha(k,m) >0$ such that, for every compact subset $U \subset  M$:  
$$ \vert {(f^n)_* [C_f]} \vert_U  =  \alpha \, d^{k n} \vert dd^c L  \vert_U + O(d^{(k-1)n})  . $$ 
\end{lem}

\proof Let us set $\kappa:=k+m-1$. Then  
$$ \vert  (f^n)_* [C_f] \vert_U =  \int_{U \times \Pj^k} (f^n)_* [C_f] \wedge [{\omega_{FS}} + dd^c_\la \vert \lambda \vert^2]^{\kappa} = \int_{U \times \Pj^k} [C_f] \wedge (f^n)^* [{\omega_{FS}} + dd^c_\la \vert \lambda \vert^2]^{\kappa} . $$
Using $\omega_{FS}^{k+1}=0$, we obtain $[{\omega_{FS}} + dd^c_\la \vert \lambda \vert^2]^{\kappa} =  \sum_{j=0}^{k} \alpha_j \, \om_{FS}^j \wedge (dd^c_\la \vert \lambda \vert^2)^{\kappa -j}$, where the $\alpha_j$'s are positive numbers. Since $\pi_M \circ f = \pi_M$, we obtain
$$ (f^n)^*  [{\omega_{FS}} + dd^c_\la \vert \lambda \vert^2]^{\kappa} = \sum_{j=0}^{k} \alpha_j \, \left(  (f^n)^*\om^j_{FS} \right) \wedge (dd^c_\la \vert \lambda \vert^2)^{\kappa -j} . $$
 Let ${\cal T} := dd_{\lambda,z}^c g + \omega_{FS}$ so that ${\cal T}^k = {\cal E}_{Green}$. By the functoriality $f^*{\cal T} = d{\cal T}$ we get
$(f^n)^*({\omega^j_{FS}}) = (d^n {\cal T} - dd^c_{\lambda,z} g \circ f^n)^j$. Now, using the fact that $g$ is bounded,  by extracting the $k$-th term of the preceding sum  we obtain:
$$(f^n)^* [{\omega_{FS}} + dd^c_\la \vert \lambda \vert^2]^{\kappa} = \alpha_k \, d^{kn} \, {\cal T}^k \wedge (dd^c_\la \vert \lambda \vert^2)^{m-1} + O(d^{(k-1)n}).$$
Then 
\begin{eqnarray*}
d^{-kn} \vert  (f^n)_* [C_f] \vert_U = \alpha_k 
\int_{U \times \Pj^k}  [C_f] \wedge {\cal T}^k \wedge ( dd^c_\la \vert \lambda \vert^2)^{m-1} + O(d^{-n})=\\
\alpha_k  \int_{U \times \Pj^k}  {\cal E}_{Green}\wedge [C_f] \wedge (\pi_B^{\star} dd^c_\la \vert \lambda \vert^2)^{m-1} + O(d^{-n})
=\alpha_k \vert dd^c L\vert_U + O(d^{-n})
\end{eqnarray*}
where the last equality comes from  Theorem \ref{thmFBB}, which asserts that $dd^c_\lambda L = \pi_{B \star} \left( \Eg\wedge [C_f] \right)$. 
We set $\alpha := \alpha_k$. This completes the proof of the lemma. \fin

\textsc{Proof of Corollary \ref{CorASW}}: By assumption, for every  $n\ge 1$ we have subsets ${\cal R}_n:=\{\rho_{n,j} \;\colon \; 1\le j\le N_n\}$ of $\cal J$ such that the 
$\rho_{n,j}(\la)$ are repelling $n$-periodic points of $f_\la$ for every $\la\in M$. Note that $\lim_n d^{-kn}N_n =1$.  We define a sequence $\left({\cal M}_n\right)_n$ of 
$\cal F$-invariant discrete probability measures on $\cal J$ by setting ${\cal M}_n:=\frac{1}{N_n} \sum_{j=1}^{N_n} \delta_{\rho_{n,j}(\la) }$.
  According to the second assertion of Proposition \ref{propgraph}, 
$\left({\cal M}_n\right)_n$  converges to an equilibrium  web $\cal M$ after taking a subsequence. Moreover, there exists a compact subset $\cal K$ of $\cal J$ such that  $\supp {\cal M}_n \subset {\cal K}$ for every $n\ge 1$.\\
Let us now prove that ${\cal M}\left({\cal J}_s\right)=0$. By the implication (A)$\Rightarrow$(B) of Theorem \ref{main2} we have  $dd^c L=0$ and then Theorem \ref{main} implies that $M$ does not contain Misiurewicz parameters.
We can now see  that for every $k \in \N$ and every $\gamma\in \supp \cal M$ one has:
$$\Gamma_\gamma \cap f^k(C_f) \ne \emptyset \Rightarrow \Gamma_\gamma\subset f^k(C_f).$$
 Indeed, if this were not  the case, by Hurwitz theorem, we could find some $\gamma' \in \cup_n \supp {\cal M}_n$ such that $\Gamma_\gamma \cap f^k(C_f) \ne \emptyset$ and $\Gamma_\gamma$ is not contained in $f^k(C_f)$.
When $k=0$ this is clearly impossible since $\gamma'(\la)$ is a repelling cycle of $f_\la$ and when $k\ge 1$, this is impossible because $M$ does not contain
Misiurewicz parameter.\\
So, fixing any $\la_0\in M$, we get 
\begin{eqnarray*}
{\cal M} \left(\{ \gamma\in {\cal J}\;\colon\; \Gamma_\gamma\cap \left(\cup_{k\ge 0} f^k(C_f)\right) \ne \emptyset\}\right)
= {\cal M} \left(\{ \gamma\in {\cal J}\;\colon\; \Gamma_\gamma\subset \left(\cup_{k\ge 0} f^k(C_f)\right) \}\right) \le\\
 {\cal M} \left(\{ \gamma\in {\cal J}\;\colon\; (\la_0,\gamma(\la_0)) \in 
 \left(\cup_{k\ge 0} f^k(C_f)\right) \}\right)=
 \mu_{\la_0}  \left(\cup_{k\ge 0} f_{\la_0}^k(C_{f_{\la_0}})\right) =0
\end{eqnarray*}
where the two last equalities come from $p_{\la_0\star}\left({\cal M}\right)=\mu_{\la_0}$ and the fact that $\mu_{\la_0}$ does not charge pluripolar sets in $\Pj^k$. The estimate ${\cal M}\left({\cal J}_s\right)=0$
 follows from the $\cal F$-invariance of $\cal M$.
Finally,  Proposition \ref{PropGE} shows that there exists an ergodic equilibrium  web ${\cal M}_0$ such that
${\cal M}_0\left({\cal J}_s\right)=0$.  \finsec

%%%%%%%%
\section{From equilibrium  webs to equilibrium laminations}\label{secFrom}
%%%%%%%%%%

Our goal here is to establish  the implication (A)$\Rightarrow$(D) in Theorem \ref{main2}. We prove the following more precise result.

\begin{thm}\label{PropExHM}
Let $M$ be a simply connected complex  manifold and $f:M\times\Pj^k\to M\times\Pj^k$ be a holomorphic family of  endomorphisms of $\Pj^k$ of degree $d\ge 2$. If the repelling $J$-cycles of $f$ move holomorphically over $M$ or if $f$ admits an acritical and ergodic equilibrium  web then there exists an equilibrium lamination $\cal L$ for $f$. Moreover, $f$ admits a unique equilibrium  web $\cal M$ and ${\cal M}\left({\cal L}_1 \Delta {\cal L}_2\right)=0$
for any pair of equilibrium laminations ${\cal L}_1,{\cal L}_2$ of $f$.
\end{thm}

Given an acritical and ergodic equilibrium  web $\cal M$ of $f$, our strategy will consist in first proving that the iterated inverse branches in $\left({\cal J},{\cal F},{\cal M}\right)$ are exponentially contracting
 and then exploit this property to extract an equilibrium lamination from the support of $\cal M$.
 By totally different methods,
Berger and Dujardin (\cite{BgDj}) have recently built measurable holomorphic motions
 in the context of polynomial automorphisms of $\mathbb{C}^2$.

%%%%%%%%%%%%
\subsection{On the rate of contraction of iterated inverse branches in $\left({\cal J},{\cal F},{\cal M}\right)$}
  %%%%%%%%%%%%
We explain here  how certain  stochastic properties of the system $\left({\cal J}, {\cal F}, {\cal M}\right)$ allow to control the rate of contraction
 of the iterated inverse branches of ${\cal F}$ (see Proposition \ref{LemSF}).  We adapt to the context of  $({\cal J},{\cal F},{\cal M})$ the tools which have been first introduced in \cite{BD1} by Briend-Duval for the case of a single holomorphic endomorphism of $\Pj^k$.
 Let us stress however that new arguments will be introduced in subsection \ref{SSestiLyap}.\\
 
Since all our statements here are local we may assume that the parameter space $M$ is a simply connected open subset of $\C^m$ which we endow with the euclidean norm.\\

 To study the inverse branches of the map $\cal F$, it is  convenient to transform the system $({\cal J},{\cal F},{\cal M})$ into an injective one. This is possible using a classical construction called the \emph {natural extension} which we now describe (we refer to \cite{CFS} page 240 for more details).

Recall that ${\cal K}:=\supp {\cal M}$ is a compact subset of $\cal J$ and that ${\cal M}\left({\cal J}_s\right)=0$.
Setting ${\cal X}:={\cal K}\setminus{\cal J}_s$, it is not difficult to check  that the map ${\cal F} : {\cal X}\to{\cal X}$ is onto. We may therefore construct the natural extension $\left(\widehat{\cal X},\widehat{\cal F},
\widehat{\cal M}\right)$ of the system $\left({\cal X}, {\cal F}, {\cal M}\right)$ in the following way.
An element of $\widehat{\cal X}$ is a bi-infinite sequence $\widehat{\gamma}:=(\cdots,\gamma_{-j},\gamma_{-j+1},\cdots,\gamma_{-1},\gamma_0,\gamma_1,\cdots)$ of elements $\gamma_{j}\in {\cal X}$ such that ${\cal F}(\gamma_{-j})=\gamma_{-j+1}$ and one defines the map $\widehat{\cal F} : \widehat{\cal X}\to\widehat{\cal X}$ by setting
\[
\widehat{\cal F}(\widehat {\gamma}):=(\cdots {\cal F}(\gamma_{-j}),{\cal F}(\gamma_{-j+1})\cdots).
\]

The map $\widehat{\cal F}$ corresponds to the shift operator and  is clearly bijective. There exists a unique measure
 $\widehat{\cal M}$ on $\widehat{\cal X}$ such that
\[(\pi_{j})_\star\left(\widehat{\cal M}\right)={\cal M}\]
for any projection $\pi_j:\widehat{\cal X}\to\widehat{\cal X}$ given by
$\pi_j(\widehat \gamma)=\gamma_{j}.$
The ergodicity of $\cal M$ implies the ergodicity of $\widehat{\cal M}$. We have thus obtained an invertible and ergodic dynamical system
$\left(\widehat{\cal X},\widehat{\cal F},
\widehat{\cal M}\right)$.

For every $\gamma \in {\cal J}$ whose graph $\Gamma_\gamma$ does not meet the critical set of $f$, we denote by $f_\gamma$ the injective map which is induced by $f$ on some neighbourhood of
 $\Gamma_\gamma$ and by $f^{-1}_{\gamma}$ the inverse branch of $f_\gamma$ which is defined on some neighbourhood of $\Gamma_{{\cal F}(\gamma)}$. Thus, given $\widehat \gamma \in \widehat{\cal X}$ and  $n\in \N$ we may define the iterated inverse branch $f^{-n}_{\widehat \gamma}$ of $f$ along
 $\widehat \gamma$ and of depth $n$ by
 $$
 f^{-n}_{\widehat \gamma}:=f^{-1}_{\gamma_{-n}}\circ\cdots \circ f^{-1}_{\gamma_{-2}}\circ f^{-1}_{\gamma_{-1}}.
$$

Let us stress that $ f^{-n}_{\widehat \gamma}$ is defined on a neighbourhood of $\Gamma_{\gamma_0}$ with values in a 
 neighbourhood of $\Gamma_{\gamma_{-n}}$. Moreover, since only a finite number of components of the grand critical orbit of $f$ are involved for 
 defining $ f^{-n}_{\widehat \gamma}$, we may always shrink the parameter space $M$ to some $\Omega\Subset M$ so  that the domain of definition of  $ f^{-n}_{\widehat \gamma}$ for  fixed $n$ and $\widehat \gamma$ 
contains a tubular neighbourhood of $\Gamma_{\gamma_0}\cap \left(\Omega\times \Pj^k\right)$ of the form 
 $$T_\Omega(\gamma_0,\eta):=\{(\la,z)\in \Omega\times\Pj^k\;\colon\; d_{\Pj^k}(z,\gamma_0(\la)) <\eta\}.$$

Our goal is to get a uniform $\eta$, independent from $n$, and to control the size of $f^{-n}_{\widehat \gamma}\left(T_{U_0}(\gamma_0,\widehat{\eta}_p (\widehat \gamma)\right))$ for suitable $\widehat{\eta}_p (\widehat \gamma)>0$ and $U_0\subset M$. We will now explain how this  boils down to estimating some kind of Lyapunov exponent.
This requires however to first introduce a few more notations.\\

To start with, we need to fix sets of holomorphic charts with bounded distorsions on $\Pj^k$. 
For any $\tau >0$, there exists a covering $\Pj^k=\cup_{i=1}^N V_i$ by open sets and a collection of holomorphic maps 
$$\psi_i : V_i \times B_{\C^k} (0,R_0) \to \Pj^k$$
such that $\psi_{i,x}:=\psi_i(x,\cdot)$ is a chart of $\Pj^k$ satisfying $\psi_{i,x}(0)=x$ and 
\begin{eqnarray}\label{LipChart}
e^{-\tau /2} \vert z-z'\vert \le d_{\Pj^k}\left( \psi_{i,x}(z), \psi_{i,x}(z')\right) \le e^{\tau/2} \vert z-z'\vert
\end{eqnarray}
for  every $(x,z)\in V_i \times B_{\C^k}(0,R_0)$ and every $1\le i\le N$.\\

We will now use these holomorphic charts to express the restrictions of $f^n$ on suitable neighbourhoods of graphs $\Gamma_\gamma$.
Let us fix $\la_0$ in $M$. Since the family $\cal K=\supp {\cal M}$ is locally equicontinuous, there exists a relatively compact open ball $W_0$ centered at  $\la_0$ in $M$  such that:
$$\forall \gamma \in {\cal K},\;\exists i\in\{1,2,\cdots,N\}\;\textrm{such that}\; \gamma(\la)\in V_i \;\textrm{for all}\; \la \in  \overline {W_0}.$$
For all $\gamma \in {\cal K}$ we set 
$$i(\gamma):= \inf \{1\le i\le N\;\colon\; \gamma(\la)\in V_i\;\textrm{for all}\; \la \in  \overline{W_0}\}.$$
Then, for every $n\ge 1$ there exists  $R_n\in ]0,R_0]$ such that the maps $F_{\gamma(\la)}^n$ given by
\begin{eqnarray}\label{chart}F_{\gamma(\la)}^n:= \left(\psi_{i({\cal F}^n\gamma),({\cal F}^n \gamma)(\la)} \right)^{-1} \circ f_{\la}^n \circ \psi_{i(\gamma),\gamma(\la)}
\end{eqnarray}
are well defined and holomorphic on on a fixed neighbourhood of $\overline{W_0}\times \overline{B_{\C^k}(0,R_n)}$ for every $\gamma \in {\cal K}$.
This follows immediately from the uniform continuity of $f^n$ on $\overline{W_0}\times \Pj^k$.\\

 As 
$F_{\gamma(\la)}^n$ is locally invertible at the origin when $\gamma\notin {\cal J}_s$, we may now define  functions $u_n$ on ${\cal X}\times \overline{W_0}$ by setting
$$u_n(\gamma,\la):=\log \Vert  (DF_{\gamma(\la)}^n(0))^{-1}\Vert.$$
Let us stress that $ (DF_{\gamma(\la)}^n(0))^{-1}$ depends holomorphically on $\la \in \overline{W_0}$.\\

From now on we consider  three open balls $U_0\Subset V_0\Subset W_0$ centered at $\la_0$ in $M$.
Let us  introduce the function $r_n$  on $\cal X$ and $\widehat{u}_n$ on $\widehat{\cal X}$ by setting
\begin{eqnarray}\label{new}
r_n(\gamma) :=e^{-2\sup_{\la \in U_0} u_n(\gamma,\la)}\;\;\textrm{and}\;\;\widehat{u}_n(\widehat{\gamma}):=\sup_{\la\in U_0} u_n(\gamma_0,\la)=-\frac{1}{2} \log r_n(\gamma_0).
\end{eqnarray}

We may now state the announced result.

\begin{prop}\label{LemSF} 
 Let $f:M\times\Pj^k\to M\times\Pj^k$ be a holomorphic family of  endomorphisms of $\Pj^k$ of degree $d\ge 2$ which  admits an acritical and ergodic equilibrium  web $\cal M$. Assume that the functions $\widehat{u}_n$ are $\widehat{\cal M}$-integrable and that
 \begin{eqnarray*}
(\star)\;\;\; \lim_n\frac{1}{n} \int_{\widehat{\chi}} \widehat{u}_n\;d\widehat{\cal M}=L\;\textrm{ for some} \; L  \le - {1 \over 2} \log  d.
 \end{eqnarray*}
   Then there exist  $p \geq 1$, a Borel subset $\widehat{\cal Y} \subset \widehat{\cal X}$
such that $\widehat{\cal M} (\widehat{\cal Y} )=1$, a measurable function 
$\widehat{\eta}_p : \widehat{\cal Y}\to ]0,1]$ and a constant $A>0$ which satisfy the following properties.

For every $\widehat{\gamma}\in \widehat{\cal Y}$ and every $n\in p\N^\star$
the iterated inverse branch  $f^{-n}_{\widehat \gamma}$ is defined on the tubular neighbourhood $T_{U_0}(\gamma_0,\widehat{\eta}_p (\widehat \gamma))$ of $\Gamma_{\gamma_0}\cap (U_0\times\Pj^k)$ and 
$$f^{-n}_{\widehat \gamma}\left(T_{U_0}(\gamma_0,\widehat{\eta}_p (\widehat \gamma))\right)\subset T_{U_0}(\gamma_{-n},e^{-nA}).$$
Moreover, the map $f^{-n}_{\widehat \gamma}$ is Lipschitz with $\mbox{Lip}\;f_{\widehat{\gamma}}^{-n} \le \widehat{l}_p (\widehat{\gamma}) e^{-nA}$ where $\widehat{l}_p(\widehat{\gamma}) \ge 1$.
\end{prop}

The proof of Proposition \ref{LemSF} follows Briend-Duval \cite{BD1} and is given in the Appendix.

%%%%%%%%%%%%%%%%
\subsection{Estimating  a Lyapunov exponent}\label{SSestiLyap}
%%%%%%%%%%%%%%%%%%

The main result of this subsection is as follows; it asserts that the assumption ($\star$)  of Proposition \ref{LemSF} is satisfied when $f$ admits an acritical and ergodic equilibrium web $\cal M$.

\begin{prop}\label{ThmLya} 
Let $f:M\times\Pj^k\to M\times\Pj^k$ be a holomorphic family of  endomorphisms of $\Pj^k$ of degree $d\ge 2$ which  admits an acritical and ergodic equilibrium  web $\cal M$. Then the functions ${\widehat u}_n$ are $\widehat{\cal M}$-integrable, 
there exists a constant $L\le - {1 \over 2} \log  d$ such that 
$$\lim_n \frac{1}{n}\int_{\widehat{\cal X}} \widehat{u}_n\;d\widehat{\cal M} 
= L $$
$\textrm{and}\;\lim_n \frac{1}{n}\widehat{u}_n(\widehat{\gamma})
=L\;\textrm{for}\;\widehat {\cal M}\textrm{-almost every}\; 
\widehat{\gamma} \in {\cal X}.$\\
\end{prop}

Note that the constant $L$ may be considered as a bound for a Lyapunov exponent of the system $({\cal J},{\cal F},{\cal M})$.
The combination of Propositions  \ref{LemSF} and \ref{ThmLya} will allow us to prove Theorem \ref{PropExHM} (see Section \ref{SecProof}).\\

We keep here the assumptions and the notations introduced  in the previous subsection.
In the next Lemma, we list some basic properties of the functions $u_n$ and $\widehat{u}_n$.

\begin{lem}\label{LemPty} Let $U_0\Subset V_0\Subset W_0$ be open balls centered at $\la_0$ in $M$. 
Let $\chi_1(\la)$ be the smallest Lyapunov exponent of the system $(J_\la,f_\la,\mu_\la)$.
The functions $u_n$ and ${\widehat u}_n$ satisfy the following properties.
\begin{itemize}
\item[1)]  $u_n(\gamma,\cdot)$ is \emph{psh} on $W_0$ for every $\gamma \in {\cal X}$.
\item[2)] The sequence $({\widehat u}_n)_n$ is subadditive on $\widehat{\cal X}$, i.e., ${\widehat u}_{m+n}\le {\widehat u}_n+ {\widehat u}_m\circ \widehat{\cal F}^n$.
\item[3)] For any fixed $\la \in W_0$, we have $\lim_n \frac{1}{n} u_n(\gamma,\la) =-\chi_1(\la)$ for ${\cal M}$-almost every $\gamma\in {\cal X}$.
\item[4)] For $\cal M$-almost every $\gamma\in {\cal X}$ we have $\lim_n \frac{1}{n} u_n(\gamma,\la)=-\chi_1(\la)$ for Lebesgue-almost every $\la\in W_0$.
\end{itemize}
\end{lem}

\proof 1) When $\gamma\in {\cal X}$ is fixed the function $u_n(\gamma,\cdot)$ is clearly continuous 
on $W_0$ and $u_n(\gamma,\la)=\sup_{\Vert e \Vert =1} \log  \Vert  (DF_{\gamma(\la)}^n(0))^{-1}\cdot e\Vert$. To see that 
$u_n(\gamma,\cdot)$ is \emph{psh} it thus suffices to notice that $\la\mapsto \log  \Vert  (DF_{\gamma(\la)}^n(0))^{-1}\cdot e\Vert$ is \emph{psh} for each unit vector
$e\in \C^k$.

2)
Let $\gamma \in {\cal X}$ and $m,n\ge 1$. It follows immediately from (\ref{chart}) that
\begin{eqnarray}\label{ChRu}
\left(DF_{\gamma(\la)} ^{m+n} (0)\right)^{-1} = \left(DF_{\gamma(\la)} ^{n} (0)\right)^{-1} \circ \left(DF_{{\cal F}^n(\gamma)(\la)} ^{m} (0)\right)^{-1}\;\;\forall \la \in W_0.
\end{eqnarray}
Thus, if $\widehat{\gamma}\in \widehat{\cal X}$ we have
$\widehat{u}_{m+n}(\widehat{\gamma})  \le  \log  \sup_{\la\in {U_0}} (\Vert (DF_{\gamma_0(\la)} ^{n} (0))^{-1}\Vert \;\Vert  (DF_{{\cal F}^n(\gamma_0)(\la)} ^{m} (0))^{-1}\Vert)   \\ 
\le  \log   \sup_{\la\in {U_0}} \Vert (DF_{\gamma_0(\la)} ^{n} (0))^{-1}\Vert  +\log  \sup_{\la\in {U_0}} \Vert  (DF_{{\cal F}^n(\gamma_0)(\la)} ^{m} (0))^{-1}\Vert 
=  {\widehat u}_n(\widehat \gamma)+ {\widehat u}_m(\widehat{\cal F}^n(\widehat \gamma)).
$

3)  By Oseledec Theorem, the subset $J_{\la,1}$ of $J_\la$ defined by 
\begin{center}
$ J_{\la,1}:=\{x\in J_\la\;\colon\; \lim_n \frac{1}{n} \log  \Vert (Df_\la^n)_x^{-1}\Vert =-\chi_1(\la)\}$
\end{center}
has full $\mu_\la$ measure. As $p_{\la\star}\left({\cal M}\right)=\mu_\la$, this implies that $\gamma(\la) \in J_{\la,1}$ for $\cal M$-almost every $\gamma$ in $\cal X$. Then the assertion follows, using (\ref{LipChart}).\\

4) Let us denote by $\cal{L}eb$ the Lebesgue measure on $M$. Let $E$ be the measurable subset of ${\cal X}\times W_0$ given by
$$E:=\{(\gamma,\la)\in {\cal X}\times W_0\;\colon\; \lim_n\frac{1}{n} u_n(\gamma,\la)=-\chi_1(\la)\}.$$
For every $\la\in W_0$ and every $\gamma\in {\cal X}$ we set
$$E^\la:=\{\gamma\in{\cal X}\;\colon\; (\gamma,\la)\in E\}\;\;\;\textrm{and}\;\; \;E_{\gamma}:=\{\la\in W_0\;\colon\; (\gamma,\la)\in E\}.$$
We have to show that ${\cal L}eb(E_\gamma)={\cal L}eb(W_0)$ for $\cal M$-almost every $\gamma \in {\cal X}$.
This immediately follows from Tonelli's theorem:
$$\int_{\cal X} {\cal L}eb(E_\gamma)\;d{\cal M}(\gamma)=\int_{W_0} {\cal M}(E^\la)\;d{\cal L}eb(\la) = {\cal L}eb(W_0)$$
since, according to the above third assertion, ${\cal M}(E^\la)=1$ for every $\la\in W_0$.\fin

Our strategy will be to transfer the  estimates known for a fixed system $(J_{\la_0},f_{\la_0},\mu_{\la_0})$ to the system $({\cal X},{\cal F},{\cal M})$.
This will be  possible because the graphs $\Gamma_{\gamma}$ for $\gamma\in {\cal X}$ cannot approach the critical set $C_f$ in a non uniform way, a phenomenon which simply
relies on the compactness of the closure of $\cal X$ and the following basic property (see the Claim in subsection \ref{ddcl}).\\
  
{\bf Fact} 
{\it There exist $0<\alpha\le 1$ such that $\sup_{V_0} \vert \varphi\vert \le \vert \varphi(\la)\vert ^{\alpha}$ for every $\la\in V_0$ and every holomorphic function
$\varphi : W_0\to\C$ such that $0<\vert \varphi\vert<1$.}\\

More specifically, the key uniformity property we need is given by the next  lemma.
In our proofs, we shall denote the smallest singular value of an invertible linear map $L$ of $\C^k$ by $\delta(L)$. 
Let us recall  that $ \delta(L)=\Vert L^{-1}\Vert^{-1}$ and that
$\delta(L) ^k  \leq  \vert \mbox{det} L \vert \leq  \delta(L) \Vert L\Vert^{k-1}$.

\begin{lem}\label{LemKey}
Let $U_0$,$V_0$,$W_0$ be 
as in Lemma \ref{LemPty}. Then there exist $\alpha >0$ and $c>0$ such that $\frac{1}{n} u_n(\gamma,\la) \le \frac{k}{\alpha} \frac{1}{n}u_n(\gamma,\la'_0) +\log  c$ for every $n\geq 1$, every $\gamma \in {\cal X}$ and every $(\la'_0, \la) \in {V_0} \times V_0$.
\end{lem}

\proof By the compactness of $\overline{\cal X}$ and $\overline{V_0}$, we get 
$c_1:=\sup_{\gamma\in {\cal X},\la\in V_0} \Vert DF_{\gamma(\la)}(0) \Vert^{k-1} < +\infty$
and thus $\vert \mbox{det} (DF_{\gamma(\la)}^1(0))\vert \le c_1 \delta  (DF_{\gamma(\la)}(0))$ 
for every $\la\in \overline{V_0}$ and every $\gamma\in {\cal X}$.\\
Then, as $\mbox{det} DF_{\gamma(\la)}^n(0)=\prod _{j=0}^{n-1} \mbox{det} DF_{{\cal F}^j(\gamma)}(0)$  and 
$ \prod _{j=0}^{n-1} \delta( DF_{{\cal F}^j(\gamma)}(0)) \le  \delta ( DF_{\gamma(\la)}^n(0))$ 
we get 
\begin{eqnarray}\label{hh}
\vert \mbox{det} DF_{\gamma(\la)}^n(0)\vert \le c_1^n  \delta ( DF_{\gamma(\la)}^n(0));\;\;\forall \gamma\in{\cal X},\;\forall \la\in\overline{V_0}.
\end{eqnarray}

Let us set $c_2:=\sup_{\la\in \overline{W_0},\gamma\in \overline{\cal X}} \vert \mbox{det} DF_{\gamma(\la)} (0)\vert.$ When $\gamma \in {\cal X}$,  the holomorphic function 
$\varphi(\la):=\frac{1}{c_2^n} \mbox{det} DF_{\gamma(\la)}^n (0)$ is non vanishing and its modulus is bounded by $1$ on $W_0$. Applying the above stated Fact  to $\varphi$, we get $0<\alpha\le 1$ 
(which only depends on $V_0$ and $W_0$) such that:
\begin{eqnarray}\label{h}
\sup_{\la'\in V_0} \vert \mbox{det} DF_{\gamma(\la')}^n(0)\vert \le c_2^{n(1-\alpha)} \vert \mbox{det} DF_{\gamma(\la)}^n(0)\vert^{\alpha};\;\;\forall n\ge 1,\;\forall \gamma\in {\cal X},\;\forall \la\in V_0.
\end{eqnarray}
Using successively (\ref{h}) and (\ref{hh}) we get for any $\la,\la'_0\in V_0$
\begin{eqnarray*}
\left[\delta(DF_{\gamma(\la'_0)}^n(0))\right]^k\le \vert \mbox{det} DF_{\gamma(\la'_0)}^n(0))\vert \le c_2^{n(1-\alpha)} \vert \mbox{det} DF_{\gamma(\la)}^n(0)\vert^\alpha
\le  c_2^{n(1-\alpha)} c_1^{n\alpha} \left[ \delta( DF_{\gamma(\la)}^n(0))\right]^\alpha.
\end{eqnarray*}

Then, applying $\log $ and multiplying by $\frac{-1}{n}$ we get 
$$k\frac{1}{n}u_n(\gamma,\la'_0) \ge \alpha \frac{1}{n} u_n(\gamma,\la) -\alpha (\log  c_1 +\frac{1-\alpha}{\alpha} \log  c_2)$$
which is the desired estimate with 
$c:=c_1 c_2^{(1-\alpha)/\alpha} $.\fin

The next Lemma gathers the properties of $(u_n)_n$ which will be crucial to end the proof.

\begin{lem}\label{PropCru} Let $U_0$,$V_0$,$W_0$ be 
as in Lemma \ref{LemPty}. Then the following properties occur.
\begin{itemize}
\item[1)] The sequence $(\frac{1}{n} u_n)_n$ is uniformly bounded from below on ${\cal X}\times V_0$.
\item[2)] The sequence $\left(\frac{1}{n} u_n(\gamma,\cdot)\right)_n$ is uniformly bounded on $V_0$ for $\cal M$-almost every $\gamma\in {\cal X}$.
\item[3)] The functions ${\widehat u}_n$ are $\widehat{\cal M}$-integrable.
\end{itemize}
\end{lem}

\proof 1) Using the properties of the smallest singular value we have
\begin{eqnarray*}
\frac{1}{n} u_n(\gamma,\la) = - \frac{1}{n} \log  \delta \left(DF_{\gamma(\la)} ^n (0)\right) \ge - \frac{1}{nk} \log  \vert \mbox{det} \left(DF_{\gamma(\la)} ^n (0)\right)\vert\\
= - \frac{1}{k} \left(\frac{1}{n} \sum_{j=0}^{n-1} \log  \vert \mbox{det} DF_{({\cal F}^j \gamma)(\la)} (0)\vert\right)
\end{eqnarray*}
and the assertion follows immediately from the definition and the continuity of $F_{\gamma(\la)}$.\\

2) 
We have just seen that $\frac{1}{n} u_n(\gamma,\cdot)$ is uniformly bounded from below on $V_0$.
By the fourth assertion of Lemma  \ref{LemPty}, for ${\cal M}$-almost every $\gamma\in{\cal X}$ there exists $\la_\gamma\in V_0$ such that $\lim_n \frac{1}{n} u_n(\gamma,\la_\gamma)=-\chi_1(\la_\gamma)$. On the other hand, by Lemma \ref{LemKey}, we have $\frac{1}{n} u_n(\gamma,\la)\le \frac{k}{\alpha} \frac{1}{n} u_n(\gamma,\la_\gamma) +\log  c$ for every $n\in \N$ and every $\la \in V_0$ and thus $\frac{1}{n} u_n(\gamma,\cdot)$ is uniformly bounded from above on $V_0$.\\

3) By the above first assertion, we know that $\widehat{u}_n$ is bounded from below. It thus suffices to show that 
$\int \widehat {u}_n(\widehat{\gamma}) \;d\widehat{\cal M} (\widehat{\gamma}) < +\infty$.  
By Lemma \ref{LemKey} we have
\begin{eqnarray*}
\int \widehat {u}_n(\widehat{\gamma}) \;d\widehat{\cal M} (\widehat{\gamma}) \le n\log  c + \frac{k}{\alpha} \int {u}_n(\pi_0({\widehat \gamma}),\la_0) \;d\widehat{\cal M} (\widehat{\gamma})
= n\log  c + \frac{k}{\alpha} \int {u}_n( \gamma,\la_0) \;d{\cal M} ({\gamma})\\
=n\log  c + \frac{k}{\alpha} \int \log  \Vert (DF_{\gamma(\la_0)}^n (0) )^{-1}\Vert \;d{\cal M} (\gamma)
= n\log  c - \frac{k}{\alpha} \int \log  \delta (DF_{\gamma(\la_0)}^n (0) ) \;d{\cal M} (\gamma).
\end{eqnarray*}
Using (\ref{hh}), we thus get
\begin{eqnarray*}
\int \widehat {u}_n(\widehat{\gamma}) \;d\widehat{\cal M} (\widehat{\gamma}) \le
 -\frac{k}{\alpha} \int \log  \vert \mbox{det} (DF_{\gamma(\la_0)}^n (0) )\vert \;d{\cal M} (\gamma) + \frac{kn}{\alpha} \log  c_1 + n\log  c\\
 =- \frac{k}{\alpha} \int \log   \vert \mbox{det} (Df_{\la_0}^n)_{\gamma(\la_0)}\vert  \;\;d{\cal M} (\gamma)  + C _n
 =- \frac{k}{\alpha} \int \log   \vert \mbox{det} (Df_{\la_0}^n)_x \vert \;(dp_{\la_0 \star}{\cal M})(x) + C_n ,
\end{eqnarray*}
the conclusion follows from the integrability of $\log   \vert \mbox{det} (Df_{\la_0}^n)_x\vert$ for  $p_{\la_0\star} {\cal M}=\mu_{\la_0}$, see \cite{DS2}. \fin

We are now ready to establish the main result of this subsection.\\
 
\textsc{Proof of Proposition \ref{ThmLya}:} 
We will apply Kingman subadditive ergodic theorem (see \cite{Arn}) to the sequence $(\widehat{u}_n)_n$. This is possible since the system $(\widehat {\cal X},\widehat {\cal F},\widehat {\cal M})$ is ergodic, the sequence $(\widehat{u}_n)_n$ is subadditive (second assertion of Lemma \ref{LemPty}) and $\widehat{u}_1\in L^1(\widehat{\cal M})$ (last assertion of Lemma \ref{PropCru}). According to this theorem, there exists $L\in \R$ such that $\lim_n \frac{1}{n}\widehat{u}_n(\widehat{\gamma})=L$ for $\widehat{\cal M}$-almost every
 $\widehat\gamma\in \widehat{\cal X}$
and $\lim_n \frac{1}{n}\int_{\widehat{\cal X}} \widehat{u_n}\;d\widehat{\cal M} =L$. It remains to show that $L\le - {1 \over 2} \log d$. \\

Taking into account the fourth assertion of Lemma \ref{LemPty} and the second assertion of Lemma \ref{PropCru}, we may  thus pick $\widehat{\gamma} \in \widehat{\chi}$ such that:

\begin{itemize}
\item[i)] $\lim_n \frac{1}{n}\widehat{u}_n(\widehat{\gamma})=L$,
\item[ii)] $\frac{1}{n} u_n(\gamma_0,\cdot)$ is uniformly bounded on $V_0$,
\item[iii)] $\lim_n \frac{1}{n} u_n(\gamma_0,\la)=-\chi_1(\la)$ for Lebesgue-almost every $\la\in V_0$.
\end{itemize}

 Assuming that $L > \frac{-\log  d}{2}$, we will reach a contradiction with the fact that $\chi_1(\la)\ge \frac{\log  d}{2}$ for all $\la$ (see \cite{BD1} or \cite{DS2}).
 Recalling that $\widehat{u}_n(\widehat{\gamma})= \sup_{\la \in U_0} u_n(\gamma_0,\la)$, there exist $\la_{n_k}\in U_0$ and $\epsilon >0$ such that $\la_{n_k}\to \la'_0\in \overline{U_0}$ and $\frac{1}{n_k} u_{n_k} (\gamma_0,\la_{n_k}) \ge \frac{-\log  d}{2} + \epsilon$.
 We may  pick $r>0$ such that $B(\la_{n_k},r)\subset V_0$ for all $k\in \N$. Then, by the subharmonicity of $u_{n_k}(\gamma_0,\cdot)$ on $V_0$ (first assertion of Lemma \ref{LemPty}) we get:
 $$\frac{-\log  d}{2} + \epsilon \le \frac{u_{n_k}(\gamma_0,\la_{n_k})}{n_k}\le \frac{1}{\vert B(\la_{n_k},r)\vert} \int_{B(\la_{n_k},r)} \frac{u_{n_k}(\gamma_0,\la)}{n_k}$$
 which, by Lebesgue dominated convergence theorem, yields $$\frac{-\log  d}{2} + \epsilon \le \frac{1}{\vert B(\la'_{0},r)\vert} \int_{B(\la'_{0},r)}  -\chi_1(\la)$$
 and contradicts the fact that $\chi_1(\la)\ge \frac{\log  d}{2}$ for every $\la$. \finsec

 %%%%%%%%%%%%%%%
 \subsection{Proof of Theorem \ref{PropExHM}}\label{SecProof}
%%%%%%%%%%%%%%%

 According to  Corollary \ref{CorASW}, we only need to consider the case where  $f$ admits an acritical and ergodic
 equilibrium  web ${\cal M}_0$. Let ${\cal K}_0:=\supp {\cal M}_0$. The proof is based on the following key property.\\
 
 \noindent{\bf Fact}:
${\cal M}_0\left(\{\gamma\in {\cal K}_0\;\colon\; \exists k\in \N, \exists \gamma'\in {\cal K}'_0\;\textrm{s.t.}\;\Gamma_{{\cal F}^k(\gamma)} \cap \Gamma_{\gamma'} \ne \emptyset\;
\textrm{and}\; {\cal F}^k(\gamma)\ne \gamma'\}\right)=0$ for every compact subset ${\cal K}'_0$ of $\cal J$.\\

We shall construct the lamination by applying this fact with ${\cal K}'_0={\cal K}_0$, the uniqueness assertions will be obtained by applying it with
${\cal K}'_0=\supp{\cal M'}_0$ for any other equilibrium web ${\cal M'}_0$ .\\

To prove the Fact, it
is sufficient  to show that for any fixed $k\in \N$  and any $\la_0 \in M$ there exists a neighbourhood $U_1$ of $\la_0$ such that
\begin{eqnarray}\label{firstgoal}
{\cal M}_0\left(\{\gamma\in {\cal K}_0\;\colon\; \exists \gamma'\in {\cal K}'_0\;\textrm{s.t.}\;\Gamma_{{\cal F}^k(\gamma)} \cap \Gamma_{\gamma'}\cap \left(U_1\times\Pj^k\right) \ne \emptyset\;
\textrm{and}\; {{\cal F}^k(\gamma)} \ne \gamma'\}\right)=0.
\end{eqnarray}

 To this purpose, we shall work with the natural extension $\left(\widehat{\cal X},\widehat{\cal F},
\widehat{{\cal M}_0}\right)$ of the system $\left({\cal X}, {\cal F}, {\cal M}_0\right)$ and apply
Proposition \ref{LemSF}. We recall that, according to Proposition \ref{ThmLya}, all the assumptions of Proposition \ref{LemSF} are satisfied.
Let $U_0$ be a neighbourhood of $\la_0$,  we may assume that $U_0$ is simply connected and that
$U_1\Subset U_0 \Subset M$. 
Let $p$ be the integer and $\widehat{\eta}_p : \widehat{\cal Y}\to ]0,1]$ be the measurable function defined on the full 
$\widehat{{\cal M}_0}$-measure
set $\widehat{\cal Y}$ given by Proposition \ref{LemSF}. 
We recall that ${\cal X}= {\cal K}_0\setminus {\cal J}_s$ and ${\cal M}_0({\cal J}_s)=0$.\\

For any $B\subset U_0$, we define the \emph{ramification} function $R_B$ by setting
\begin{eqnarray*}
R_B(\gamma):=\sup_{\gamma'\in  {\mathcal K}'_0 : \Gamma_{\gamma'\vert B} \cap \Gamma_{\gamma\vert B} \neq \emptyset}
\sup_B d_{\Pj^k}\left(\gamma(\la),\gamma'(\la)\right), \;\;
 \forall \gamma\in {\mathcal J}.
\end{eqnarray*}
Let $\widehat{\cal Y}_\epsilon:=\{\widehat{\gamma}\in{\widehat {\cal Y}}\;\colon\; R_{U_0}(\gamma_k)>\epsilon\}$,
it then  suffices to prove that $\widehat{{\cal M}_0}\left( \widehat{\cal Y}_\epsilon\right)=0$ for every $\epsilon >0$ as it follows from the following observation:
 \begin{eqnarray*}
 {\cal M}_0\left(\{\gamma\in {\cal K}_0\;\colon\; \exists \gamma'\in {\cal K}'_0\;\textrm{s.t.}\;\Gamma_{{\cal F}^k(\gamma)} \cap \Gamma_{\gamma'}\cap\left(U_0\times\Pj^k\right) \ne \emptyset\;
\textrm{and}\; {\cal F}^k(\gamma)\ne \gamma'\}\right)\\
={\cal M}_0\left(\{\gamma\in{\cal K}_0\;\colon\; R_{U_0}({\cal F}^k(\gamma))>0\}\right)
 = {\cal M}_0\left(\{\gamma\in{\cal X}\;\colon\; R_{U_0}({\cal F}^k(\gamma))>0\}\right)
\\
=\widehat{{\cal M}_0}\left(\{\widehat{\gamma}\in{\widehat {\cal Y}}\;\colon\; R_{U_0}(\gamma_k)>0\}\right)=
\widehat{{\cal M}_0}\left(\cup_{s\in \N^*} \widehat{\cal Y}_{\frac{1}{s}}\right).
\end{eqnarray*}
Let us proceed by contradiction and assume that $\widehat{{\cal M}_0}\left( \widehat{\cal Y}_\epsilon\right) >0$ for some $\epsilon>0$.
Then, after reducing $\epsilon >0$, we may assume that 
$\widehat{{\cal M}_0}\left(\{ \widehat \gamma \in \widehat{\cal Y}_\epsilon\;\colon\;\widehat{\eta}_p(\widehat {\cal F}^k(\widehat {\gamma})) >\epsilon \}\right) >0$.
In the sequel we shall denote $\widehat{\gamma}_k:= \widehat {\cal F}^k(\widehat {\gamma})$.
Owing to the equicontinuity of ${\cal X}\cup {\cal K}'_0$ (we recall that ${\cal X} \subset {\cal K}_0$)  we may cover $U_1$ with  finitely many open sets $B_i \subset U_0$, say with $1\le i\le N$, such that
\begin{eqnarray}\label{equic}
\forall (\gamma,\gamma')\in {\cal X}\times {\cal K}'_0,\;\forall \la_1\in B_i\;:\; \gamma(\la_1)=\gamma'(\la_1) \Rightarrow \sup_{\la\in B_i} d_{\Pj^k}\left(\gamma(\la),\gamma'(\la)\right) < \epsilon.
\end{eqnarray}
As $R_{U_1}(\gamma)=0$ when $\max_{1\le i\le N} R_{B_i}(\gamma) =0$ (by analyticity we have $\gamma=\gamma'$ on $U_1$ if 
$\gamma=\gamma'$ on some $B_i$), there exist $1\le j\le N$ and $\alpha>0$ such that:
\begin{eqnarray*}
\widehat{{\cal M}_0}\left(\{ {\widehat \gamma}\in \widehat{\cal Y}_\epsilon\;\colon\;  \widehat{\eta}_p({\widehat \gamma}_k) >\epsilon\;\textrm{and}\;R_{B_j}(\gamma_k)>\alpha\}\right) >0.
\end{eqnarray*}
Let us set $\widehat{\cal Y}_{\epsilon,j,\alpha}:=\{\widehat \gamma\in \widehat{\cal Y}_\epsilon\;\colon\; \widehat{\eta}_p({\widehat \gamma}_k) >\epsilon\;\textrm{and}\; R_{B_j}(\gamma_k)>\alpha\}.$ Applying Poincar\'e recurrence theorem to $\widehat{\cal F}^{-p}$, we  find $\widehat\gamma\in \widehat{\cal Y}_{\epsilon,j,\alpha}$ and an increasing sequence of integers $(n_q)_q$ with $n_q\in p\N$ such that $\widehat{\cal F}^{-n_q} (\widehat \gamma) \in \widehat{\cal Y}_{\epsilon,j,\alpha}$
for every $q\in \N$. In particular $\widehat\gamma\in \widehat{\cal Y}_{\epsilon,j,\alpha}$ and $R_{B_j} (\gamma_{k-n_q})>\alpha$ for every $q\in \N$.
We will reach a contradiction by establishing that
\begin{eqnarray}\label{Lim0}
\lim_{m\to+\infty} R_{B_j} (\gamma_{k-m p})=0, \;\;\forall \widehat\gamma\in \widehat{\cal Y}_{\epsilon,j,\alpha}.
\end{eqnarray}
To this purpose we shall use Proposition \ref{LemSF} to show that $ R_{B_j} (\gamma_{k-n}) \le e^{-nA}$ when $n\in p\N$ and $\widehat\gamma\in \widehat{\cal Y}_{\epsilon,j,\alpha}$. Let $\gamma'\in {\cal K}'_0$ such that $\gamma'(\la_1)=\gamma_{k-n} (\la_1)$ for some $\la_1\in B_j$. Then
$({\cal F}^n\gamma')(\la_1)=\gamma_k(\la_1)$ and thus, according to (\ref{equic}), 
$\sup_{\la\in B_j} d\left(({\cal F}^n\gamma')(\la),\gamma_k(\la)\right) < \epsilon < \widehat{\eta}_p({\widehat \gamma}_k)$. This means that
\begin{eqnarray}\label{tube}
\Gamma_{{\cal F}^n\gamma'} \cap \left(B_j\times \Pj^k\right)\subset T_{B_j}\left(\gamma_k,\widehat{\eta}_p({\widehat \gamma}_k)\right).
\end{eqnarray}
Now, by Proposition \ref{LemSF}, the inverse branch $f_{{\widehat \gamma}_k}^{-n}$ of $f^n$ is defined on the tube 
$T_{U_0}\left(\gamma_k,\widehat{\eta}_p({\widehat \gamma}_k)\right)$ and maps it biholomorphically into $T_{U_0}\left(\gamma_{k-n},e^{-nA}\right)$.
As $B_j \subset U_0$, this yields:
\begin{eqnarray}\label{tubeagain}
f_{{\widehat \gamma}_k}^{-n}\left(
T_{B_j}\left(\gamma_k,\widehat{\eta}_p({\widehat \gamma}_k)\right) \right) \subset T_{B_j}\left(\gamma_{k-n},e^{-nA}\right).
\end{eqnarray}
By construction,  we have $f_{{\widehat \gamma}_k}^{-n}\left(
\Gamma_{\gamma_k}\right) =\Gamma_{\gamma_{k-n}}$ and therefore 
$f_{{\widehat \gamma}_k}^{-n}\left(({\cal F}^n\gamma')(\la_1)\right) =f_{{\widehat \gamma}_k}^{-n}\left(\gamma_k(\la_1)\right) =\gamma_{k-n}(\la_1)=\gamma'(\la_1)$. This implies that $f_{{\widehat \gamma}_k}^{-n}\left(\Gamma_{{\cal F}^n\gamma'}\right)=\Gamma_{\gamma'}$ which in turns, by 
(\ref{tube}) and (\ref{tubeagain}), implies that $\sup_{\la\in B_j} d_{\Pj^k}(\gamma'(\la),\gamma_{k-n}(\la)) \le e^{-nA}$. Then (\ref{Lim0}) follows and thus (\ref{firstgoal}) and the Fact are proved.\\

Let us now establish the existence of an equilibrium lamination ${\cal L}_0$. Consider the set
\[{\cal L}_0^+:=\{\gamma\in {\cal K}_0\setminus {\cal J}_s\;\colon\; \forall \gamma'\in {\cal K}_0, \forall k \in \N,  
\Gamma_{{\cal F}^k(\gamma)}\cap \Gamma_{\gamma'} \ne \emptyset \Rightarrow {\cal F}^k(\gamma)=\gamma'\}.\]

The Fact, applied with ${\cal K}'_0={\cal K}_0$ yields  ${\cal M}_0\left({\cal L}_0^+\right)=1$ and, by construction,  ${\cal L}_0^+$ satisfies the following 
properties:
\begin{itemize}
\item[1)] ${\cal L}_0^+\subset {\cal K}_0\setminus {\cal J}_s\subset{\cal J}\setminus {\cal J}_s$,
\item[2)] ${\cal F}\left({\cal L}_0^+\right)\subset{\cal L}_0^+$,
\item[3)] $ \forall \gamma, \gamma'\in {\cal L}_0^+\; : \; \Gamma_\gamma \cap \Gamma_{\gamma'} \ne \emptyset \Rightarrow \gamma=\gamma'$.
\end{itemize}

The set ${\cal L}_0:=\cup_{m\ge 0} {\cal F}^{-m} \left({\cal L}_+\right)$ also satisfies the  properties (2) and (3). It is also relatively compact in $\cal J$, indeed
${\cal L}_0^+\subset {\cal K}_0$ and taking inverse branches cannot destroy the equicontinuity (a result of Ueda \cite[Theorem 2.1]{U}).
Moreover ${\cal F} : {\cal L}_0\to {\cal L}_0$ is $d^k$-to-1.\\

We finally prove the uniqueness assertions. Let ${\cal M}'_0$ be an equilibrium web for $f$
(or, more generally, a compactly supported probability measure on ${\cal J}$ such that
$p_{\la_0\star} {\cal M}'_0=\mu_{\la_0}$ for some $\la_0\in M$). We set ${\cal K}'_0:=\supp
{\cal M}'_0$. Let us fix $\la_0\in M$ and recall that $\mu_{\la_0}=p_{\la_0\star} {\cal M}'_0=
p_{\la_0\star} {\cal M}_0$. Then, for any Borel subset $\cal A$ of $\cal J$, we have $\mu_{\la_0}\left(p_{\la_0}({\cal K}'_0\cap {\cal A})\right)= {\cal M}'_0 \left(p_{\la_0}^{-1}\left(p_{\la_0}({\cal K}'_0\cap {\cal A})\right)\right)\ge  {\cal M}'_0 \left({\cal K}'_0\cap {\cal A}\right)=
{\cal M}'_0 \left({\cal A}\right)$ and thus
\begin{eqnarray*}
{\cal M}_0 \left( \{\gamma\in {\cal K}_0\;\colon\; \exists \gamma'\in {\cal K}'_0\cap {\cal A}\;\textrm{s.t.}\;
\gamma(\la_0)=\gamma'(\la_0)\}\right) ={\cal M}_0 \left(p_{\la_0}^{-1}\left(p_{\la_0}({\cal K}'_0\cap {\cal A})\right)\right)\\
=\mu_{\la_0}\left(p_{\la_0}({\cal K}'_0\cap {\cal A})\right)\ge 
{\cal M}'_0 \left({\cal A}\right).
\end{eqnarray*} 
But, according to  the Fact we have
$${\cal M}_0 \left( \{\gamma\in {\cal K}_0\;\colon\; \exists \gamma'\in {\cal K}'_0\cap {\cal A}\;\textrm{s.t.}\;
\gamma(\la_0)=\gamma'(\la_0)\}\right) ={\cal M}_0\left({\cal K}'_0\cap{\cal K}_0\cap {\cal A}\right)={\cal M}_0\left({\cal K}'_0\cap {\cal A}\right)$$ and therefore 
${\cal M}_0\left({\cal K}'_0\cap {\cal A}\right)\ge {\cal M}'_0 \left({\cal A}\right)$.
This implies that ${\cal M}_0\left({\cal K}'_0\right)=1$ and that ${\cal M}_0\ge{\cal M}'_0$. As both ${\cal M}_0$ and ${\cal M}'_0$ are probability measures, we have proved that 
${\cal M}_0={\cal M}'_0$.\\

Let ${\cal L}'$ be an arbitrary equilibrium lamination for $f$. Let us pick $\la_0\in M$ and set ${\cal L}'_{\la_0}:= p_{\la_0}\left({\cal L}'\right)$. Using ${\cal M}_0\left({\cal L}_0\right)=1$, $\mu_{\la_0}=p_{\la_0\star} {\cal M}_0$ and $\mu_{\la_0}\left({\cal L}'_{\la_0}\right)=1$ yields  ${\cal M}_0\left(\{\gamma\in {\cal L}_0\;\colon\; \gamma(\la_0)\in {\cal L}'_{\la_0}\}\right)=1$. On the other hand,  the Fact implies that ${\cal M}_0\left(\{\gamma\in {\cal L}_0\;\colon\; \gamma(\la_0)\in {\cal L}'_{\la_0}\}\right)={\cal M}_0\left({\cal L}_0\cap{\cal L}'\right)$. This shows that ${\cal M}_0\left({\cal L}_0\Delta{\cal L}'\right)=0$.
 \finsec

%%%%%%%%%%%%%%%%%%%%%
\section{ Siegel discs and bifurcations} \label{ContJ}
%%%%%%%%%%%%%%%%%%%%%

As it is well known,  the Julia sets of any holomorphic family of rational maps of $\Pj^1$  depends continuously on the parameter for the Hausdorff topology if and only if the family is stable. It is worth emphasizing that discontinuities can be explained by the appearance of Siegel discs, see \cite{Do}. We investigate this in higher dimension and,
as a consequence, show that the existence of virtually repelling Siegel periodic points in the Julia set (see Definitions \ref{defiSieg} and \ref{defiVR}) is an obstruction to the existence of an equilibrium  web. We finally exploit this fact to end the proof of  Theorem \ref{main2}.

%%%%%%%%%%%%%%%%%%%%%%%%%%%%%
 \subsection{Siegel discs as obstructions to stability}
%%%%%%%%%%%%%%%%%%%%%%%%%%%%%%
 
 We define a notion of Siegel disc for endomorphisms of $\Pj^k$ and investigate how they behave with respect to Julia sets.
In this subsection, we endow $\C^k$ with the norm $\Vert z\Vert :=\sup_i \vert z_i\vert$ and set $1 \leq q \leq k-1$. We write $z=:(z',z'')$
where $z':=(z_1, \cdot\cdot\cdot, z_{k-q})\in \C^{k-q}$ and $z'':=(z_{k-q+1}, \cdot\cdot\cdot, z_k)\in \C^q$. We also set $k':=k-q$, $e^{i\theta_0}:=(e^{i\theta_{0,k'+1}},\cdots e^{i\theta_{0,k}})$ and $e^{i\theta_0}\cdot z'':=(e^{i\theta_{0,k'+1}} z_{k'+1},\cdots, e^{i\theta_{0,k}} z_{k})$.

\begin{defn}\label{defiSieg}
Let $f_0$ be  a holomorphic endomorphism of $\Pj^k$. One says that $z_0\in \Pj^k$ is a \emph{Siegel fixed point} for $f_0$  if 
$f_0$ is holomorphically linearizable at $z_0$ and its differential at $z_0$ is of the form $\left(A_0z',e^{i\theta_0}\cdot z''\right)$
where $A_0$ is an expanding linear map on $\C^{k'}$ and $\pi, \theta_{0,k'+1}, \cdots,\theta_{0,k}$ are linearly independent over $\Q$.
In other words, there
exists a local holomorphic chart $\psi_0: B_R\to \Pj^k$ such that $\psi_0(0)=z_0$ and 
\begin{center}
$\psi_0^{-1} \circ f_0\circ \psi_0=\left(A_0z',e^{i\theta_0}\cdot z''\right)$
\end{center}
where $\theta_0$  and $A_0$ are as above.
Any set of the form $\psi_0\left( \{0'\}\times B_\rho\right)$ where $\rho < R$ and $B_\rho$ is a ball centered at the origin in $\C^q$ is called a \emph{local Siegel $q$-disc of $f_0$ centered at $z_0$}.
\end{defn}

Let us consider a holomorphic family $f$ of endomorphisms of $\Pj^k$. If $f_{0}$ admits a Siegel fixed point $z_0$ then, by the implicit function theorem, there exists a unique holomorphic map $z(\la)$ defined on some neighbourhood of $0$ in $M$ such that $z(0)=z_0$ and $z(\la)$ is fixed by $f_\la$. Moreover, since $\theta_{0,k'+1},\cdots,\theta_{0,k}$ are pairwise distinct, there exist  holomorphic functions $w_{j}(\la)$ such that $w_j(0)=e^{i\theta_{0,j}}$  and $w_j(\la)$ is an eigenvalue of $d_{z(\la)}f_\la$ for $k'+1\le j\le k$. In this context, we coin the following definition.  

\begin{defn}\label{defiVR}
The Siegel fixed point $z_0$ is called  \emph{virtually repelling} if there exist a holomorphic disc $\sigma : \Delta_{\epsilon_0} \to M$ and positive constants $c_j$ such that $\sigma(0)=0$
and $\vert w_j\circ \sigma (t) \vert =1+ c_j t+ o(t)$ for $k'+1\le j\le k$ and $-t_0 < t< t_0$. If, moreover, $z \circ \sigma(t) \in J_{\sigma(t)}$ for $-t_0 < t< t_0$
the Siegel fixed point $z_0$ is called  \emph{virtually $J$-repelling}.
\end{defn}

Let us observe that if $J_\la$ is continuous at $\la_0$ and if $f_{\la_0}$ has a virtually repelling Siegel periodic point outside $J_{\la_0}$, then $\la_0$ must be accumulated by parameters $\la$ for which $f_{\la}$ has periodic repelling points
 outside $J_{\la}$. Examples of such repelling points have been given by Hubbard-Papadopol \cite[section 6, example 2]{HP} (see also Fornaess-Sibony \cite[section 4.1]{FS2}). The following proposition discusses the position of Siegel discs with respect to Julia sets. Note that the second item will only be used in Remark \ref{RemContSieg}.

\begin{prop}\label{SiegJul}
Let $f:M\times \Pj^{k}\to M\times \Pj^{k}$ be a holomorphic family of endomorphisms of $\Pj^k$ such that $f_{\la_0}$ admits a virtually repelling Siegel fixed point $z_0$.
\begin{itemize}
\item[1)] If $f$ admits an equilibrum  web then every local Siegel $q$-disc centered at $z_0$  is contained in $\Pj^k\setminus J_{\la_0}$. In particular $z_0\notin J_{\la_0}$.
\item[2)] When $q=1$, if $z_0\in J_{\la_0}$ and if  $\la \mapsto J_\la$ is continuous at $\la_0$ then  any local Siegel $q$-disc centered at $z_0$ is contained in $J_{\la_0}$.
\end{itemize}
\end{prop}

The first item of the preceding proposition immediately yields the following result.

\begin{cor}\label{corSiegBif}
Let $f:M\times \Pj^{k}\to M\times \Pj^{k}$ be a holomorphic family of endomorphisms of $\Pj^k$. Let $U_0$ be any neighbourhood of $\la_0$ in $M$.
If the restriction of $f$ to $U_0\times\Pj^k$ admits an equilibrium web then $f_{\la_0}$ has no virtually repelling Siegel periodic point in $J_{\la_0}$.
\end{cor}

The proofs of items 1) and 2) of Proposition \ref{SiegJul} respectively use items 1) and 2) of the following lemma. We shall also need the fact that $\mu_\la$ gives no mass to analytic sets.

\begin{lem}\label{techlem}
Let $g :\Delta_{r_0}\times B_R\to \Delta_{r_0}\times B_{R'}$ be a holomorphic map such that
$g(\la,z)= \left(\la, g_\la(z)\right)$,  $g_\la(0)=0$ and  
$g_0(z)=\left(A_0 ^{-1} \cdot z',e^{- i\theta_0}\cdot z''\right)$ where $A_0$ is an expanding linear map on $\C^{k'}$.
Assume that $\frac{\partial g_{\la,j}}{\partial z_i} (0)=0$ for $k'+1\le j\le k$ and $i\ne j$. Assume moreover that
there exists $\vert u_0 \vert =1$ and $c_{j}>0$ such that $\vert \frac{\partial g_{ t u_0 ,j}}{\partial z_j} (0) \vert =1+ c_j t + o(t)$ for $k'+1\le j\le k$ and $-r_0 < t< r_0$. Then, after taking $r_0$ and $R$ smaller, the following properties occur.
\begin{itemize}
\item[1)] There exists arbitrarily small $\la$ such that $\Vert g_\la(z)\Vert \le \alpha_0 \Vert z\Vert$ on $B_R$ with $0<\alpha_0<1$.
\item[2)] Assume $k'=k-1$. For any $0<\rho<R_1<R_2<R$, there exists arbitrarily small $\la$ such that, for every $a\in B_{R_1}$ which does not belong to the local stable manifold $S_\la$ of $g_\la$, there exists $n_0$ such that 
$g_\la^{n_0} (a) \in \{\Vert z'\Vert < \rho\} \times \{R_1 < \Vert z''\Vert < R_2\}$ and $g_{\la}^k(a)\in B_{R_1}$ for $0\le k\le n_0-1$.
\end{itemize}
\end{lem}

\proof We will  exploit the form of  the Taylor expansion of order one of $g$.  Let us write $g_\la:=\left(g_{\la,j}\right)_{1\le j\le k}$ as
\begin{eqnarray*}
&g_{\la,j}&=\sum_{i=1}^{k'} \left(a_{ij} +\la \mu_{ij}(\la)+\la q_{ij}(\la,z)\right)z_i + \la\sum_{i=k'+1}^k s_{ij}(\la,z)z_i \ \ \textrm{for}\; 
1\le j\le k'\\
&g_{\la,j}&= \left(e^{i\theta_j} +\la \mu_{jj}(\la)+\la q_{jj}(\la,z)\right)z_j+ \la \sum_{i\ne j}  s_{ij}(\la,z)z_i \ \ \textrm{for}\; 
k'+1\le j\le k
\end{eqnarray*}
where $\mu_{ij}$, $q_{ij}$ and $s_{ij}$ are holomorphic on $\Delta_{r_0}\times B_R$ and satisfy $q_{ij}(\la,0)=q_{jj}(\la,0)=0$. By assumption, we also have $s_{ij}(\la,0)=0$ for  $k'+1\le j\le k$ and $i\ne j$.\\
By shrinking $r_0$ and $R$, there exists $0<\alpha_1<1$ such that 
\begin{eqnarray}\label{hlaj}
\sup_{1\le j\le k'} \vert g_{\la,j} (z)\vert \le \alpha_1\Vert z\Vert\;\textrm{on}\; \Delta_{r_0}\times B_R.
\end{eqnarray}

Let us set $\la_t:= t u_0$ where $-r_0<t<r_0$ and
$Q_{jt} (z):= e^{i\theta_j} + \la_t \mu_{jj}(\la_t) +  \la_t q_{jj}(\la_t,z)$ and $R_{jt}(z) :=\vert \la_t\vert \sum_{i\ne j} \vert s_{ij}(\la_t,z)\vert$
 for $k'+1\le j\le k$. Then, by our assumptions and after taking $r_0$ and $R$ smaller, we have 
\begin{eqnarray}\label{Atz}
\vert Q_{jt}(z)\vert \le 1+\frac{c_jt}{2}\;\textrm{for}\; -r_0<t<0\;\textrm{and}\; z\in B_R
\end{eqnarray}
\begin{eqnarray}\label{Rtz}
 R_{jt}(z) \le \frac{c_j\vert t\vert}{4}\;\textrm{for}\; -r_0<t<r_0\;\textrm{and}\; z\in B_R
\end{eqnarray}
\begin{eqnarray}\label{Atz2}
1+\frac{c_jt}{2} \le \vert Q_{jt}(z)\vert \le 1+2c_jt\;\textrm{for}\;  0<t<r_0\;\textrm{and}\; z\in B_R.
\end{eqnarray}
It follows from (\ref{Atz}) and (\ref{Rtz}) that $\vert g_{\la_t,j}(z)\vert \le (1+\frac{t c_j}{4}) \Vert z\Vert$ for $k'+1\le j\le k$, $-r_0<t<0$ and $z\in B_R$.
This and (\ref{hlaj}) yields the first assertion of the lemma.\\

Let us now establish the second one. Fix $0<t<t_0$ so small that $(1+\frac{9t c_k}{4}) R_1<R_2$. Let $a \in B_{R_1}$ be outside the
local stable manifold of $g_{\la_t}$. Assume that one cannot find $n_0$ such that $g_{\la_t}^p (a)\in B_{R_1}$ for $0\le p\le n_0-1$ and
$g_{\la_t}^{n_0}(a) \in \{\Vert z'\Vert <\rho\}\times \{\Vert z''\Vert >R_1\}$.
Then, according to (\ref{hlaj}),  the sequence $a_n:=g_{\la_t}^n (a)$ is well defined and $\Vert a'_n\Vert \to 0$. From (\ref{Rtz}) and (\ref{Atz2}) one gets $\vert a_{n+1,k}\vert \ge (1+\frac{c_k t}{2})\vert a_{n,k}\vert -\frac{t c_k}{4}\Vert a'_n\Vert$ and therefore, since  $\Vert a'_n\Vert\to 0$ and $(a_{n,k})_n$ is bounded,  $\limsup_n \vert a_{n,k}\vert =0$. Then $a_n$ tends to the origin and  this contradicts   the fact that $a$ does not belong to the local stable manifold of $g_{\la_t}$. Thus $n_0$ exists and it remains to check that 
$\Vert a_{n_0}''\Vert <R_2$. From (\ref{Rtz}) and (\ref{Atz2}) one gets 
$\vert a_{n_0,k}\vert \le (1+2c_kt)\vert a_{n_0-1,k}\vert + \frac{t c_k}{4} \Vert a'_{n_0-1}\Vert \le (1+\frac{9t c_k}{4})R_1< R_2$.\fin

\textsc{Proof of Proposition \ref{SiegJul}:} 
We may assume that $M=\Delta_{\epsilon_0}$ and $\la_0=0$ so that $z_0$ is a virtually repelling Siegel fixed point of $f_0$. Thus there exists a biholomorphism $\psi_0 : B_R\to \psi_0\left(B_R\right)$ such that $\psi_0(0)=z_0$
and $\psi_0^{-1}\circ f_0\circ \psi_0=\left(A_0\cdot z',e^{i\theta_0}\cdot z''\right)$ where  $A_0$ is linear and expanding on $\C^{k'}$ and
$\pi, \theta_{0,k'+1}, \cdots,\theta_{0,k}$ are linearly independent over $\Q$.
Let $z(\la)$ be the fixed point of $f_\la$ obtained by the implicit function theorem.

The mapping $\psi_0^{-1}\circ f_\la^{-1}\circ \psi_0$ is well defined on $\Delta_{\epsilon_0}\times B_R$ after taking $R$ and $\epsilon_0$ smaller.
Since the $e^{i\theta_{0,j}}$ are pairwise distinct for $k'+1\le j\le k$, we may find $q$ linearly independent vectors $v_{k'+1}(\la),\cdots, v_k(\la)$
in $\C^k$ and $q$ scalars $w_{k'+1}(\la),\cdots, w_k(\la)$ which depend holomorphically on $\la \in \Delta_{\epsilon_0}$ and such that
\begin{eqnarray}\label{reduction}
d_{\psi_0^{-1} (z(\lambda))} \left(\psi_0^{-1}\circ f_\la^{-1}\circ \psi_0\right) (v_j(\la) )= w_j(\la) v_j(\la)\;\;\textrm{for}\; k'+1\le j\le k.
\end{eqnarray}

Using basis like $\left(v_1,\cdots,v_{k'}, v_{k'+1}(\la),\cdots, v_k(\la)\right)$ we may perform change of coordinates of the form $\left(\la,A(\la,z)\right)$ where $A(\la,\cdot)$ is affine on $\C^k$ which, conjugate by $\psi_0$, yield
 biholomorphisms $\psi_\la : B_R\to \psi_\la\left(B_R\right)$ such  that $g_\la:=\psi_\la^{-1}\circ f_\la^{-1}\circ \psi_\la$
satisfies the assumptions of Lemma \ref{techlem}. The condition $\frac{\partial g_{\la,j}}{\partial z_i} (0)=0$ for $k'+1\le j\le k$ and $i\ne j$  indeed follows from (\ref{reduction}) and the condition $\vert \frac{\partial g_{ t u_0 ,j}}{\partial z_j} (0) \vert =1+ c_j t + o(t)$ follows from the fact that $z_0$ is virtually repelling. By an  abuse of notation we shall assume that $f_\la^{-1}=g_\la$.
%, in particular denote $J_\la$ the set $\psi_\la^{-1} \left(J_\la\cap \psi_\la(B_R)\right)$.\\

1)  We proceed by contradiction and assume that $(0',z_0'')\in J_0$ for $0< \Vert z_0''\Vert <r<R$. Since $f$ admits an equilibrium web,  Lemma \ref{lemMCR} ensures that  there exists a holomorphic map $\gamma : \Delta_{\epsilon_0} \to \Pj^k$ such that $ \gamma (0)=(0' , z_0'')$ and $(\FF^n \cdot \gamma)_n$ is normal on $\Delta_{\epsilon_0}$. 
%We may assume that $\tilde{\gamma}(\la):=\psi_\la^{-1} \left( \gamma (\la)\right)$ is well defined on 
%$\Delta_{\epsilon_0}$. 
Since $ f_0^n \left({\gamma}(0)\right)=(0', e^{i n\theta_0}\cdot z_0'')$ and $(\FF^n \cdot \gamma)_n$ is normal,
after reducing $\epsilon_0$, we may suppose that 
\begin{eqnarray}\label{hta}
\Vert  f_\la^n \left({\gamma}(\la)\right)\Vert \le r\;\textrm{on}\;\Delta_{\epsilon_0}\;\textrm{for}\;n\ge 1.
\end{eqnarray}

Let us recall that $g_\la = f_\la^{-1}$. By Lemma \ref{techlem}, there exists $\la_p\to 0$ and $0<\alpha_p<1$ such that $\Vert g_{\la_p}(z)\Vert \le \alpha_p \Vert z\Vert$ on $B_R$. We may thus find a sequence
$n_p\to \infty$ such that
\begin{eqnarray}\label{trz}
\Vert g_{\la_p}^{n_p}(z)\Vert \le \frac{1}{p}\Vert z\Vert\;\textrm{on}\; B_r.
\end{eqnarray}
From (\ref{hta}) and (\ref{trz}) one gets 
\begin{eqnarray}
\Vert {\gamma} (\la_p)\Vert =\Vert g_{\la_p}^{n_p} \circ f_{\la_p}^{n_p} \left({\gamma}(\la_p)\right)\Vert \le \frac{r}{p}
\end{eqnarray}
which is impossible since $\lim_p \Vert {\gamma} (\la_p)\Vert =\Vert {\gamma}(0)\Vert=\Vert z_0''\Vert >0$.\\
So far we have shown that the punctured $q$-disc $\{0'\}\times \{0<\Vert z''\Vert <R\}$ is contained in $J_0^c$. Since $J_0$ is totally invariant and 
$g_0=\left(A_0^{-1}\cdot z',e^{-i\theta_0}\cdot z''\right)$ where  $A_0$ is linear and expanding, this implies that $B_R\setminus \{z\in B_R\;\colon\; z''=0\} \subset J_0^c$. Finally, as $\mu_0$ does not give mass to analytic sets, we get $B_R\subset J_0^c$.
 
2) We have to show that $(0',z_{0k})\in  J_0$ if $0<\vert z_{0k} \vert <R$. Assume, to the contrary, that 
$(0',z_{0k})\notin  J_0$ for some  $0<\vert z_{0k}\vert <R$. Then one may pick a neighbourhood $V_0$ of $(0',z_{0k})$ such that $V_0\subset \left(J_0\right)^c$ and which is of the form
$$ V_0:=\{\Vert z' \Vert <\rho\}\times \{R_1<\vert z_k\vert <R_2\;\textrm{and}\; \vert\arg z_k-\arg z_{0 k}\vert <\eta\}.$$
Let us now denote by $T_{\rho,R_1,R_2}$ the tube 
\begin{center}
$T_{\rho,R_1,R_2}:=\{\Vert z'\Vert <\rho\}\times \{R_1<\vert z_k\vert <R_2\}.$
\end{center}
Since $A_0$ is contracting and $\theta_0 / \pi$ irrational, for any $z\in T_{\rho,R_1,R_2}$ there exists an integer $n$ such that $g_0^n(z)\in V_0$. By the invariance of Julia sets we thus have $T_{\rho,R_1,R_2}\subset \left(J_0\right)^c$.
Let us shrink the tube $T_{\rho,R_1,R_2}$. By assumption, $J_\la$ is $u.s.c$ at $0$ and therefore 
\begin{center}
$T_{\rho,R_1,R_2}\subset \left(J_\la\right)^c$ when $\la$ is close enough to $0$.
\end{center}
 On the other hand, according to the second assertion of Lemma \ref{techlem}, we may find parameters $\la$ which are arbitrarily close to $0$ and such that $B_{R_1}\setminus S_\la\subset \cup_n \left(g_\la^n \right)^{-1} T_{\rho,R_1,R_2}$ where $S_\la$ denotes the stable manifold of $g_\la$.
As $\mu_\la$ gives no mass to analytic sets, this and the inclusion  $T_{\rho,R_1,R_2}\subset \left(J_\la\right)^c$ implies the existence of a sequence of parameters 
$\la_p\to 0$ such that $B_{R_1}\subset \left(J_{\la_p}\right)^c$.  This contradicts the lower semi continuity of $J_\la$ at $0$ since  
$0\notin \left(J_{\la_p}\right)_{\frac{R_1}{2}}$ but $0\in J_0$ by our assumption. \finsec

%%%%%%%%%%%%%
\subsection{Holomorphic motion of  all repelling $J$-cycles of $f$ when $f$ admits an equilibrum  web}
%%%%%%%%%%%%%%%

Our main goal here is to  end  the proof of Theorem \ref{main2}. We only have to establish the implication 
(C)$\Rightarrow$(A). To this purpose,
we shall use Corollary \ref{corSiegBif} and show how a Siegel disc may appear when a repelling $J$-cycle fails to move holomorphically.

\begin{prop}\label{thbouc} Let $f : M \times \Pj^k \to M \times \Pj^k$ be a holomorphic family such that $k=2$ or $M$ is a simply connected open subset of ${\cal H}_d(\Pj^k)$. If $f$ admits an equilibrium  web then all repelling $J$-cycles of $f$  move holomorphically. 
\end{prop}

Let us recall that a periodic point is said to be \emph{resonant} if its multipliers $w_1,\cdots,w_k$ satisfy a relation of the form $w_1^{m_1}\cdots w_k^{m_k}-w_j=0$ where the $m_j$ are integers and $m_1+\cdots + m_k\ge 2$. Note that when $w_j=e^{i\theta_j}$ for $1\le j \le s$ and some $s \le k$ then the absence of resonances forces $\pi, \theta_{1}, \cdots,\theta_{s}$ to be linearly independent over $\Q$.\\

We shall use the following Lemmas.

\begin{lem}\label{LemS}
Let $f : M \times \Pj^k \to M \times \Pj^k$ be a holomorphic family of degree $d\geq 2$,
where $M$ is an open subset of ${\cal H}_d(\Pj^k)$ and $n\geq1$
a fixed integer.
There exist an analytic subset $Y$ of $M$ and  an analytic subset $Y_0$
of $M\setminus Y$
such that
 every point $\lambda_0\in M\setminus (Y\cup Y_0)$
admits a neighbourhood $U_0$ and an analytic subset $Y'$ of $U_0$
such that:
all $n$-periodic points of $f_{\la}$ are given by holomorphic parametrizations $\zeta(\la)$ over $U_0$ 
and the eigenvalues of the Jacobian matrix of $f_\la^n$ at $\zeta(\la)$ 
give rise to a parametrization $w: \lambda\mapsto (w_1 (\la), \dots, w_k(\la))$ which satisfies the following properties:
\begin{enumerate}
\item $w$ is holomorphic,
\item the rank of the Jacobian matrix of $w$ on $U_0 \setminus Y'$ is maximal, i.e., equal to $k$
(notice that $\dim M = (k+1)(k+d)!/(k!d!)>k$),
\item $w_i (\la)\neq w_j (\la)$ for every $\la\in U_0\setminus Y'$ and $i\neq j$.
\end{enumerate}
\end{lem}

\proof
We shall use the notation $[z]$ for the elements of $\Pj^k$,  in particular we shall
denote by $f(\la,[z])$ the holomorphic family under consideration.
We only have to treat the case $M={\cal H}_d(\Pj^k)$. However,  all we shall actually need is that 
maps of the form $ [z_0^d + \eps_0 z_j^d: \dots : z_j^d:\dots : z_k^d + \eps_{k-1} z_j^d]$, where $1\le j\le k$ and 
$\eps:=(\eps_0,\dots,\eps_{k-1})$ belongs to some neighbourhood $V_0$ of $0\in\C^k$,  are parametrized by $M$.
More precisely, we will consider  the families $g^j:V_0\times\Pj^k\to V_0\times\Pj^k$ given by  
\[
g^j(\eps, [z]):= (\eps, [z_0^d + \eps_0 z_j^d: \dots : z_j^d:\dots : z_k^d + \eps_{k-1} z_j^d])
\]
as sub-families of $f$ and use them to investigate the perturbations of the map
$ f_0:=[z_0^d: \dots : z_k^d]$.\\

Given $n$ a fixed integer, let us
set
$
X_n :=\left\{ (\la,[z])\in M\times \Pj^k \colon f^n (\la,[z]) = (\la,[z]) \right\}.
$
The number of $n$-periodic points of $f_\lambda$ is bounded above by a constant which does not depend on $\lambda \in M$ (but depends on $k,d,n$). Hence the canonical projection $\pi:X_n \to M$
is a ramified covering of a certain degree $D_n$;
let  $Y\subsetneq M$ be  the analytic subset such that
 $\pi: X_n \setminus \pi^{-1} (Y) \to M\setminus Y$
 is a covering. Every $\la_0\in M\setminus Y$
 admits a neighbourhood $U_0\subset M\setminus Y$
 and holomorphic
 maps $\zeta_q: U_0\to \Pj^k$ ($1\leq q\leq D_n$)
 such that $\left\{\zeta_q (\la)\colon 1\leq q\leq D_n\right\}$
 is the set of  fixed points of $f^n_\la$.
We denote by $A_q(\la)$
the Jacobian matrix of $f_\la^n$ at $\zeta_q (\la)$, the dependence 
$\la\mapsto A_q(\la)$ is holomorphic on $U_0$.
We then denote by $\left\{w_{q,1} (\la), \dots, w_{q,k} (\la)\right\}$
the set of eigenvalues of $A_q (\la)$.

The parameter $0$ (corresponding to   $f_0= [z_0^d: \dots : z_k^d]$) does not belong to $Y$ and  therefore there exists a  neighbourhood $U_0$ of $0$ in $M$ on which all the maps $\zeta_q $  for  $1\le q\le D_n$ are defined. The key observation is as follows. 
Each map $\zeta_q$, when restricted to a suitable family $g^j(\eps, \cdot)$,
is of  the form
$\zeta_q (\eps) = (\zeta_{q,0} (\eps_0), \dots, \zeta_{q,k-1} (\eps_{k-1}))$
and the corresponding matrix-map $A_q(\epsilon)$ is diagonal with  non-constant holomorphic entries 
of the form
$w_{q,1} (\eps_0), \dots, w_{q,k}(\eps_{k-1})$.
Assume indeed that $\zeta_q(0)\in \{z_k\ne 0\}$
and  consider the family $g^k(\eps, [z])$.
In the standard chart  $\{z_k\ne 0\}$, we can write $g^k$ as
$
g^k(\eps, z):= (\eps, z_0^d + \eps_0 , \dots, z_{k-1}^d + \eps_{k-1} ).
$
Then $\zeta_q (\eps)$ takes the announced form
and the functions $w_{q,j}(\eps)$  correspond to multipliers of $n$-periodic points of the polynomial $z^d+\epsilon$. They cannot be constant since such periodic points uniformly tend to infinity with $\epsilon$ and the multiplier is given by the product of the function $dz^{d-1}$ evaluated along the points of the cycle. Moreover $w_{q,j}$ can not be persistently equal to $w_{q,j'}$ for $j \neq j'$ since these multipliers can be chosen independently from each other. \\

By the above observation, for a generic $\la \in M$ the matrix $A_q(\la)$ have $k$ distinct eigenvalues. More precisely, 
there exists an analytic subset $Y_0$
of $M\setminus Y$ (the union over $q$ of the zero sets of the discriminants of the characteristic polynomials of $A_q(\la)$)
such that, for every $\la\in M\setminus Y$ and $U_0$ as above,
the maps $\la\mapsto w_{q,j} (\la)$ (with $1\leq q\leq D_n$ and $1\leq j\leq k$)
are well defined, holomorphic and take pairwise distinct values on $U_0\setminus  Y_0$. 
We can thus define the (holomorphic) maps
\[
\psi_q : \la \mapsto
\left(
w_{q,i_1} (\la), \dots, w_{q,i_k}(\la)
\right),
\]
on some neighbourhood of any $\la_0\in M\setminus (Y\cup Y_0)$.
Although the definition  of $\psi_q$
depends on the ordering of the eigenvalues,  this is not the case for the rank $r_q(\la)$ of its Jacobian
with respect to $\la$
which is thus a well defined function on $M\setminus (Y\cup Y_0)$.
For every $\la_0 \in M\setminus (Y\cup Y_0)$
and every neighbourhood $U_0\subset M\setminus (Y\cup Y_0)$ of  $\la_0$
we have the 
following dichotomy, valid for every fixed $1\leq q\leq D_n$:
\begin{enumerate}
\item $r_q < k$ on $U_0$; or
\item $r_q = k$ outside a (proper) analytic subset $Y_q$ of $U_0$.
\end{enumerate}
We will establish that the first possibility is impossible, thus proving the lemma. Owing to analytic continuation, it is sufficient to prove the existence of 
$\la_0\in M\setminus (Y\cup Y_0)$ such that $r_q(\la_0) = k$. By our key observation, the map $\psi_q$ restricted to some convenient family $g^j(\eps, \cdot)$ is of the form 
$\left(w_{q,1} (\eps_0), \dots, w_{q,k}(\eps_{k-1})\right)$ where the $w_{q,i}$ are non-constant holomorphic functions of one variable. It then suffices to choose $\eps$ such that
the $w_{q,j}(\eps_{j-1})$ are pairwise distinct (and thus $\eps\notin Y_0$)
and the derivatives 
$w'_{q,j} (\eps_{j-1})$ non zero (so that $r_q (\eps)=k$). The desired parameter $\la_0$ is given by 
$f_{\la_0}=g^k(\eps,\cdot)$. \fin

\begin{lem}\label{newL}
Let $w_1,\cdots,w_k : D(0,R) \to \C$ be holomorphic functions. Assume that $w_j(0) \ne 0$ and that there exists $\lambda_n \to 0$ such that $\min_{1 \le j \le k} \vert w_j(\la_n) \vert  >1$. Assume moreover that there exists $1\le N\le k$ such that

- $\vert w_j(0) \vert =1$ and $w_j'(0)\ne 0$ for $1\le j\le N$,

- $\vert w_j(0) \vert \neq 1$ for $N+1 \le j \le k$.

Then, after renumbering,  there exist an integer $1\le q \le k$, a disc $D(\la_0,r)\subset D(0,R)$ and a partition $D(\la_0,r) = D^+(\la_0,r) \cup C \cup D^-(\la_0,r)$ where $C$ is a real analytic arc through $\la_0$ and $D^+$ and $D^-$ are open connected subsets of $D$
 such that
\begin{enumerate}
\item $\vert w_j \vert > 1$ on $D^+(\la_0,r)$ and $\vert w_j \vert < 1$ on $D^-(\la_0,r)$ for $k-q+1\le j\le k$,
\item $\vert w_j\vert =1$ and $w_j'\ne 0$ on $C$  for $k-q+1\le j\le k$,
\item $\vert w_j\vert >1$ on $D(\la_0,r)$ for $1\le j\le k-q$ if $q\le k-1$.
\end{enumerate}
\end{lem}

\proof In the sequel we allow to shrink $R$ without specifying it. Let us set $C_j:=\{\vert w_j\vert =1\}$ and
$U_j^{+}:=\{\vert w_j\vert >1\}$, $U_j^{-}:=\{\vert w_j\vert <1\}$. Since we can assume that $w_j'(0)\ne 0$ when $\{ \vert w_j \vert =1 \} \neq \emptyset$ the subset $C_j$ is either empty or is a real-analytic arc through $0$ in $D(0,R)$ on which $w_j'\ne 0$. In particular we have
\begin{center}
$C_j=C_l$ \ \textrm{ if }  $C_j\cap C_l$\;\textrm{is strictly bigger than}\;$\{0\}$.
\end{center} 
Let us set $U^+:=\cap_{j=1}^k U_j^+$. By assumption, $0\in \overline{U^+}$ and therefore $U^+$ is a non-empty open subset of $D(0,R)$. It is clear that
$\partial U^+ \subset \partial D(0,R) \cup \left(\cup_{j=1}^k C_j\right)$. On the other hand, we can not have $\partial U^+ \subset \{0\} \cup \partial D(0,R)$ since otherwise $U^+ = D(0,R) \setminus \{ 0 \}$ and the subharmonic function $\psi(\la):= \max_{1\le j\le k} \vert w_j(\la)\vert^{-1}$ would violate the maximum principle (recall that $\psi(0)\ge 1$). We may thus pick $\la_0\ne 0$ such that $\la_0\in C_{j_0} \cap \partial U^+$ for some $1\le j_0\le k$. Observe that $\la_0 \notin U_i^-$ for $1\le i\le k$.

If $C_i\ne C_{j_0}$ for some $1\le i\le k$ then $\la_0\notin C_i$ and thus $\la_0\in U_i^+$. After renumbering we may therefore find $1\le q \le k-1$ such that
\begin{center}
$\la_0\in C_{k-q+1}=C_{k-q+2}=\cdots = C_k:=C \ \textrm{and} \ \la_0 \in U_{1}^+\cap\cdots\cap U_{k-q}^+$.
\end{center}
For $r>0$ sufficiently small we have $D(\la_0,r) \subset \cap_{1}^{k-q} U_i^+$ and $D(\la_0,r)\setminus C$ has two connected components $\Omega_1$ and $\Omega_2$. For each $k-q+1\le i\le k$, one has $\Omega_1 \subset U_i^+$ and $\Omega_2 \subset U_i^-$ or
$\Omega_1 \subset U_i^-$ and $\Omega_2 \subset U_i^+$. Assume for instance that $\Omega_1 \subset U_{k-q+1}^+$. Then, since $\la_0\in \partial U^+$,
we must have $\Omega_1 \subset U_i^+$ and $\Omega_2 \subset U_i^-$ for every $k-q+1\le i\le k$ and we set
$D(\la_0,r)^+:=\Omega_1$ and $D(\la_0,r)^-:=\Omega_2$. \fin

\noindent \textsc{Proof of Proposition \ref{thbouc}}:  By assumption, $f$ admits an equilibrium web. Let $\la_a\in M$ and assume that
 $z_a$  belongs to some $p$-periodic repelling $J$-cycle of $f_{\la_a}$.
It suffices to show that the map $\gamma : M\to \Pj^k$, which is the  element of $\mathcal J$ given by Lemma \ref{lemMCR} ($3$), enjoys the property that $\gamma(\la) \in J_\lambda$ is repelling for every
$\la\in M$. 
%Let us observe that $\gamma(\la)$ is not persistently resonant and not persistently undiagonalizable for any $\la \in M$.

Since $M$ is connected, we have to show that the subset $\{\, \la\in M \; \colon \; \gamma(\la)\;\textrm{is repelling}\,\}$ is closed in $M$. Assume, to the contrary, that this is not true. Then, for arbitrarily small $\epsilon_0$,  one finds a new holomorphic map $\gamma_0:B_{\epsilon_0} \to \Pj^k$ such that $\gamma_0(\la) \in J_\lambda$ is fixed by $f_\la^p$ for all $\la \in B_{\epsilon_0}$ and $\gamma_0(0)$ is not repelling but $\gamma_0(\la_0)$ is repelling for some $\la_0\in B_{\epsilon_0}$. Our aim below is to find $\la_0'\in B_{\epsilon_0}$ such that $\gamma(\la_0')$ is a virtually repelling Siegel fixed point of $f_{\la_0'}^p$.  Corollary \ref{corSiegBif} then yields a contradiction. \\

Reducing $\epsilon_0$ allows to use charts and replace $\Pj^k$ by $\C^k$. 
Let us denote by $w_1(\la),\cdots,w_k(\la)$ the 
eigenvalues of $A(\la):=\left(f_\la^p\right)'(\gamma(\la))$. There exists
a proper analytic subset $Z$ of $B_{\epsilon_0}$ such that $w_1,\cdots,w_k$ are  holomorphic on 
$B_{\epsilon_0}\setminus Z$.
For every $n\in \N$ we define a function 
$\omega_n$ on $\C^l$
\begin{center}
$\omega_n(w):=\min_{2\le \vert m\vert\le n\;,\;1\le j\le l} \; \vert w_1^{m_1}\cdots w_l^{m_l} -w_j\vert$
\end{center}
where $\vert m\vert:=m_1+\cdots+m_l$ for any $m:=(m_1,\cdots,m_l)\in \N^l$
and then  define a function ${\cal B}_l$ on $\C^l$ by setting
\begin{center}
${\cal B}_l(w):=\sum_{n\ge 0} \frac{1}{2^n} \log  \omega_{2^{n+1}} (w).$
\end{center}
We shall  set ${\cal B}(\la):={\cal B}_k(w_1(\la),\cdots, w_k(\la))$.
The interest of this function is that, according to Brjuno's theorem (see \cite{Br}),  $f_\la^p$
is holomorphically linearizable at $\gamma(\la)$ if ${\cal B}(\la)>-\infty$ and  $A(\la)$ is diagonalizable.
We shall also use the fact that each function ${\cal B}_l$ is finite on a dense subset of the real torus $T^l:=\{\vert w_1\vert=\cdots=\vert w_l\vert =1\}$.\\

Let us denote by $\Delta_{\epsilon_0}$ the disc in $\C$ obtained by intersecting $B_{\epsilon_0}$ with the complex line through $0$ and $\la_0$.
We may move a little bit $\la_0$ so that  the set $Z\cap \Delta_{\epsilon_0}$ is discrete. 
%and $\gamma_0(\la)$ is not persistently undiagonalizable on $\Delta_{\epsilon_0}$. 
In particular, there exists a discrete subset $Z_0$ of $\Delta_{\epsilon_0}$
such that on $\Delta_{\epsilon_0}\setminus Z_0$ 
%the cycle $\gamma_0(\la)$ is diagonalizable and 
the functions $w_1,\cdots,w_k$ are either constant or holomorphic, non-vanishing and with non-vanishing derivatives. 
When $M$ is an open subset of  ${\cal H}_d(\Pj^k)$, we shall need the following more precise fact which immediately follows from Lemma \ref{LemS}.\\

{\bf Fact}: {\it If $M$ is an open subset of ${\cal H}_d(\Pj^k)$, we may move $\la_0$ so that there exists a set $Z'_0 \supset   Z_0$ for which
$\Delta_{\epsilon_0}\setminus Z'_0$ is open, path-connected and dense in $\Delta_{\epsilon_0}$, the functions $w_1,\cdots,w_k$ take pairwise distinct values on  
$\Delta_{\epsilon_0}\setminus Z'_0$ and the maps $\la \mapsto (w_{i_1}(\la),\cdots,w_{i_k}(\la))$ are well defined holomorphic submersions on some neighbourhood of any 
$\la_0\in \Delta_{\epsilon_0}\setminus Z'_0$ in $M$.}\\

Let us set 
$$ \varphi(\la):=\min \left(\vert w_1(\la)\vert,\cdots,\vert w_k(\la)\vert\right) . $$
This is a continuous function on $\Delta_{\epsilon_0}$. Moreover $\varphi(0)\le 1$ and $\varphi(\la_0)>1$, in particular $\varphi$ is not constant.  By the maximum principle applied to the subharmonic function $ \varphi^{-1}$, there exists $\la_1 \in \Delta_{\epsilon_0}\setminus Z_0$ such that $\varphi(\la_1)<1$. %Indeed, if $\varphi \ge 1$ on $\Delta_{\epsilon_0}\setminus Z_0$, then $\varphi\ge 1$ on $%\Delta_{\epsilon_0}$ and therefore the subharmonic function $\psi := \varphi^{-1}$ violates the maximum %principle (indeed $\psi\le 1=\psi(0)$ and this function is not constant). 
Considering a continuous path connecting $\la_0$ to $ \la_1$ in $\Delta_{\epsilon_0}\setminus Z_0$, one finds $\la_2\in \Delta_{\epsilon_0}\setminus Z_0$ and $\tilde\la_p\to \la_2$ such that $\varphi(\la_2)=1$ and $\varphi(\tilde\la_p)>1$.
Let us pick  a small disc $D(\la_2,R)$ contained in  $\Delta_{\epsilon_0}\setminus Z_0$. Then (after renumbering) the functions $w_1,\cdots,w_k$ satisfy the assumptions of Lemma
\ref{newL} on $D(\la_2,R)$. Let $q$ be the integer and $C$ be the real analytic arc in $D(\la_3,r)\subset D(\la_2,R)$ provided by this Lemma. Since
$\vert w_j \vert < 1$ on $D^-(\la_3,r)$ for $k-q+1\le j\le k$ and $\gamma(\la)\in J_\la$, we must have $1\le q\le k-1$.

Let us stress that when $M$ is an open subset of ${\cal H}_d(\Pj^k)$, we may use the Fact and repeat the above argument with $Z'_0$ replacing $Z_0$. This gives us a point
$\la_3$ with a neighbourhood $V_{\la_3}$ on which $A(\la_3)$ has pairwise distinct eigenvalues and  the map $\la \mapsto (w_1(\la),\cdots,w_k(\la))$ is a holomorphic submersion.
We shall denote by $S_q$ the piece of real manifold defined by 
$$S_q:=\{\la\in V_{\la_3}\;\colon\; \vert w_{k-q+1}(\la)\vert=\cdots= \vert w_{k}(\la)\vert=1\}.$$
 Note that, if $r$ is small enough,  the arc $C$ is contained in $S_q$.\\

Let us explain how to end the proof. Recall that we seek a parameter $\la'_0$ such that 
$\gamma(\la_0')$ is a virtually repelling Siegel fixed point of $f_{\la_0'}^p$. By definition,   if $\la'_0 \in S_q$ and $f_{\la'_0}^p$ 
 is holomorphically linearizable at $\gamma(\la'_0)$ then $\gamma(\la'_0)$ is a Siegel  fixed point of $f_{\la'_0}^p$ (provided that $\pi$ and the angles of rotation are linearly independent over $\Q$). Let us now observe that, since $\la_3$ was produced by  Lemma \ref{newL}, such a point $\gamma(\la'_0)$ is virtually repelling when $\la'_0$  is sufficently close to $\la_3$. We thus have to find  points on $S_q$ where Brjuno theorem applies and which are  arbitrarily close to $\la_3$. 
 This is easier when $k=2$, in higher dimension we need to assume that
$M$ is an open subset of ${\cal H}_d(\Pj^k)$ and use the above Fact.\\

We first justify the existence of such a parameter $\la_0'$ in dimension $k=2$. In that case, $A(\la)$ is diagonalizable  when $\la \in C$ since  only one of its two eigenvalues
has modulus $1$. Let us inspect the quantities $\vert w_1(\la)^{m_1} w_2(\la)^{m_2} -w_j(\la)\vert$
in the definition of $\omega_n(\la)$. Since $\vert w_1\vert >1$ and $\vert w_2 \vert =1$ on $C$,
they all are larger than some uniform constant $a>0$ on $C$ except  those of the form
$\vert w_1(\la) w_2(\la)^{m_2} -w_1(\la)\vert$ with $m_2\ge 1$ or 
$\vert w_2(\la)^{m_2} -w_2(\la)\vert$ with $m_2\ge 2$. It follows that 
${\cal B}(\la) >-\infty$ if ${\cal B}_1(e^{i\theta(\la)}) >-\infty$ 
 where $\theta$ is a non constant real analytic function on $C$ such that $w_2(\la)=e^{i\theta(\la)}$.
The existence of $\la'_0 \in C$ is thus deduced from the fact that ${\cal B}_1$ is finite on a dense subset of the torus $T^1$.\\

We finally consider the case where $M$ is an open subset of ${\cal H}_d(\Pj^k)$. 
As we already saw, in that situation, the point $\la_3\in S_q$ enjoys the following two properties: $A(\la_3)$ has pairwise distinct eigenvalues  and  the map $\la \mapsto (w_1(\la),\cdots,w_k(\la))$ is a holomorphic submersion on a neighbourhood $V_{\la_3}$ of $\la_3$. Let us write $w_{j}(\la) =:e^{i\theta_j(\la)}$ for $k-q+1\le j\le k$ and $\la\in S_q$.

By the submersion property, the sets $\{\la\in S_q\;\colon\; \vert w_1(\la)\vert^{m_1}\cdots
\vert w_{k-q}(\la)\vert^{m_{k-q}} =\vert w_j(\la)\vert\}$ are either empty or real hypersurfaces in $S_q$. By Baire
theorem, we may therefore move a little bit $\la_3$ in $S_q$ so that 
\begin{eqnarray}\label{NRM}
\vert w_1(\la_3)\vert^{m_1}\cdots
\vert w_{k-q}(\la_3)\vert^{m_{k-q}} \ne \vert w_j(\la_3)\vert \; \textrm{for}\; 1\le j\le n \; \textrm{and}\; m_1+\cdots +m_{k-q} \ge 1.
\end{eqnarray}
On the other hand, since $\vert w_j(\la_3)\vert >1$ for $1\le j\le k-q$, we may shrink $V_{\la_3}$ and  find $A\in \N$ and $c>0$ such that
\begin{eqnarray}\label{NR}
\vert w_1(\la)^{m_1}\cdots
 w_{k}(\la)^{m_{k}} -  w_j(\la)\vert \ge c  \; \textrm{on}\; V_{\la_3}\; \textrm{for}\; 1\le j\le n\;\textrm{and}\; m_1+\cdots+ m_{k-q} \ge A.
\end{eqnarray}
From (\ref{NRM}) we get $\vert w_1(\la_3)^{m_1}\cdots
 w_{k}(\la_3)^{m_{k}} -  w_j(\la_3)\vert \ge \vert \vert w_1(\la_3)^{m_1}\cdots
 w_{k-q}(\la_3)^{m_{k-q}} \vert -\vert w_j(\la_3)\vert \vert >0$
 $\textrm{for}\; 1\le j\le n \; \textrm{and}\; m_1+\cdots +m_{k-q} \ge 1$  and thus, after shrinking $V_{\la_3}$ again
 and taking $c$ smaller, one has 
$\vert w_1(\la)^{m_1}\cdots
 w_{k}(\la)^{m_{k}} -  w_j(\la)\vert \ge c$ on $V_{\la_3}$ for   $1\le m_1+\cdots+ m_{k-q} \le A$ and $1\le j\le n$. Taking (\ref{NR}) into account we thus have 
 \begin{eqnarray*}
\vert w_1(\la)^{m_1}\cdots
 w_{k}(\la)^{m_{k}} -  w_j(\la)\vert \ge c \; \; \textrm{on}\; V_{\la_3}\;\textrm{for}\; 1\le j\le n\;\textrm{and}\;   m_1+\cdots+ m_{k-q}
 \geq 1.
\end{eqnarray*}

It follows that, for $\la \in V_{\la_3}\cap S_q$,  
${\cal B}(\la) >-\infty$ if ${\cal B}_q(e^{i\theta_{k-q+1}(\la)},\cdots,e^{i\theta_{k}(\la)}) >-\infty$.
The existence of $\la'_0 \in S_q$ satisfying ${\cal B}(\la'_0) >-\infty$ is thus obtained by using the fact that ${\cal B}_q$ is finite on a dense subset of $T^q$ and, once again, the submersion property. Note that if $\la'_0$ is close enough to $\la_3$ then
$A(\la'_0)$ is diagonalizable since its eigenvalues are pairwise distinct. Note also that
$\theta_{k-q+1}(\la'_0),\cdots,\theta_k(\la'_0), \pi$ are linearly independent over $\Q$ since
 ${\cal B}_q(e^{i\theta_{k-q+1}(\la'_0)},\cdots,e^{i\theta_{k}(\la'_0)}) >-\infty$.
\fin

\begin{rem}\label{RemContSieg}
It would be interesting to know if  the continuity of the map $\la \mapsto J_\la$ on some open subset of the parameter space  is equivalent to the existence of an equilibrium  web, as it occurs  when $k=1$. The second item of Proposition \ref{SiegJul} should be useful to study this question.
\end{rem}

\textsc{Proof of Theorem \ref{main2}}:
Let $f:M\times \Pj^{k}\to M\times \Pj^{k}$
be a holomorphic family of endomorphisms where $M$ is either a simply connected open subset of  ${\cal H}_d(\Pj^k)$ or any simply connected complex manifold when $k=2$.\\
In subsection \ref{sspartof} we saw that
 (A)$\Rightarrow$(B)$\Leftrightarrow$(E). Theorem \ref{main} yields (B)$\Rightarrow$(C'), where (C') is the assertion : "the restriction $f_{B\times \Pj^k}$ admits an equilibrium  web for any sufficiently small ball $B$". Assume now that (C') is satisfied. By Proposition \ref{thbouc},  the repelling $J$-cycles locally move holomorphically.  This implies that the set 
\begin{eqnarray*}
\{  (\la,z)\in M  \times \Pj^k \; \colon \; z \; \textrm{belongs to some} \; \textrm{$n$-periodic and repelling $J$-cycle of} \; f_\la  \}
\end{eqnarray*}
 is an unramified covering of $M$. As $M$ is simply-connected, we thus get  that
the repelling $J$-cycles move holomorphically over $M$, hence (C')$\Rightarrow$(C). Proposition \ref{thbouc} also yields (C)$\Rightarrow$(A), 
and therefore the properties  (A), (B) and (C) are equivalent.\\
If (D) is satisfied then by definition  any element $\gamma$ of the equilibrium lamination belongs to $\cal J$ and satisfies $\Gamma_{\gamma} \cap PC_f =\emptyset$.
Then the first assertion of Proposition \ref{propgraph} shows that $f$ admits an equilibrium  web. We thus have (D)$\Rightarrow$(C). 
Finally, since  (A)$\Rightarrow$(D) by Theorem \ref{PropExHM}  and (B)$\Leftrightarrow$(F) by Theorem \ref{main}, 
 the proof of Theorem \ref{main2} is completed. \finsec

%%%%%%%%%%%%%
\section{Bifurcation loci}
%%%%%%%%%%%%%%%

In view of Theorem \ref{main2}, we define the bifurcation locus and  current as follows.

\begin{defn}
Let $f:M\times \Pj^{k}\to M\times \Pj^{k}$
be a holomorphic family of endomorphisms of $\Pj^k$ of degree $d \geq 2$. 
Let $L(\la)$ be the sum of Lyapunov exponents of $f_\la$ with respect to its equilibrium measure.
The closed positive current $dd^c_\la  L$
 is called the \emph{bifurcation current} of the family, its support is the \emph{bifurcation locus} of the family.
\end{defn}

We will exploit here our results  to get some informations on these loci.

  %%%%%%%%%%%%%%%%
  \subsection{Remarkable elements in bifurcation loci}
  %%%%%%%%%%%%%%%%
  
 Theorem \ref{main2} and the proof of Proposition \ref{thbouc} immediately yield the following result.

\begin{thm}\label{theoAdS}
 A degree $d \geq 2$ endomorphism of $\Pj^k$ belongs to the bifurcation locus in ${\cal H}_d(\Pj^k)$ if and only if it is accumulated by endomorphisms which admit a virtually $J$-repelling Siegel periodic point or a repelling cycle outside the Julia set which becomes a repelling $J$-cycle after an arbitrarily small perturbation. 
\end{thm}

The next theorem shows that \emph{isolated} Latt\`es maps belong to the bifurcation locus.
We refer to the articles \cite{Di1}, \cite{Du1} for an account on Latt\`es maps of $\Pj^k$.

\begin{thm}\label{thLa}
Let $f:M\times \Pj^{k}\to M\times \Pj^{k}$ be a  holomorphic family of endomorphisms of $\Pj^k$.
If the family is stable (i.e. $dd^c_\la  L = 0$ on $M$) and $f_{\la_0}$ is a Latt\`es map for some $\la_0 \in M$ then $f_\la$ is a Latt\`es map for every $\la \in M$. 
\end{thm}

\proof 
By a Theorem of Briend-Duval \cite{BD1} we have $L \ge k \frac{\log d}{2}$. The articles of 
Berteloot, Dupont and Loeb \cite{BL}, \cite{BDu} and \cite{Du2} show that $L(\la) = k \frac{\log d}{2}$ if and only if $f_\la$ is a Latt\`es map.
If the family is stable,  then the function $L$ is pluriharmonic on $M$. By the maximum principle  (applied to the harmonic function $-L$) we thus have $L(\la)=L(\la_0)= k\frac{\log d}{2}$ for all $\la \in M$ and the conclusion follows.\finsec

%%%%%%%%%%%%%%%%%%%%%%%%%%%%
\subsection{On the interior of  bifurcation loci}
%%%%%%%%%%%%%%%%%%%%%%%%%%%%
In his work on the persistence of homoclinic tangencies, Buzzard \cite{Bu} found open subsets of the space of degree $d$ endomorphisms of $\Pj^2$ (for $d$ large enough) in which the maps having infinitely many sinks are dense. This lead us to believe that the bifurcation locus may have a non-empty interior when $k\ge 2$.
We investigate here the relations between the presence of open subsets in  the support of $dd^c_\la  L$ and the existence of parameters for which the postcritical set is dense in $\Pj^k$
and in particular prove Theorem \ref{thmopbif}.\\

Let $f:M\times \Pj^{k}\to M\times \Pj^{k}$ be a holomorphic family of endomorphisms of $\Pj^k$. Let $C$ denote the critical set of $f$ and let $C_\la$ denote the critical set of $f_\la$. We set
  \begin{eqnarray*}
  \overline{C^+} := \overline{\cup_{n\ge 1} f^n (C)} \ \ \ \ \textrm{and}  \ \ \ \ \ 
  \overline{C_\lambda^+}  := \overline{\cup_{n\ge 1} f_\lambda^n (C_\la) } \ \ \textrm{for every} \  \la \in M. 
  \end{eqnarray*}
We define $( \overline{C^+})_\la:= \left(\{\la\}\times \Pj^k\right)  \cap \overline{C^+}$, observe that $\{\la\} \times  \overline{C_\lambda^+}  \subset ( \overline{C^+})_\la$.
Our aim is to show that if $\supp dd^c_\la  L$ contains an open subset $\Omega$, then $\{\la \in \Omega \;\colon\;  \overline{C_\lambda^+} =\Pj^k\}$ contains a $G_\delta$-dense subset of $\Omega$.\\

Note that this result sheds new light on the well-known fact that stable parameters are dense for holomorphic families of rational maps. For such families the bifurcation locus is known to coincide with $\supp dd^c_\la  L$ (\cite{dM}).
 
 \begin{cor}
 Let $f:M\times \Pj^{1}\to M\times \Pj^{1}$ be a holomorphic family of rational maps. Then $\supp dd^c_\la  L$ has empty interior.
 \end{cor} 
 
 \proof Every  $\la_0 \in \supp dd^c_\la  L$ can be approximated by parameters $\la$ for which $f_\la$ has 
an attracting basin, see \cite[section 4.3.1]{Be}, which is an open condition in $M$. On the other hand, as the critical set is finite, the set $\overline{C_\lambda^+}$ can not be  equal to $\Pj^1$ when $f_\la$ has an attracting basin.
 According to Theorem \ref{thmopbif}, this implies that $\supp dd^c_\la  L$ has empty interior. \finsec

\begin{rem} We raise the question, for $k \geq 2$, of the existence of holomorphic families for which $\supp dd^c_\la  L$ has non empty interior. Note that Theorem \ref{thmopbif} could be useful for finding families for which  $\supp dd^c_\la  L$ has empty interior. 
\end{rem} 

The proof of Theorem \ref{thmopbif} relies on a Baire's category argument based on the continuity properties of  
$\la\mapsto \overline{C_\lambda^+}$ and $\la\mapsto (\overline{C^+})_\la$.  The notion of semi continuity with respect to the Hausdorff topology has been discussed in subsection \ref{ssContJ}. We have the following properties, the upper semi continuity can be found in \cite[Proposition 2.1]{Do}, we give the argument for sake of completeness.

\begin{lem}\label{doudou}
The maps $\lambda \mapsto ( \overline{C^+})_{\la}$ and $\lambda \mapsto \overline{C_\lambda^+}$ from $M$ to ${\textsc Comp}^{\star}\left(\Pj^k\right)$ are respectively upper and lower semi continuous.
\end{lem}

\proof
By definition $\{ (\lambda, z ) \in M \times \Pj^k \, , \,  z \in (\overline{C^+})_{\la} \}$ is equal to $\overline{C^+}$, hence is closed in $M \times \Pj^k$.
In particular, for every $\lambda_0 \in M$  and $\epsilon > 0$, the set $F :=  \{  (\lambda , z ) \in  \overline{C^+} \ , \  d_{\Pj^k}( z , (\overline{C^+})_{\lambda_0} ) \geq \epsilon \}$ is a closed subset of $\overline{C^+}$. Let us show that $\pi_M(F)$ is closed in $M$. Indeed, if $\lambda_n \in \pi_M(F)$ converges  to $\lambda \in M$ one may pick $z_n \in 
(\overline{C^+})_{\la_n}$ such that $d_{\Pj^k}(z_n , (\overline{C^+})_{\la_0}) \geq \epsilon$ and $(z_n)_n$ converges to some $z \in \Pj^k$ after taking a subsequence. Then $(\lambda_n , z_n) \in \overline{C^+}$ converges to $(\lambda,z) \in \overline{C^+}$ satisfying  $d_{\Pj^k}( z , (\overline{C^+})_{\lambda_0} ) \geq \epsilon$ and thus $\lambda \in \pi_M(F)$ as desired. Since $\la_0 \notin \pi_M(F)$ it follows that $M \setminus \pi_M(F)$ contains an open ball $B$ centered at $\lambda_0$ such that $d_{\Pj^k}( z , (\overline{C^+})_{\lambda_0} ) < \epsilon$ for every $z \in (\overline{C^+})_{\lambda}$ with $\lambda \in B$. This proves the upper semi continuity. 

Let us now prove the lower semi continuity of the map $\lambda \mapsto \overline{C_\lambda^+}$. Assume to the contrary that it is not \emph{l.s.c} at $\la_0 \in M$. Then there exist $\epsilon >0$, a sequence $(\la_n)_n$ converging to $\la_0$ and a sequence $(z_n)_n$  in $\overline{C_{\lambda_0}^+}$ such that $d_{\Pj^k}(z_n, \overline{C_{\lambda_n}^+} ) \ge \epsilon$. After taking a subsequence  $(z_n)_n$ converges to $z_0\in \overline{C_{\lambda_0}^+}$. Pick $\xi_0 \in C_{\la_0}$ and $p_0 \ge 1$ such that $d_{\Pj^k}(z_0,f_{\la_0}^{p_0} (\xi_0)) < \frac{\epsilon}{4}$. Let also $\xi_n \in C_{\la_n}$ such that $\xi_n\to \xi_0$. Then $d_{\Pj^k}(z_n, \overline{C_{\lambda_n}^+} ) \le d_{\Pj^k}(z_n,f_{\la_n}^{p_0} (\xi_n)) <\frac{\epsilon}{2}$ for $n$ large, contradicting $d_{\Pj^k}(z_n, \overline{C_{\lambda_n}^+} )  \geq \epsilon$. \fin

Lemma \ref{doudou} allows us to prove: 

\begin{prop}\label{propbifcrit}
Let $f:M\times \Pj^{k}\to M\times \Pj^{k}$ be a holomorphic family 
  of endomorphisms of $\Pj^k$. If $\la_0 \in \supp dd^c_\la L$  then $( \overline{C^+})_{\la_0} =\Pj^k$.
\end{prop}

\proof Assume that $B(z_0,r) \cap (\overline{C^+})_{\la_0} =\emptyset$ and let us show that $\la_0\notin \supp dd^c_\la  L$.  Since $\la\mapsto (\overline{C^+})_{\la}$ is upper semi continuous we deduce that $B(z_0,\frac{r}{2}) \cap (\overline{C^+})_{\la} =\emptyset$ when $\la$ is sufficiently close to $\la_0$. In particular, the constant graph $\Gamma_0:=\{(\la,z_0)\; \colon \; \la\in B(\la_0,\epsilon)\}$ does not meet 
$\cup_{n\ge 1} f^n (C)$ for $\epsilon$ small enough. By the first assertion of Proposition \ref{propgraph} and Proposition \ref{thmHMR}, we get
$dd^c_\la  L  = 0$ on $B(\la_0,\epsilon)$. \fin

\noindent\textsc{Proof of theorem \ref{thmopbif} :}  The lower semi continuity of  $\la \mapsto \overline{C_\lambda^+}$ implies that $$I(B):=\{\la \in M\;\colon\; \overline{C_\lambda^+} \cap B\ne \emptyset\}$$
is an open subset of $M$ for every open ball $B \subset \Pj^k$.
 Now let $\Omega$ be an open subset of $M$ which is contained in the bifurcation locus. Let us show that $I(B)$ is dense in $\Omega$. We may assume that $\Omega$ is a ball in $\C^m$.
 Let $\la_0\in \Omega$ and $\epsilon >0$. Since $\la_0 \in \supp dd^c_\la  L$,  Proposition \ref{propbifcrit} implies that
 $( \overline{C^+})_{\la_0} \cap B =B$. Thus $\big(\cup_{n\ge 1} f^n (C)\big) \cap \big( B(\la_0,\epsilon)\times B\big) \ne \emptyset$ and there exists $(\la_1,z_1)\in f^{n_1}(C)\cap \big( B(\la_0,\epsilon)\times B\big)$. This shows that $\la_1\in I(B)\cap B(\la_0,\epsilon)$ and thus $I(B)$ is open and dense in $\Omega$. Now consider a countable collection $B_i:=B(\z_i,r_i)$ of balls in $\Pj^k$ whose centers are dense in $\Pj^k$ and whose radii tend to $0$. According to Baire's theorem $M' := \cap_{i\ge 1} I(B_i)$ is a dense $G_\delta$-subset of $\Omega$. We also have $\overline{C_\lambda^+}=\Pj^k$ for every $\la \in M'$. \finsec

\appendix

%%%%%%%%%%%%%%%%%%%%%
\section{}
%%%%%%%%%%%%%%%%%%%%%

%%%%%%%%%%%%%%%%%%%%%%%%%%%
\subsection{Hyperbolic sets and holomorphic motions}\label{HMHS}
%%%%%%%%%%%%%%%%%%%%%%%%%%%

\begin{defn}\label{defihypmot} Let $f : B \times \Pj^k \to B \times \Pj^k$ be a holomorphic family of endomorphisms 
where $B$ is a ball centered at the origin in $\C^m$.
 Let $E_0$ be an $f_0$-invariant subset of $\Pj^k$. A \emph{holomorphic motion} of $E_0$ over $B_\rho \subset B$ is  a continuous map $h : B_\rho \times E_0 \to \Pj^k$ such that :
 \begin{enumerate}
\item $\lambda \mapsto h_\lambda(z)$ is holomorphic on $B_\rho$  for every $z \in E_0$.
\item $z \mapsto h_\lambda(z)$ is injective on $E_0$ for every $\lambda \in B_\rho$.
\item $h_\lambda  \circ f_0 = f_\lambda \circ h_\lambda$ on $E_0$ for every $\lambda \in B_\rho$.
\end{enumerate}
One says that $E_0$ is a \emph{hyperbolic set} for $f_0$ if it is $f_0$-invariant and if there exists $K>1$ such that $\vert (d f_0 )^{-1} \vert^{-1} \geq K$ on $E_0$.
\end{defn}

\begin{thm} \label{hypmot} Let $f : B \times \Pj^k \to B \times \Pj^k$ be a holomorphic family of endomorphisms. Let $E_0 \subset \Pj^k$ such that $\vert (d f_0 )^{-1} \vert^{-1} \geq K > 3$ on $E_0$. Then there exists a holomorphic motion $h : B_\rho \times E_0 \to \Pj^k$ which preserves repelling cycles.
\end{thm}

The proof is based on classical arguments, we refer to \cite[chapter 3, section 2.d]{dMvS} for the one dimensional case. To simplify the exposition we assume that the dilation is larger than $3$ on the hyperbolic set.

\proof  Let  $\varphi (z) := \inf_{\lambda \in B_\rho} \vert (d_z f_\lambda)^{-1} \vert ^{-1}$, with the convention $\vert (d_z f_\lambda)^{-1} \vert ^{-1} = 0$ if  $z \in C_{f_\lambda}$. This is a continuous function on $\Pj^k$. By taking a smaller $\rho$, we may assume that 
\begin{equation}\label{pinch}
   \varphi \geq K' > 3 \ \textrm{ on a $\tau$-neighbourhood } \left(E_0\right)_\tau . 
\end{equation}
We shall mainly use the lower estimate on $E_0$ itself, the lower bound on  $\left(E_0\right)_\tau$ appears at  the end of the proof. Let $\delta = \delta(\rho) := \min \{   (1+\sup_{\la \in B_\rho } \norm{f_\lambda}_{C^2})^{-1} , \tau \}$. 

\begin{lem} \label{til} For every $(\lambda,z) \in B_\rho \times E_0$,  
\begin{enumerate}
\item $d_{\Pj^k}(f_\lambda(z),f_\lambda(w)) \geq (K'-1) d_{\Pj^k}(z,w)$ for every $w \in \bar B(z,\delta)$,
\item $f_\lambda (B(z, c\delta))\supset B(f_\lambda(z),c\delta)$ for every $0 \leq c \leq 1$,
\item  if $g_{\lambda,z} :B(f_\lambda(z),\delta) \to B(z, \delta)$ is the inverse map of $f_\lambda$, then $\Lip g_{\lambda,z} \leq (K'-1)^{-1}$.
\end{enumerate}
\end{lem}

\proof
Assertions 2 and 3 follow from the first one (use Jordan's theorem and $K' > 3$ for the second one). So let us prove Assertion 1. We work in local coordinates. For $(\lambda,z) \in B_\rho \times E_0$ and $w \in \bar B(z,\delta)$ we have
\begin{align} \label{ineg}
\vert \Id_{\C^k}- (d_z f_\lambda)^{-1} \circ d_w f_\lambda \vert & \leq \vert(d_z f_\lambda)^{-1} \vert \cdot \vert  d_z f_\lambda - d_w f_\lambda \vert \notag \\
                                               & \leq \vert (d_z f_\lambda)^{-1} \vert \cdot \vert z-w \vert \cdot \delta^{-1}  \leq 1/K' \notag .
\end{align}
That implies $\textrm{Lip} \, (\Id - (d_z f_\lambda)^{-1} \circ f_\lambda) \leq 1/K'$ on $\bar B(z,\delta)$, which gives in turn  
\[ \vert (d_z f_\lambda)^{-1} (f_\lambda (z) - f_\lambda (w)) \vert   \geq (1-1/K') \vert z-w \vert  \]
for every $w \in \bar B(z,\delta)$. Hence $\vert f_\lambda (z) - f_\lambda (w) \vert   \geq (K'-1) \vert z-w \vert$ as desired. \fin

\begin{lem}\label{typ}
For every $(\lambda,z) \in B_\rho \times E_0$, we have $B(f_\lambda(z),\delta) \supset B(f_0(z),\delta/2)$ and the inverse map $g_{\lambda,z}: B \left(f_\la(z),\delta\right) \to B(z,\delta)$ given by Lemma \ref{til}  satisfies:

\begin{enumerate}
\item $g_{\lambda,z}$ is well defined on $B(f_0(z),\delta/2)$,
\item it satisfies  $\Lip g_{\lambda,z} \leq (K'-1)^{-1}$ on $B(f_0(z),\delta/2)$,
\item $g_{\lambda,z}(B(f_0(z),\delta/2)) \subset B(z,\delta/2)$.
\end{enumerate}
\end{lem}

\proof Let $Q := \max \, \{ \, \norm{d_\lambda f_\lambda(z)} \, , \, (\lambda,z) \in B_\rho \times E_0 \, \}$. 
As $\delta$ is a continuous function of $\rho$ and $\delta(0) > 0$, we may assume $\delta \geq 2Q\rho$ by taking $\rho$ small enough. For every $\lambda \in B_\rho$ and $z \in E_0$, $d(f_\lambda(z),f_0(z)) \leq Q \rho \le\delta/2$.  That yields $B(f_\lambda(z),\delta) \supset B(f_0(z),\delta/2)$. Items 1 and 2 are obvious from lemma \ref{til}. For item 3, we use $g_{\lambda,z}(B(f_0(z),\delta/2)) \subset g_{\lambda,z} (B(f_\lambda(z),\delta))$, which  is included in $B(z,\delta/2)$ by using lemma \ref{til}(3) and $K' > 3$.   \fin

Let us end the proof of theorem \ref{hypmot}. For $(\lambda , z ) \in B_\rho \times E_0$ we set  $z_n := f_0^n(z)$ and   
\[  g_{\lambda,z}^n := g_{\lambda,z} \circ \ldots \circ g_{\lambda,z_{n-1}} . \]
This is an inverse branch of $f_\lambda ^n$. Since $z_{1},\cdots,z_{n-1} \in E_0$, lemma \ref{typ} yields by induction
\[ g_{\lambda,z}^n  : B(z_n ,\delta/2) \to  B(z,\delta/2) \;\textrm{and}\; \Lip g_{\lambda,z}^n \leq (K'-1)^{-n} \textrm{ on } B(z_n ,\delta / 2) . \]
For $(\lambda,z)\in  B_\rho\times E_0 $ let us define  \[ h_n(\lambda,z) := g_{\lambda,z}^n \circ f_0^n(z) = g_{\lambda,z}^n(z_n). \]
The map $h_n$ is continuous in $(\lambda,z)$, holomorphic in $\lambda$ and $h_n(\lambda,z) \in B(z,\delta/2)$. Moreover
\begin{equation}\label{com} 
f_\lambda \circ h_n(\lambda,z) = h_{n-1}(\lambda,f_0(z)) . 
\end{equation}
The sequence $(h_n)_n$ is uniformly Cauchy on $B_\rho \times E_0$. Indeed 
$h_{n+1}(\lambda,z) - h_n(\lambda,z) = g_{\lambda,z}^n \circ g_{\lambda,z_n} (z_{n+1}) - g_{\lambda,z}^n (z_n)$ and we get
$\norm{h_{n+1} - h_n }_{B_\rho \times E_0} \leq  (\delta/2)  \cdot (K'-1)^{-n}$ 
since $g_{\lambda,z_n} (z_{n+1}) \in B(z_n, \delta/2)$ by Lemma \ref{typ}(3).
We  define  $h_\la(z)$ for $(\la,z) \in B_\rho\times E_0$ by 
\[ h_\lambda(z) := \lim _{n} h_n(\lambda,z) = \lim_n g_{\lambda,z}^n \circ f_0^n(z) . \]
The map $h$ is continuous in $(\lambda,z)$, holomorphic in $\lambda$ and $h_\lambda(z) \in \bar B(z,\delta/2)$. It also follows from (\ref{com}) that
\begin{equation}\label{comnn} 
  f_\lambda \circ h_\lambda = h_\lambda \circ f_0  . 
  \end{equation}
Let us now check that $h_\lambda$ is injective. Assume $h_\lambda(z) = h_\lambda(z')$. Iterating (\ref{comnn}) yields  $h_\lambda (f_0^n(z)) = h_\lambda (f_0^n(z'))$.
As $h_\lambda(w) \in \bar B(w,\delta/2)$ for $w \in E_0$,
we get $d(f_0^n(z) , f_0^n(z')) \leq \delta$. Then, since $d(f_0^n(z) , f_0^n(z')) \geq (K'-1)^n d(z,z')$ by Lemma \ref{til}(1), we must have $z=z'$.\\
 Finally,  $h_\la$ preserves cycles (see (\ref{com})) and any periodic $h_\la(z)$ must be repelling since $h_\lambda(z) \in \bar B(z,\delta/2) \subset \left(E_0\right)_\tau$ and  $\vert (d f_\lambda)^{-1} \vert ^{-1}  > 3$ on $\left(E_0\right)_\tau$ (see (\ref{pinch})). This completes the proof. \finsec

%%%%%%%%%%%%%%%%%%%%%%%%%
\subsection{Proof of Proposition \ref{LemSF}}
%%%%%%%%%%%%%%%%%%%%%%%%%

We work with the notations of Section \ref{secFrom}. 
Let $\tau , \epsilon >0$ such that $-\frac{\log  d}{2} +\tau + 2\epsilon <0$. Recall that the distortion of the  charts is controlled by $\tau$, see  Equation (\ref{LipChart}). 
Let $p \geq 1$  and $r_p(\gamma) = \inf_{\la\in {U_0}} \Vert (DF_{\gamma(\la)}^p(0))^{-1}\Vert^{-2}$, see Equation (\ref{new}). 
The next lemma shows that $r_p$ measures the size of tubular neighbourhoods of $\Gamma_\gamma$ on which $f^p$ is invertible and contracting. 

\begin{lem}\label{Lemr} 
For every small $\epsilon >0$ there exists $C_p(\epsilon) >0$ such that  for any $\gamma\in {\cal X}$ the map $f^p$ admits an inverse branch
$(f^p)_{\gamma}^{-1}$ on the tube $ T_{U_0}({\cal F}^p(\gamma),C_p(\epsilon)r_p(\gamma))$ which maps 
$\Gamma_{{\cal F}^p(\gamma)}\cap(U_0\times\Pj^k)$ to $\Gamma_\gamma\cap(U_0\times\Pj^k)$ and satisfy
$\mbox{Lip} (f^p)_{\gamma}^{-1} \le e^{\tau + \epsilon/3} r_p(\gamma)^{-1/2}$ .
\end{lem}

\proof We use a quantitative version of the inverse mapping theorem,   see \cite[Lemme 2]{BD2}. This version is more precise than Lemma \ref{til}. Let $M :=\sup_{\la\in {U_0},\gamma\in{\cal X}} \Vert F_{\gamma(\la)}^p\Vert_{{\cal C}^2,\overline{B(0,R_p)}}$ and let $\delta_p(\epsilon) := R_p(1-e^{-\epsilon/3})/M$. Then for every $(\gamma,\la) \in {\cal X}\times {U_0}$:
 \begin{enumerate} 
\item[$\cdot$]  $(F_{\gamma(\la)}^p)^{-1}$ is defined on $B_{\C^k} \left(0,\delta_p(\epsilon)\Vert (DF_{\gamma(\la)}^p(0))^{-1}\Vert^{-2} \right)$,
\item[$\cdot$] $\mbox{Lip}(F_{\gamma(\la)}^p)^{-1} \le e^{\frac{\epsilon}{3}} \Vert (DF_{\gamma(\la)}^p(0))^{-1}\Vert$.
\end{enumerate}
To complete  the proof of the Lemma, we have to consider the distortion due to the charts. Replacing $\delta_p(\epsilon)$ by a smaller constant $C_p(\epsilon)$ and recalling that $\tau$ controls this distortion, we obtain for every $\lambda \in M$:
 \begin{enumerate} 
\item[$\cdot$] $(f^p)_{\gamma(\la)}^{-1}$ is defined on $B_{\Pj^k} \left(f^p_\la(\gamma(\la)), C_p(\epsilon)\Vert (DF_{\gamma(\la)}^p(0))^{-1}\Vert^{-2} \right)$,
\item[$\cdot$] $\mbox{Lip} (f^p)_{\gamma(\la)}^{-1} \le e^{\tau + \frac{\epsilon}{3}} \Vert (DF_{\gamma(\la)}^p(0))^{-1}\Vert$. \fin
\end{enumerate}

Let us now
prove Proposition \ref{LemSF}. We recall that $\widehat{u}_p (\widehat{\gamma}) = - \frac{1}{2}  \log  r_p(\gamma_0)$. By assumption  $\lim_n \frac{1}{n}\int_{\widehat{\cal X}} \widehat{u_n}\;d\widehat{\cal M} =L$ with $L\le -\frac{\log  d}{2}$. Let $p\geq 1$ such that $\frac{1}{p}\int_{\widehat{\cal X}} \widehat{u}_p\;d\widehat{\cal M} =:L' \le L +\epsilon$.  By applying Birkhoff Ergodic Theorem there exists $\widehat{\cal Y} \subset\widehat{\cal X}$ such that $\widehat{\cal M} (\widehat{\cal Y})=1$ and
 \begin{eqnarray}\label{Bir}
\forall \widehat{\gamma}\in \widehat{\cal Y} \ , \  \lim_n \frac{1}{n} \sum_{j=1}^n \widehat{u}_p\left(\widehat{\cal F}^{-j} (\widehat{\gamma})\right) = \int_{\widehat{\cal X}} \widehat{u}_p\;d\widehat{\cal M} =pL' .
   \end{eqnarray}
Since $\widehat{u}_p(\widehat{\cal F}^{-n}  (\widehat{\gamma})) = - \frac{1}{2}  \log  r_p(\gamma_{-n})$ we deduce from (\ref{Bir}) that $\lim_n \frac{1}{n} \log  r_p(\gamma_{-n}) =0$. In particular there exists a measurable function $\widehat{r}_p : \widehat{\cal Y} \to ]0,1]$ such that
 \begin{eqnarray*}\label{R}
 C_p(\epsilon)r_p(\gamma_{-n})  \ge \widehat{r}_p (\widehat{\gamma}) e^{-(n-1)\epsilon/2} .
  \end{eqnarray*}
 We also deduce from (\ref{Bir}) that there exists $\widehat{l}_p : \widehat{\cal Y} \to [1,+\infty[$ such that 
  \begin{eqnarray*}\label{RR}
\prod_{j=1}^n \left(r_p(\gamma_{-j})\right)^{-1/2} \leq \widehat{l}_p(\widehat{\gamma}) e^{npL'+n\epsilon/6} .
 \end{eqnarray*}
Now, setting $\widehat{\eta}_p := \widehat{r}_p / \widehat{l}_p$, one can verify by induction:
 \begin{itemize}
\item[$\cdot$]$(f^p)_{\widehat{\gamma}}^{-n}$ is defined on $T_{U_0}(\gamma_0, \widehat{\eta}_p(\widehat{\gamma}))$,
 \item[$\cdot$]$\mbox{Lip} (f^p)_{\widehat{\gamma}}^{-n} \le \widehat{l}_p (\widehat{\gamma}) e^{n(pL'+\tau + \epsilon/2)}$,
 \item[$\cdot$]$(f^p)_{\widehat{\gamma}}^{-n}\left[T_{U_0}(\gamma_0, \widehat{\eta}_p(\widehat{\gamma})) \right] \subset T_{U_0}(\gamma_{-n}, C_p(\epsilon) r_p(\gamma_{-(n+1)}))$. \fin
\end{itemize}

\end{document}